\documentclass[]{article}
\pdfoutput=1

\usepackage{arxiv}

\usepackage{xspace}
\usepackage{graphicx}
\usepackage{float}

\usepackage[dvipsnames]{xcolor}
\usepackage{dsfont}
\usepackage{amsmath}
\usepackage{amssymb}
\usepackage{empheq}
\usepackage{enumitem}
\usepackage{bm}

\usepackage{booktabs}

\usepackage{hyperref}

\usepackage{caption}
\usepackage{multirow}

\usepackage{array}
\usepackage[super]{nth}

\usepackage{textcomp}

\renewcommand{\vec}[1]{\bm{#1}}
\newcommand{\mat}[1]{\bm{#1}}

\newcommand{\block}[1]{\mathbf{#1}}

\newcommand{\statevec}[1]{\vec{#1}}

\newcommand{\face}[1]{\mathfrak{#1}}

\newcommand{\dens}{\rho}
\newcommand{\pres}{p}

\newcommand{\velv}{\vec{v}}

\newcommand{\Ener}{\mathcal{E}}

\newcommand{\naturals}{\mathbb{N}}
\newcommand{\reals}{\mathbb{R}}

\newcommand{\idmat}{\mat{\mathds{1}}}
\newcommand{\bigO}{\mathcal{O}}

\newcommand{\bspace}[0]{\;}
\newenvironment{lcases}
{\left\lbrace\bspace\begin{aligned}}
{\end{aligned}\right.}

\newcommand{\fvdiffmat}{\mat{\Delta}}

\newcommand{\dgMassMat}{\mat{\mathcal{M}}}
\newcommand{\dgDiffMat}{\mat{\mathcal{D}}}
\newcommand{\dgBoundaryMat}{\mat{\mathcal{B}}}

\newcommand{\projMat}{\mat{P}}
\newcommand{\recoMat}{\mat{R}}

\newcommand{\dgVoluOp}{\mat{D}}
\newcommand{\dgSurfOp}{\mat{B}}

\newcommand{\dgfvVoluOp}{\hat{\dgVoluOp}}
\newcommand{\dgfvSurfOp}{\hat{\dgSurfOp}}

\newcommand{\pforest}{\textsc{p4est}\xspace}

\newcolumntype{C}[1]{>{\centering\arraybackslash$}m{#1}<{$}}
\newlength{\mycolwd}                                         
{\noindent \small \textbf{typedef} \emph{#1}\begin{itemize}[label={},topsep=1mm,leftmargin=1\parindent]}
{\end{itemize}\textbf{end typedef}}

{\small \begin{itemize}[label={},topsep=1mm,itemsep=0.5mm,leftmargin=1\parindent]}
{\end{itemize}}

\newcommand{\compilewithplots}{1}

\newenvironment{changed}{\color{black}}{}

\newcommand{\email}[1]{\href{mailto:#1}{#1}}

\begin{document}


\title{A Sub-Element Adaptive Shock Capturing Approach for Discontinuous Galerkin Methods}



\author{Johannes Markert\\
\email{jmarkert@math.uni-koeln.de}, Homepage: \url{www.jmark.de}\\
Department for Mathematics and Computer Science,\\
University of Cologne, Weyertal 86-90, 50931, Cologne, Germany
\And
Gregor Gassner \\
Department for Mathematics and Computer Science;\\
Center for Data and Simulation Science,\\
University of Cologne, Weyertal 86-90, 50931, Cologne, Germany
\And
Stefanie Walch \\
Universität zu Köln, Zülpicher Str. 77, I. Physikalisches Institut;\\
Center for Data and Simulation Science, 50937 Köln, Germany
}


\maketitle

\begin{abstract}
In this paper, a new strategy for a sub-element based shock capturing for
discontinuous Galerkin (DG) approximations is presented. The idea is to
interpret a DG element as a collection of data and construct a hierarchy of low
to high order discretizations on this set of data, including a first order
finite volume scheme up to the full order DG scheme. The different DG
discretizations are then blended according to sub-element troubled cell
indicators, resulting in a final discretization that adaptively
blends from low to high order within a single DG element. The goal is to retain as
much high order accuracy as possible, even in simulations with very strong
shocks, as e.g. presented in the Sedov test. The framework retains the locality
of the standard DG scheme and is hence well suited for a combination with
adaptive mesh refinement and parallel computing. The numerical tests
demonstrate the sub-element adaptive behavior of the new shock capturing approach
and its high accuracy. 
\end{abstract}


\section{Introduction}
\label{sec:intro}

For non-linear hyperbolic problems, such as the compressible Euler equations,
there are two major sources of instabilities when applying discontinuous
Galerkin (DG) methods as a high order spatial discretization. (i) Aliasing is
caused by under-resolution of e.g. turbulent vortical structures and can lead
to instabilities that may even crash the code. As a cure, de-aliasing
mechanisms are introduced in the DG methodology based on e.g. filtering
\cite{Kenevsky:2006qa,hesthaven2008}, polynomial de-aliasing
\cite{Kopriva2017,spiegel2015aliasing,kopriva2019free} or analytical
integration \cite{Kirby2003,Gassner:2013qf}, split forms of the non-linear
terms that for instance preserve kinetic energy
\cite{winters2018compareSplit,Gassner:2016ye,gassner_kepdg,flad2017}, and
entropy stability
\cite{carpenter_esdg,wintermeyer2018entropy,pazner2019analysis,chen2017entropy,Parsani2016,murman2016,chan2018,Parsani:2015:ESD:2784985.2785151,Liu2018,Bohm2018,Gassner2018,gassner_skew_burgers,Friedrich2018}.
(ii) The second major source of instabilities is the Gibbs phenomenon, i.e.
oscillations at discontinuities. It is well known that solutions of non-linear
hyperbolic problems may develop discontinuities in finite time, so called
shocks. While aliasing issues can be mostly attributed to variational crimes,
i.e. a bad design and implementation of the numerical discretization with
simplifications that affects aliasing (e.g. collocation), oscillations at
discontinuities have their root deeper in the innards of the high order DG
approach and are unfortunately an inherent property of the high order
polynomial approximation space. Oscillations are fatal when physical
constraints and bounds on the solution exist that then break because of over-
and undershoots, resulting in nonphysical solutions such as negative density or
pressure. 

In this work, we focus on a remedy for the second major stability issue, often
referred to as shock capturing. There are many shock capturing approaches
available in DG literature, for instance based on slope limiters
\cite{Cockburn1998199,CockburnHou&Shu,Kuzmin2012} or (H)WENO limiters
\cite{Guo2015,Qiu2010,Zhu2009}, filtering \cite{Bohm2019aa}, and artificial
viscosity  \cite{Persson2006,Ching2019}. Slope limiters or (H)WENO limiters are
always applied in combination with troubled cell indicators that detect shocks
and only flag the DG elements that are affected by oscillations. Only in these
elements, the DG polynomial is replaced by an oscillation free reconstruction
using data from neighboring elements. Element based slope limiting, e.g., with
TVB limiters is effective for low order DG discretizations only, as its
accuracy strongly degrades when increasing the polynomial order $N$, as in some
sense the size of the DG elements $\Delta x$ gets larger and larger with higher
polynomial degree $N$. The accuracy of the slope limiters is not based on the
DG resolution $\sim\frac{\Delta x}{N}$, but only on the element size $\Delta
x$. Considering that a finite volume (FV) discretization with high order
reconstructions and limiters typically resolve a shock within about 2-3 cells,
the shock width for high order DG schemes with large elements would be very
wide when relying on element based limited reconstructions. The same behavior
is still an issue for high order reconstruction based limiters such as the
(H)WENO methodology. While these formally use high order reconstructions, its
leading discretization parameter is still $\Delta x$ and not $\frac{\Delta
x}{N}$ and hence the shock width still scales with $\Delta x$. 

Sub-element resolution, i.e. a numerical shock width that scales proportional
to $\frac{\Delta x}{N}$, can be achieved with e.g. artificial viscosity. The
idea is to widen the discontinuity into a sharp, but smooth profile such that
high order polynomials can resolve it. Again, some form of troubled cell
indicator is introduced to not only flag the element that contains the shock
(or that oscillates) but also determine the amount of necessary viscosity,
e.g., based on local entropy production \cite{Zingan2013}. It is interesting
that less viscosity is needed, i.e. sharper shock profiles are possible, the
higher the polynomial degree of the DG discretization. Shocks can be captured
in a single DG element if the polynomial order is high enough
\cite{Persson2006}. An issue of artificial dissipation is that a high amount of
artificial viscosity is needed for very strong shocks, which makes the overall
discretization very stiff with a very small explicit time step restriction.
Local time stepping can be used to reduce this issue \cite{Gassner2015aa} or
specialized many stage Runge-Kutta schemes with optimized coefficients for
strongly dissipative operators \cite{Klckner2011}. 

An alternative idea to achieve sub-element resolution is to actually replace
the elements by sub-elements (lower order DG on finer grid) or maybe even
finite volume (FV)  subcells. A straight forward, maybe even a natural idea in
the DG methodology is to use hp-adaptivity, where the polynomial degree is
reduced and simultaneously the grid resolution is increased close to
discontinuities. When switching to the lowest order DG, i.e a first order DG,
which is nothing but a first order FV method, it is assumed that the inherent
artificial dissipation of the first order FV scheme is enough to clear all
oscillations at even the strongest shocks. The low accuracy of the first order
FV method is compensated by a respective (very strong) increase in local grid
resolution through h-adaptivity. An obvious downside to the hp-adaptivity
approach is the strong dynamic change of the approximation space during the
simulation and the associated strong change of the underlying data structures. 

Thus, as a simpler alternative,  subcell based FV discretizations with a first
order, second order TVD or high order WENO reconstruction was introduced, e.g.,
in \cite{Vilar2019,Sonntag2014,Dumbser2014,Zanotti2015,Dumbser2016}. The idea
is to completely switch the troubled DG element to a robust FV discretization
on a fine  subcell grid. This approach keeps the change of the data local to
the affected DG element and helps to streamline the implementation. The robust
and essential oscillation free FV scheme on the finer  subcell grid in the DG
element is constructed based on the available DG data. This helps to keep the
data structures feasible and reduces the shock capturing again to an element
local technique instead of changing the global mesh topology and approximation
space. Consequently, this again allows to combine the approach with the idea of
a troubled cell mechanism to identify the oscillatory DG elements and replace
the high order DG update with the robust FV update. As noted, there are many
different variants for the FV  subcell scheme with the high order WENO
\cite{Vilar2019,Sonntag2014,Dumbser2014,Zanotti2015,Dumbser2016} variant being
especially interesting, as it alleviates the 'stress' on the troubled cell
indicator. When accidentally switching the whole high order DG element to the
subcell FV scheme in smooth parts of the solution, accuracy may strongly
degrade. A first order FV scheme, even on the fine  subcell grid, would wash
away all fine structure details that the high order DG approach simulated.
Thus, from this point of view, a high order  subcell WENO FV scheme is highly
desired in this context. A downside of a subcell WENO FV scheme is the
non-locality of the data dependence due to the high order reconstruction
stencils. This is especially cumbersome close to the interface of the high
order DG elements, where the reconstruction stencils reach into the
neighboring DG elements and thus fundamentally change the data dependency
footprint of the resulting implementation, that consequently has a direct
impact on the parallelization of the code. 

\begin{figure}[htb!]
\centering
\ifnum \compilewithplots = 1
\includegraphics[width=0.5\textwidth]{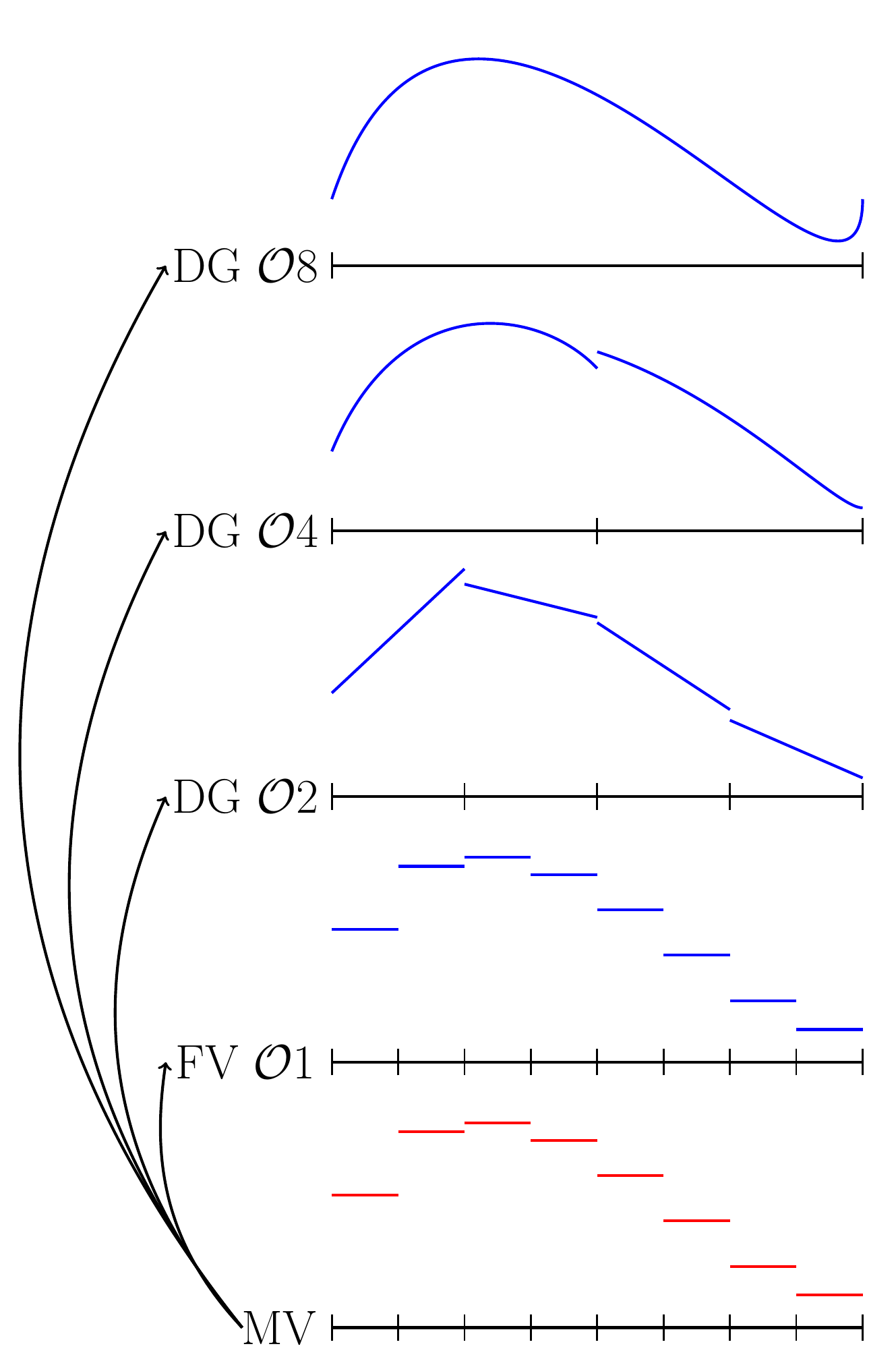}
\fi
\caption{Sketch of the main idea in a one-dimensional setting. Starting with
eight mean values (red), we interpret them as a collection of available data.
We construct four different approximation spaces (blue): piecewise constant on
a subcell grid, piecewise linear and piecewise cubic on a sub-element grid, and
the full 7th degree polynomial on the element.}
\label{fig:mainidea}
\end{figure}

In this paper, the goal is to augment a high order DG method with sub-element
shock capturing capabilities to allow the robust simulation of highly
supersonic turbulent flows which feature strong shocks, as e.g. in
astrophysics, without changing the data dependency footprint. Instead of
flagging an element and completely switching from the high order DG scheme to
the subcell FV scheme, we aim to smoothly combine different schemes in the
flagged element. To achieve this, we firstly interpret the DG element as a
collection of available mean value data. As depicted in Fig.\ref{fig:mainidea},
local reconstructions allow to define a hierarchy of approximation spaces, from
pure piecewise constant approximations (subcell FV) up to a smooth global high
order polynomial (DG element), with all piecewise polynomial combinations
(sub-element DG) in between. We then associate each possible approximation
space with the corresponding FV or DG discretization to define a hierarchy of
schemes that is available for the given collection of data. 

Instead of switching between these different schemes, we aim to smoothly blend
them in order to achieve the highest possible accuracy in every DG element.
Assuming that the first order subcell FV approximation has enough inherent
dissipation to capture shocks in an oscillation free manner, the goal is to
give the FV discretization a high enough weight  at a shock and then smoothly
transition to higher order (sub-element) DG away from the shock. In contrast to
the common approaches, where one indicator is used for the whole DG element, we
aim to introduce sub-element indicators that are adaptively adjusted inside the
DG element. A difficulty that arises when using sub-element indicators and an
adaptive blending approach is to retain conservation of the final
discretization for arbitrary combinations of blending weights. 

\begin{figure}[htb!]
\centering
\ifnum \compilewithplots = 1
\includegraphics[width=0.5\textwidth]{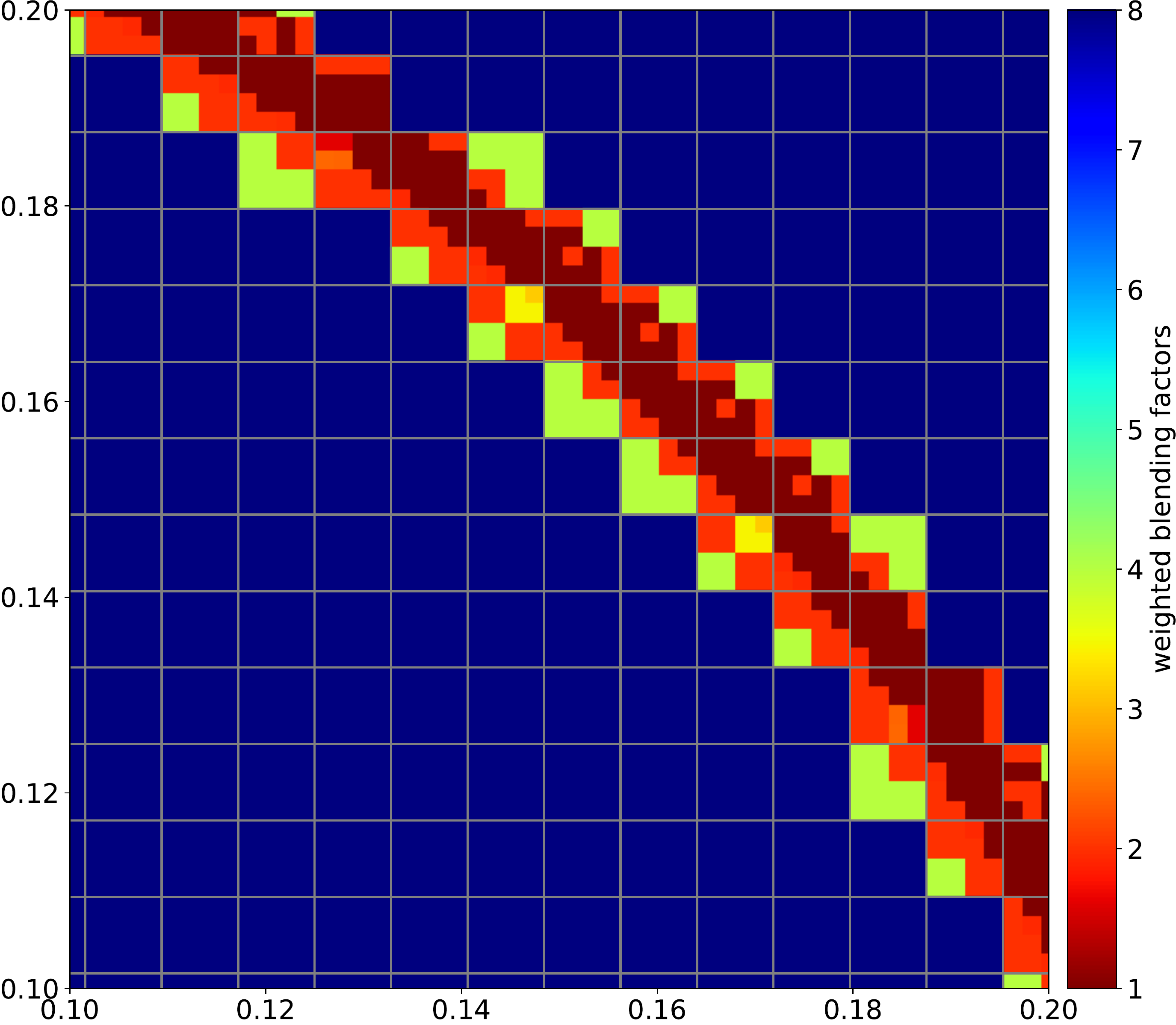}
\fi
\caption{Visualization of the sub-element adaptive blending approach for the 2D
Sedov blast wave, see Sec. \ref{sec:5} for a detailed discussion. The gray
grid lines represent the 8th order DG elements. The color represents the
blending process, where dark red (value $1$) corresponds to pure first order FV
and dark blue (value $8$) corresponds to the pure 8th order DG scheme, with all
combination in between. Note, that the plotted value is only for illustration
and does not correspond to an actual mixed mathematical order of the final
discretization. But the value nicely illustrates how even inside a DG element,
the discretization is adaptively adjusted.}
\label{fig:subcelladaptive}
\end{figure}

The final discretization is a blended scheme, where the weights of the different
discretizations may vary inside the DG element. To demonstrate this idea, the
resulting behavior of the new approach is depicted in Fig.
\ref{fig:subcelladaptive} for the example of a Sedov blast wave in 2D. 
In this particular setup, we start with 8th order DG
as our baseline high order discretization. The grey grid lines indicate the 8th order DG element borders.
Re-interpreting the 2D 8th order DG
element as $8^2$ mean value data, we can construct either a first order FV scheme on $8^2$
subcells, a fourth order DG scheme on $2^2$, or a second order DG scheme on
$4^2$ sub-elements. The blending of the four available discretizations is
visualized via a weighted blending factor (introduced in equation
\eqref{eq:order-estimate}) in Fig. \ref{fig:subcelladaptive}
ranging from $1$ (dark red, basically pure first order FV) up to $8$ (dark blue,
basically pure 8th order DG) and all combinations in between. 

Clearly, the new approach adaptively changes on a sub-element level with a very
localized low order near the shock front, while a substantial part of the high
order DG approximation is recovered away from the shock. We refer to the
validation and application sections, Sec. \ref{sec:5} and Sec. \ref{sec:7}, for
a more detailed discussion.    

The remainder of the paper is organized as follows: in Sec. \ref{sec:2} we
introduce the numerical scheme with the blending of the different
discretizations. Sec. \ref{sec:5} include numerical
validations and Sec. \ref{sec:7} shows an application of the presented
approach with the conclusion given in the last section.       

\section{The Sub-Element Adaptive Blending Approach}
\label{sec:2}

\begin{changed}
The derivations in sections
\ref{sec:finite-volume}-\ref{sec:multilevel-blending} are done in 1D, where the
domain $\Omega$ is divided into $Q$ non-overlapping \textbf{elements}. Each element
$q$ with midpoint $x_q$ and size $\Delta x_q$ is mapped to the reference space as
\begin{equation}
x(\chi) = x_q + \chi\,\Delta x_q\,,\quad \chi \in \left[-\frac{1}{2},\frac{1}{2}\right].
\end{equation}
Each element holds a total number of degrees of freedom (DOF) $N$, where we require $N
\in \{2^r,r\in\mathbb{N}\}$. For the FV method, the element is
divided into $N$ \textbf{subcells} of size $\frac{\Delta x_q}{N}$ leading to a
uniform grid with midpoints $\mu_i$ and interfaces $ \mu_{i\pm\frac{1}{2}}$ in the reference space at
\begin{equation}
\label{eq:subcell-grid}
    \mu_i = -\frac{1}{2} + \frac{i}{N} - \frac{1}{2\,N} \quad \text{and} \quad
    \mu_{i\pm\frac{1}{2}} = \mu_i \pm \frac{1}{2\,N},\quad i = 1,\ldots,N.
\end{equation}
For the sub-element DG scheme of order $n=2^r < N, r\in\mathbb{N}$, the element is
divided into $\frac{N}{n}$ uniform \textbf{sub-elements} of size $\frac{n\Delta
x_q}{N}$, and the mapping of each sub-element $s$ to its (sub-)reference space
reads as
\begin{equation}
\label{eq:sub-element}
    \chi(\xi)_s^{\bigO(n)} = -\frac{1}{2} + \frac{s\,n}{N} - \frac{n}{2\,N} + \frac{n}{N}\xi,
        \quad s = 1,\dots,\frac{N}{n},\; \xi \in \left[-\frac{1}{2},\frac{1}{2}\right].
\end{equation}
Fig. \ref{fig:multilevel-blending} in appendix
\ref{sec:appendix-visu-multilevel-blending} provides a visual representation of
the reference element and its hierarchical decomposition into sub-elements and
subcells.
\end{changed}

\subsection{Finite Volume Method}
\label{sec:finite-volume}

Consider the general 1D conservation law
\begin{equation}
\label{eq:conservation-law}
\partial_t u + \partial_x f(u) = 0.
\end{equation}
We approximate the exact solution $u(x)$ in element $q$ with $N$ mean values
\begin{equation}
    \vec{\bar{u}}_q = (\bar{u}_1,\ldots,\bar{u}_N)_q^T \; \in \reals^N,
\end{equation}
residing on the uniform subcell grid \eqref{eq:subcell-grid}.
Then the 1D FV semi-discretization on subcell $i$ reads
\begin{equation}
\label{eq:fv-simple-form}
\big(\partial_t \bar{u}_{i}\big)_q 
    = -\frac{N}{\Delta x_q} \, \Big(f^*\big((\bar{u_{i}})_q,(\bar{u}_{i+1})_q\big) 
    - f^*\big((\bar{u}_{i-1})_q,(\bar{u}_{i})_q\big)\Big), \quad i \in 1,\ldots,N,
\end{equation}
where $f^*$ denotes a consistent numerical two point flux. The adjacent values
from neighboring elements $q-1$ and $q+1$ are given by $(\bar{u}_0)_q :=
(\bar{u}_N)_{q-1}$ and $(\bar{u}_{N+1})_q := (\bar{u}_1)_{q+1}$. We define the
\emph{Finite Volume Difference Matrix}
\begin{equation}
\label{eq:fv-diff-matrix}
\fvdiffmat := \begin{pmatrix} 
-1      &  1        & 0        & \cdots & 0 \\
0       & -1        & 1        & \ddots & \vdots\\
\vdots  &  \ddots   & \ddots   & \ddots & 0 \\
0       &  \cdots   & 0        & -1     & 1
\end{pmatrix} \in \reals^{N \times (N+1)},
\end{equation}
and rewrite \eqref{eq:fv-simple-form} in compact matrix form
\begin{equation}
\label{eq:fv-matrix-form}
\partial_t \vec{\bar{u}}_q = -\frac{N}{\Delta x_q} \, \fvdiffmat\,\vec{f}_q^{*},
\end{equation}
with the numerical flux vector $\vec{f}_q^* \in \reals^{N+1}$. For later use we split
$\fvdiffmat$ into a \emph{volume} and a \emph{surface} operator:
\begin{equation}
\label{eq:fv-vol-surf-matrices}
\fvdiffmat = \fvdiffmat^{(V)} + \fvdiffmat^{(S)} =
\begin{pmatrix} 
0       &  1        & 0        & \cdots & 0 \\
0       & -1        & 1        & \ddots & \vdots\\
\vdots  &  \ddots   & \ddots   & \ddots & 0 \\
0       &  \cdots   & 0        & -1     & 0
\end{pmatrix} +
\begin{pmatrix} 
-1      &  0        & 0        & \cdots & 0 \\
0       &  0        & 0        & \ddots & \vdots\\
\vdots  &  \ddots   & \ddots   & \ddots & 0 \\
0       &  \cdots   & 0        &  0     & 1
\end{pmatrix}.
\end{equation}
\noindent\textbf{Remark:} 
The FV volume operator vanishes when summed along columns, i.e.
\begin{equation}
\label{eq:fv-vol-mat-property}
\vec{1}^T \, \fvdiffmat^{(V)} = \vec{0}^T,
\end{equation}
where $\vec{1} = (1,\ldots,1)^T \in \reals^N$ is the vector of ones.

\subsection {Discontinuous Galerkin Method}
\label{sec:disc-galerkin}
In this section we briefly outline the Discontinuous Galerkin Spectral Element
Method (DGSEM) \cite{Black1999,kopriva2009implementing}.  To derive the $n$-th
order DG scheme within sub-element $s$, (eq.  \eqref{eq:sub-element}) residing
inside element $q$, we first rewrite the conservation law
\eqref{eq:conservation-law} into variational form with a test function $\phi$
and apply spatial integration-by-parts to the flux divergence
\begin{equation}
\label{eq:weak-form}
\int_{\Omega_s} \partial_t u_q\,\phi \, d\xi =
    \frac{N}{n\,\Delta x_q}\left(\int_{\Omega_s}  f(u_q)\,\partial_{\xi}\phi \, d\xi - f^*(u_q)\,\phi \big|_{\partial \Omega_s} \right).
\end{equation}
We choose the test functions $\phi$ as Lagrange functions
\begin{equation}
\label{eq:lagrange-function}
\ell_i(\xi) = \prod^n_{k=1,k\neq i} \frac{\xi-\xi_k}{\xi_i-\xi_k}, \quad i = 1,\dots,n,
\end{equation}
spanned with the collocation nodes $\xi_i \in
\left[-\tfrac{1}{2},\tfrac{1}{2}\right]$, and use the polynomial ansatz
\begin{equation}
\label{eq:polynomial-ansatz}
    u(t;\xi) \approx \sum_{i=1}^{n} \tilde{u}_{i}(t) \, \ell_i(\xi) = \vec{\tilde{u}}^T(t)\,\vec{\ell}(\xi),
\end{equation}
where
\begin{equation}
    \vec{\tilde{u}}(t) = \big(\tilde{u}_1(t),\ldots,\tilde{u}_n(t)\big)^T \quad \text{and}\quad \vec{\ell}(\xi) = \big(\ell_1(\xi),\ldots,\ell_n(\xi)\big)^T
\end{equation}
are the vectors of nodal coefficients and Lagrange polynomials. We use
collocation of the non-linear flux functions with the same polynomial
approximation. Numerical quadrature is collocated as well. The nodes
are the Gauss-Legendre quadrature nodes. These quadrature nodes and the interpolation polynomials are also illustrated in Fig. \ref{fig:multilevel-blending} in the appendix
\ref{sec:appendix-visu-multilevel-blending}.

Inserting everything into \eqref{eq:weak-form} gives the semi-discrete
weak-form DG scheme
\begin{equation}
\label{eq:semi-discrete-form-dg}
\dgMassMat \, \partial_t \vec{\tilde{u}}_s
    = \frac{N}{n\,\Delta x_q}\left(\dgDiffMat^T\dgMassMat\,\vec{\tilde{f}}_s - \dgBoundaryMat\,\vec{f}_s^* \right).
\end{equation}
The diagonal mass matrix is
\begin{equation}
\left(\dgMassMat\right)_{ij} = \delta_{ij}\,\omega_i,
\end{equation}
with associated numerical quadrature weights $\omega_i$, and the
differentiation matrix is
\begin{equation}
\left(\dgDiffMat\right)_{ij} = \partial_{\xi} \ell_j(\xi)\big|_{\xi_i}.
\end{equation}
The DG volume flux $\vec{\tilde{f}}_s$ is computed directly with the nodal values
$\vec{\tilde{u}}_s$ of the polynomial ansatz
\eqref{eq:polynomial-ansatz}
\begin{equation}
\label{eq:dg-volume-flux}
    \vec{\tilde{f}}_s = f\left(\vec{\tilde{u}}_s\right) := \Big(f\big((\tilde{u}_1)_s\big),\ldots,f\big((\tilde{u}_n)_s\big)\Big)^T\; \in \reals^{n},
\end{equation}
and the surface flux vector
\begin{equation}
\label{eq:dg-surface-flux}
\vec{f}_s^* = \Big(f^*(\underbrace{\tilde{u}_{s-1}^+,\tilde{u}_s^-}_{\mathclap{\substack{\text{left sub-element}\\\text{interface}}}}),0,\ldots,0, f^*(\underbrace{\tilde{u}_{s}^+,\tilde{u}_{s+1}^-}_{\mathclap{\substack{\text{right sub-element}\\\text{interface}}}})\Big)^T \; \in \reals^n
\end{equation}
is constructed with the numerical flux evaluated at DG sub-element boundaries
\begin{equation}
\label{eq:dg-boundary-values}
\tilde{u}_s^+ = (\vec{b}^+)^T\,\vec{\tilde{u}}_s \quad \text{and} \quad \tilde{u}_s^- = (\vec{b}^-)^T\,\vec{\tilde{u}}_s,
\end{equation}
where $\vec{b}^\pm$ are the boundary interpolation operators
\begin{equation}
\label{eq:boundary-interpolation-operators}
\vec{b}^\pm = \Big(\ell_1\big(\pm \tfrac{1}{2}\big),\ldots,\ell_n\big(\pm \tfrac{1}{2}\big)\Big)^T \; \in \reals^{n}.
\end{equation}
The boundary evaluation operator can be compactly written as 
\begin{equation}
\label{eq:boundary-evalution-matrix}
\dgBoundaryMat = \begin{pmatrix} 
-\ell_1(-\tfrac{1}{2}) & 0 & \cdots & 0 & \ell_1(\tfrac{1}{2}) \\
 \vdots  & \vdots  &   \ddots     &  \vdots & \vdots  \\
-\ell_n(-\tfrac{1}{2}) & 0 & \cdots & 0 & \ell_n(\tfrac{1}{2}) \\
\end{pmatrix} \; \in \reals^{n \times n}.
\end{equation}
We rewrite \eqref{eq:semi-discrete-form-dg} and arrive at
\begin{equation}
\label{eq:semi-discrete-dg}
\partial_t \vec{\tilde{u}}_s = \frac{N}{n\,\Delta x_q}\left(\dgVoluOp \,\vec{\tilde{f}}_s - \dgSurfOp \, \vec{f}_s^{*}\right)
\end{equation}
with $\mat{D} := \dgMassMat^{-1}\,\dgDiffMat^T\,\dgMassMat$ and $\dgSurfOp := \dgMassMat^{-1}\dgBoundaryMat$.\\

\noindent\textbf{Remark:} 
The DG volume operator $\dgVoluOp$ vanishes when contracted with the
vector of quadrature weights $\vec{\omega}^T := (\omega_1,\ldots,\omega_n)^T =
\vec{1}^T \, \dgMassMat$.
\begin{equation}
\label{eq:dg-vol-mat-property}
    \vec{\omega}^T \, \dgVoluOp = \vec{0}^T
\end{equation}

\noindent\emph{Proof:} 
\begin{equation*}
\vec{\omega}^T \, \dgVoluOp 
    = \vec{1}^T\,\dgMassMat \, \dgMassMat^{-1}\,\dgDiffMat^T\,\dgMassMat 
    = \vec{1}^T\,\dgDiffMat^T\,\dgMassMat = (\dgDiffMat\,\vec{1})^T\,\dgMassMat  = \vec{0}\\[2mm]
\end{equation*} \hfill$\square$ \\
Here, we use the fact that $\dgDiffMat\,\vec{1} = \vec{0}$ is equivalent to taking the
derivative of the constant function $u(x) = 1$, i.e. $\partial_x u(x) = \partial_x 1 =
0$.

\subsection{Projection Operator}
\label{sec:projection}
\begin{changed}
We consider a given polynomial of degree $n-1$ with its nodal data
$\vec{\tilde{u}} \in \reals^{n}$ in a (sub-)element, which we want to project
to $n_{\mu}$ mean values on regular subcells.  We define the projection
operator $\projMat^{(n \rightarrow n_{\mu})} \in \reals^{n_{\mu} \times n}$
component-wise via the mean value of the polynomial in the interval of
subcell $\mu_i$,
\begin{equation}
\bar{u}_i = n_{\mu}\int_{\mu_{i-\frac{1}{2}}}^{\mu_{i+\frac{1}{2}}} \tilde{u}(\xi) \, d\xi
    = \sum_{j=1}^n \tilde{u}_j  \underbrace{n_{\mu}\int_{\mu_{i-\frac{1}{2}}}^{\mu_{i+\frac{1}{2}}} \ell_j(\xi)\,d\xi}_{:= P_{ij}^{(n\rightarrow n_{\mu})}}, \quad i = 1,\dots,n_{\mu}.
\end{equation}
The integration of the Lagrange polynomial is done exactly with an appropriate
quadrature rule.\\

\noindent\textbf{Remark:} 
Contracting the projection operator $\projMat^{(n \rightarrow n_{\mu})}$ with
$\vec{1}_{n_{\mu}} := (1,\ldots,1)^T \in \reals^{n_{\mu}}$ gives:
\begin{equation}
\label{eq:projection-mat-property}
    \vec{1}_{n_{\mu}}^T\,\projMat^{(n \rightarrow n_{\mu})} = n_{\mu}\,\vec{\omega}^T.
\end{equation}\\
\noindent\emph{Proof:}
\begin{equation}
   \sum_{i=1}^{n_{\mu}} P_{ij} = n_{\mu}\sum_{i=1}^{n_{\mu}} \int_{\mu_{i-\frac{1}{2}}}^{\mu_{i+\frac{1}{2}}} \ell_j(\xi)\,d\xi = n_{\mu}\int_{\mu_{1-\frac{1}{2}}}^{\mu_{n_{\mu}+\frac{1}{2}}} \ell_j(\xi)\,d\xi = n_{\mu}\int_{-\frac{1}{2}}^{\frac{1}{2}} \ell_j(\xi)\,d\xi = n_{\mu} \omega_j
\end{equation} \hfill$\square$

\noindent\textbf{Remark:} When $n = n_{\mu}$ we succinctly write $\projMat
:= \projMat^{(n \rightarrow n)}$.
\end{changed}

\subsection{Reconstruction Operator}
\label{sec:reconstruction}
The inverse of the quadratic non-singular projection matrix $\projMat^{(n
\rightarrow n)}$ can be interpreted as a reconstruction operator from $n$ mean
values to $n$ nodal values of a polynomial with degree $n-1$:
\begin{equation}
\label{eq:trafo-mv-nodal}
\vec{\tilde{u}} = \projMat^{-1} \,\vec{\bar{u}} := \recoMat \, \vec{\bar{u}}.
\end{equation}
However, the na\"ive reconstruction \eqref{eq:trafo-mv-nodal} can lead to
invalid data such as negative density or energy. Considering a convex set of
permissible states for a given conservation law \cite{zhang2012positivity} we
want to have a limiting procedure which recovers permissible states while still
conserving the element mean value. Provided that the element average is a valid
state,
\begin{equation}
\label{eq:element-average}
    \langle\vec{u}\rangle := \frac{1}{n} \, \vec{1}^T \vec{\bar{u}} \;\; \in \big\{\text{permissible states}\big\},
\end{equation}
we expand the reconstruction process \eqref{eq:trafo-mv-nodal} by the following limiting procedure
\begin{equation}
\label{eq:positivity-limiter}
\recoMat^{(\beta)} = \frac{1-\beta}{n} \, \vec{1} \otimes \vec{1}^T + \beta \, \recoMat,
\end{equation}
with the dyadic product $\vec{1} \otimes \vec{1}^T \in \reals^{n \times n}$ and
the ``squeezing'' parameter $\beta \in [0,1]$ chosen such that
\begin{equation}
\label{eq:positivity-limiter-2}
\big\{
\underbrace{\vphantom{\big(\vec{b}^+\big)^T\,\vec{\tilde{u}}^{(\beta)}}\vec{\tilde{u}}^{(\beta)}}_{\text{node values}},\;
\underbrace{\big(\vec{b}^+\big)^T\,\vec{\tilde{u}}^{(\beta)}}_{\text{left boundary}}, \;
\underbrace{\big(\vec{b}^-\big)^T\,\vec{\tilde{u}}^{(\beta)}}_{\text{right boundary}}\big\} \subset \big\{\text{permissible states}\big\},
\end{equation}
where
\begin{equation}
\label{eq:positivity-limiter-applied}
\vec{\tilde{u}}^{(\beta)} = \recoMat^{(\beta)}\,\vec{\bar{u}}
    = \left(\frac{1-\beta}{n}\, \vec{1} \otimes \vec{1}^T + \beta \, \recoMat\right) \vec{\bar{u}}
    = (1-\beta)\,\vec{1}\,\langle\vec{u}\rangle + \beta \, \vec{\tilde{u}},
\end{equation}
which shows that the extended reconstruction operator preserves the original element average, i.e.
\begin{equation*}
\langle\vec{\tilde{u}}^{(\beta)}\rangle = \langle\vec{\tilde{u}}\rangle.
\end{equation*}

\noindent\textbf{Remark:} 
The expanded reconstruction operator fulfills the following relationships:
\begin{equation}
\label{eq:reconstruction-mat-property}
    \vec{\omega}^T\,\recoMat^{(\beta)} = \vec{\omega}^T\,\recoMat = \frac{1}{n}\,\vec{1}^T.
\end{equation}
\noindent\emph{Proof:} We begin with the second relation. From
\eqref{eq:projection-mat-property} we know:
\begin{align*}
    \vec{1}^T\,\projMat &= n\,\vec{\omega}^T \\
    \vec{1}^T\,\projMat \, \projMat^{-1} &= n\,\vec{\omega}^T\projMat^{-1} = n\,\vec{\omega}^T\recoMat\\
    \frac{1}{n}\,\vec{1}^T &= \vec{\omega}^T\recoMat.
\end{align*}
The first relation follows as
\begin{equation*}
    \vec{\omega}^T\,\recoMat^{(\beta)}
        = \frac{1-\beta}{n}\, \underbrace{\vec{\omega}^T\, \vec{1}}_{= 1} \otimes \vec{1}^T + \beta \, \underbrace{\vec{\omega}^T\, \recoMat}_{= \frac{1}{n}\,\vec{1}^T}
        = \frac{1}{n}\,\vec{1}^T  = \vec{\omega}^T\,\recoMat.
\end{equation*} \hfill$\square$\\

\begin{changed}
The specific calculation of the parameter $\beta$ depends on the permissibility
constraints of the conservation law at hand. In section \ref{sec:gov-equations}
we solve the compressible Euler equations where the density $\rho$ and
pressure $p$ are positive values. Based on \cite{zhang2012positivity} we
calculate the squeezing parameter $\beta$ as
\begin{equation}
\label{eq:calc-squeezing-parameter}
    \beta = \min\left(\frac{\langle\rho\rangle-\epsilon}{\langle\rho\rangle - \rho_{\text{min}}},\frac{\langle p \rangle-\epsilon}{\langle p \rangle - p_{\text{min}}},1\right),
\end{equation}
where $\langle\rho\rangle$, $\langle p \rangle$ are the element averages
\eqref{eq:element-average}, $\rho_{\text{min}}$, $p_{\text{min}}$ are the
minimum values of the unlimited polynomial given in
\eqref{eq:positivity-limiter-2} and $\epsilon =
\min(10^{-20},\langle\rho\rangle,\langle p \rangle)$.

\noindent\textbf{Remark:}
We observed that the indispensable positivity limiting can have a potential
negative impact on the robustness of the high order DG scheme, especially for
high polynomial degrees. It seems that huge artificial jumps and gradients can
be generated at element interfaces.  To alleviate this problem it is advisable
to only allow small squeezing margins between 5\% to 10\%. If for example the
squeezing parameter $\beta$ falls below 0.95 the reconstructed high order
polynomial is decided un-salvageable and the next lower order scheme is considered.
\end{changed}

\subsection{Single-level Blending}
\label{sec:convex-blending}
We first focus on the case where only two different discretizations per
element $q$ are blended: a single-level blending of the low order FV
scheme and the high order DG scheme
\begin{equation*}
\partial_t u_{q} = (1-\alpha_{q})\,\partial_t u_{q}^{(\text{low})} + \alpha_{q}\,\partial_t u_{q}^{(\text{high})}.
\end{equation*}
Note that we do not blend the solutions, but directly the discretizations
themselves, i.e. the right hand sides which we denote by $\partial_t u_{q}$. If
both schemes operate with different data representations, appropriate
transformations ensure compatibility during the blending process. In our case
we use the projection operator $\projMat^{(N \rightarrow N)}$ introduced in
section \ref{sec:projection} to transfer the nodal output of the DG operator to
subcell mean values. The input for the DG scheme on the other hand, i.e. the
nodal coefficients $\vec{\tilde{u}}_q \in \reals^N$, are reconstructed from the
given mean values $\vec{\bar{u}}_q \in \reals^N$ as described in section
\ref{sec:reconstruction}. The na\"ive single-level blended discretization reads
as
\begin{equation*}
\partial_t \vec{\bar{u}}_{q} = (1-\alpha_{q})\;\partial_t \vec{\bar{u}}_{q}^{(\text{FV})} + \alpha_{q}\;\projMat\,\partial_t \vec{\tilde{u}}_{q}^{(\text{DG})}.
\end{equation*}
Inserting \eqref{eq:fv-matrix-form} and \eqref{eq:semi-discrete-dg} gives
\begin{equation}
\label{eq:dgfv-scheme-1d-nonconservative}
\partial_t \vec{\bar{u}}_{q}
= (1-\alpha_{q}) \, \frac{-N}{\Delta x_q} \, \fvdiffmat\,\vec{f}_{q}^{*(\text{FV})}
    + \alpha_{q} \,\frac{1}{\Delta x_q}\projMat\,\left(\dgVoluOp\,\vec{\tilde{f}}_{q} - \dgSurfOp\,\vec{f}_{q}^{*(\text{DG})} \right).
\end{equation}
\begin{figure}[htb!]
\centering
\ifnum \compilewithplots = 1
    \includegraphics[width=0.8\textwidth]{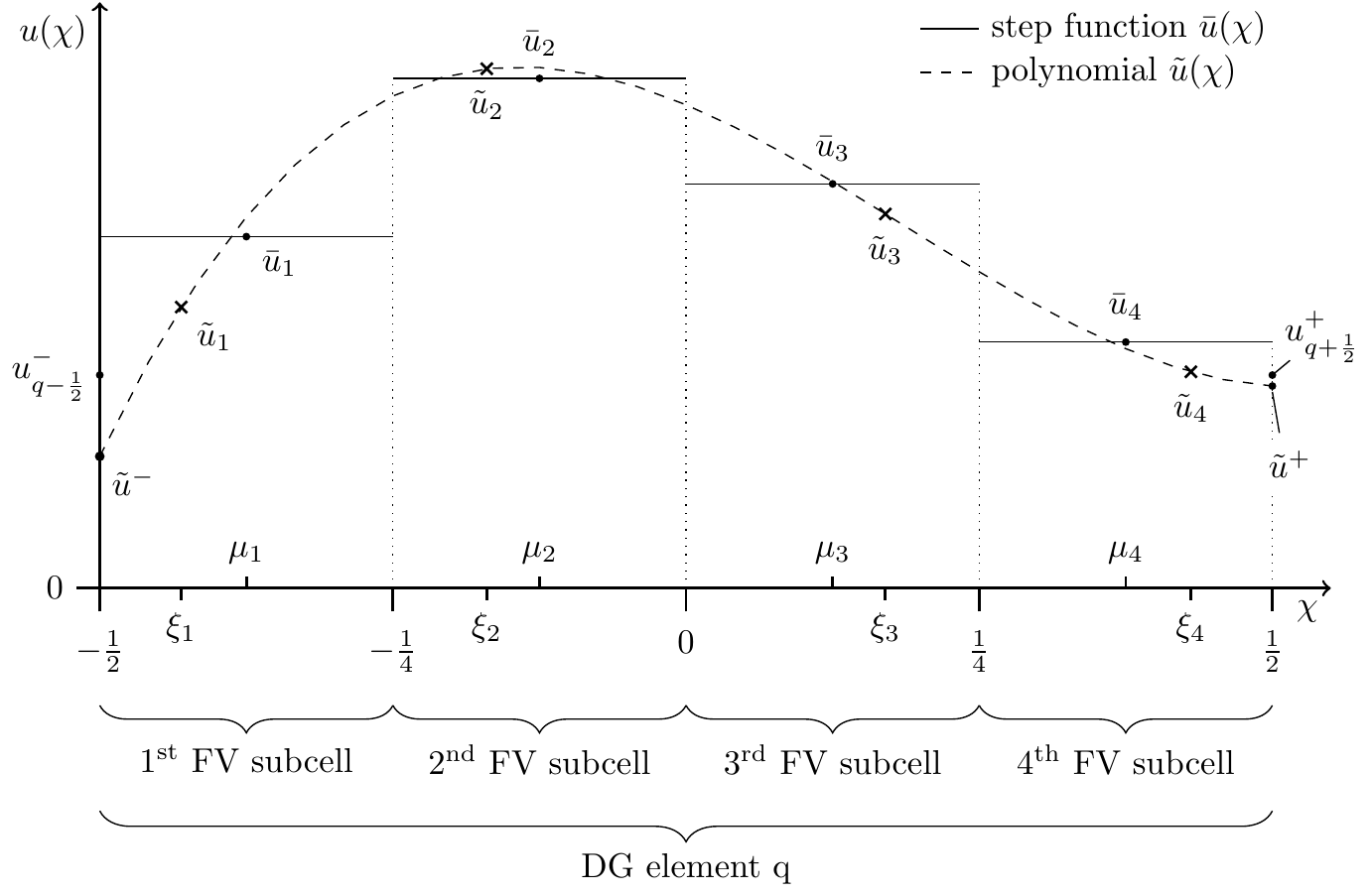}
\fi
\caption{Schematic of four ($N = 4$) mean values and the reconstructed
polynomial on the reference element. The cell boundary values $u_{q-\frac{1}{2}}^-$ and
$u_{q+\frac{1}{2}}^+$ were calculated with \eqref{eq:inner-surface-prolongation} where we
assumed a blending factor of $\alpha= 0.75$.}
\label{fig:dgfv-single-blend}
\end{figure}
The fundamental issue of the na\"ive single-level blending discretization
\eqref{eq:dgfv-scheme-1d-nonconservative} is that it is not conservative if the blending
parameter $\alpha_{q}$ varies from element to element. It turns out
that the direct blending of the surface fluxes of both discretizations, namely
$\vec{f}_{q}^{*(\text{FV})}$ and $\vec{f}_{q}^{*(\text{DG})}$, leads to a
non-conservative balance of the net fluxes across element interfaces. To remedy
this issue the goal is to define unique surface fluxes common to both, the low
order and the high order discretization, such that they can be directly
blended. We thus compute a blending of the reconstructed interface values to
evaluate a unique common surface flux  $\vec{f}_{q}^{*}$. First, we define two
interface vectors of length $N+1$
\begin{equation}
\label{eq:assemble-boundary-values}
\vec{u}_{q}^+ := \big(\underbrace{u_{q-\frac{1}{2}}^+, (\bar{u}_{1})_q, \ldots, (\bar{u}_{N-1})_q, u_{q+\frac{1}{2}}^+}_{\substack{\text{prolongated states at} \\ \text{FV subcell interfaces}}} \big)^T
\quad \text{and} \quad
\vec{u}_{q}^- := \big(\underbrace{u_{q-\frac{1}{2}}^-, (\bar{u}_{2})_q, \ldots, (\bar{u}_{N})_q, u_{q+\frac{1}{2}}^-}_{\substack{\text{prolongated states at} \\ \text{FV subcell interfaces}}} \big)^T,
\end{equation}
where the outermost interface values are given by
\begin{equation}
\label{eq:inner-surface-prolongation}
u_{q-\frac{1}{2}}^- = (1-\alpha_{q})\,(\bar{u}_{1})_q + \alpha_{q} \,\tilde{u}_{q}^-
\quad \text{and} \quad
u_{q+\frac{1}{2}}^+ = (1-\alpha_{q})\,(\bar{u}_{N})_q + \alpha_{q} \,\tilde{u}_{q}^+.
\end{equation}
A concrete example for the blended boundary values $u_{q\pm\frac{1}{2}}^{\pm}$ is shown in
Fig. \ref{fig:dgfv-single-blend}.

The common surface flux $\vec{f}_q^{*}$ is then evaluated with the interface
vectors $\vec{u}_{q}^+$ and $\vec{u}_{q}^-$ as
\begin{equation}
\label{eq:common-surface-flux}
\vec{f}_{q}^* = f^*\left(\vec{u}_{q}^+,\vec{u}_{q}^-\right) \; \in \reals^{N+1}.
\end{equation}
To make the DG surface operator $\dgSurfOp$ compatible with $\vec{f}_{q}^*$, we
adapt the notation by inserting an additional column of zeros
\begin{equation}
\label{eq:tweak-boundary-evaluation-matrix}
\dgSurfOp^{(0)} = \begin{pmatrix} 
B_{11} & 0 & \cdots & 0 & \cdots & 0 & B_{1N} \\
B_{21} & 0 & \cdots & 0 & \cdots & 0 & B_{2N} \\
             \vdots &     \vdots & \ddots & \vdots  & \ddots & \vdots& \vdots  \\
B_{N1} & 0 & \cdots & 0 & \cdots & 0 & B_{NN} \\
\end{pmatrix}  \; \in \reals^{N \times (N+1)}.
\end{equation}
\textbf{Corollary:} 
Contracting the DG surface operator $\dgSurfOp^{(0)}$ with the quadrature
weights $\vec{\omega}$ is equivalent to the contraction of the FV surface
operator, i.e.
\begin{equation}
\label{eq:dg-surf-mat-property}
    \vec{\omega}^T \, \dgSurfOp^{(0)} = \vec{1}^T\,\fvdiffmat^{(S)}.
\end{equation}

\noindent\emph{Proof:} We expand the left side of
\eqref{eq:dg-surf-mat-property} and write
\begin{align*}
\vec{\omega}^T \, \dgSurfOp^{(0)}
    &= \left( -\sum_{i=1}^{N} \omega_i \frac{\ell_i(-\tfrac{1}{2})}{\omega_i},0,\ldots,0,\sum_{i=1}^{N} \omega_i \frac{\ell_i(\tfrac{1}{2})}{\omega_i}\right)\\
    &= \left( -\sum_{i=1}^{N} \ell_i(-\tfrac{1}{2}),0,\ldots,0,\sum_{i=1}^{N} \ell_i(\tfrac{1}{2})\right) = \left( -1,0,\ldots,0,1\right) = \vec{1}^T\,\fvdiffmat^{(S)}.
\end{align*}\hfill$\square$

We replace $\vec{f}_{q}^{*(\text{FV})}$ and $\vec{f}_{q}^{*(\text{DG})}$ in
\eqref{eq:dgfv-scheme-1d-nonconservative} with $\vec{f}_{q}^{*}$, and additionally,
write the DG operators compactly, including the projection operator, as 
\begin{equation}
\label{eq:dgfv-operators}
\dgfvVoluOp = \projMat\,\dgVoluOp \quad \text{and} \quad \dgfvSurfOp = \projMat\,\dgSurfOp^{(0)}.
\end{equation}
The final form of the single-level blended discretization is
\begin{equation}
\label{eq:dgfv-scheme-1d}
\partial_t \vec{\bar{u}}_{q}
    = (1-\alpha_{q}) \, \frac{-N}{\Delta x} \, \fvdiffmat\,\vec{f}_{q}^*
    + \alpha_{q} \,\frac{1}{\Delta x}\,\left(\dgfvVoluOp\,\vec{\tilde{f}}_{q} - \dgfvSurfOp\,\vec{f}_{q}^* \right).
\end{equation}

\textbf{Remark:} It is easy to see that with $\alpha_{q} = 0$ for all elements
$q$, the pure subcell FV discretization is recovered and with $\alpha_{q} = 1$
for all elements the blending scheme recovers the pure high order DG method.\\ 

\textbf{Lemma:} Given arbitrary blending factors $\alpha_{q} \in \reals$ for
each element $q$ the blending scheme \eqref{eq:dgfv-scheme-1d} is conservative.
\\

\emph{Proof:} We discretely integrate the single-level blended discretization over all elements
$q$ with a total number of $Q$
\begin{equation*}
\label{eq:mass-conservation}
\sum_{q}^Q \Delta x_q \, \vec{1}^T \partial_t \vec{\bar{u}}_q
    = \sum_{q}^Q \Delta x_q \, \vec{1}^T \Big\{(1-\alpha_q) \, \frac{-N}{\Delta x_q} \, \fvdiffmat\,\vec{f}^*_q
        + \alpha_q \,\frac{1}{\Delta x_q}\left(\dgfvVoluOp\,\vec{\tilde{f}}_q - \dgfvSurfOp\,\vec{f}^*_q \right)\Big\}.
\end{equation*}
With \eqref{eq:fv-vol-surf-matrices} and \eqref{eq:dgfv-operators} we get
\begin{equation*}
\sum_{q}^Q \Delta x_q \, \vec{1}^T \partial_t \vec{\bar{u}}_q
    = \sum_{q}^Q \vec{1}^T \Big\{-(1-\alpha_q) \, N \, \big(\fvdiffmat^{(V)} + \fvdiffmat^{(S)}\big)\,\vec{f}^*_q
    + \alpha_q \,\left(\projMat \, \dgVoluOp\,\vec{\tilde{f}}_q - \projMat \, \dgSurfOp^{(0)}\,\vec{f}^*_q \right)\Big\}.
\end{equation*}
With properties \eqref{eq:fv-vol-mat-property}, \eqref{eq:dg-vol-mat-property}
and \eqref{eq:projection-mat-property} the volume terms vanish, i.e.
\begin{equation*}
\vec{1}^T \, \fvdiffmat^{(V)} = \vec{0}^T \quad \text{and} \quad \vec{1}^T \, \projMat \, \dgVoluOp = N\,\vec{\omega}^T \, \dgVoluOp = \vec{0}^T.
\end{equation*}
We reformulate the DG surface term with \eqref{eq:projection-mat-property} and \eqref{eq:dg-surf-mat-property} as
\begin{equation*}
\vec{1}^T \, \projMat \, \dgSurfOp^{(0)} = N \, \vec{\omega}^T \, \dgSurfOp^{(0)} = N \, \vec{1}^T \fvdiffmat^{(S)}.
\end{equation*}
The expansion of the telescopic sum only leaves the outermost surface fluxes
\begin{equation*}
    \sum_{q}^Q \Big[-(1-\alpha_q) \, N \, \vec{1}^T \fvdiffmat^{(S)} - \alpha_q \,N \, \vec{1}^T \fvdiffmat^{(S)}\Big] \, \vec{f}^*_q
    = \sum_{q}^Q -N \, \vec{1}^T \fvdiffmat^{(S)} \, \vec{f}^*_q\\
\end{equation*}
\begin{equation*}
    = -N\,\Big(\big(-f^*_{1}\big)_{1} + \big(f^*_{N+1}\big)_{Q}\Big) = 0,
\end{equation*}
which represent the change due to physical boundary conditions. For instance in case of periodic boundary conditions, these two fluxes would cancel to zero. \hfill$\quad\square$\\

This concludes the description of the single-level blending scheme.

\subsection{Multi-level Blending}
\label{sec:multilevel-blending}
As mentioned in the introduction, the general idea is the
construction of a \emph{hierarchy} of discretizations, with lower order DG schemes on
sub-elements. In principle, it is thus possible to define a multi-level
blending discretization, where the low order FV scheme is blended with DG
variants of different approximation order. This multi-level extension is driven
by the desire to retain as much 'high order DG accuracy' as possible,
especially on a sub-element level. 

For the discussion, we consider a specific example setup with a DG element with
polynomial degree $N_p=7$. The data of this DG element is collected in form of
$N=8$ mean values on a regular subcell grid. Our goal is to blend a first order
subcell FV scheme, a second and fourth order sub-element DG scheme and the full
eight order DG scheme. The construction of the individual schemes follows the
discussion in the previous sections \ref{sec:finite-volume},
\ref{sec:disc-galerkin}, and \ref{sec:convex-blending}. In this section, the
extension to a multi-level blending approach is presented. 

We get the first data interpretation and its corresponding discretization by
directly using the mean values of element $q$ as a first order approximation
space with no reconstruction at all
\begin{equation*}
    \vec{\bar{u}}_q^{\bigO(1)} := \vec{\bar{u}}_q = (\bar{u}_1,\bar{u}_2,\ldots,\bar{u}_8)_q^T.
\end{equation*}
Next, we interpret the eight mean values as four second order DG sub-elements.
The reconstruction matrix $\mat{R}^{\bigO(2)} \in \reals^{2 \times 2}$
transforms two adjacent mean values to two nodal values
\begin{equation}
\label{eq:reco-2nd-order}
\vec{\tilde{u}}_q^{\bigO(2)} = 
\Big(\underbrace{\left(\tilde{u}_1\right)_1^{\bigO(2)},\left(\tilde{u}_2\right)_1^{\bigO(2)}}_{\bigO(2)-\text{sub-element}},\ldots,\underbrace{\left(\tilde{u}_1\right)_4^{\bigO(2)},\left(\tilde{u}_2\right)_4^{\bigO(2)}}_{\bigO(2)-\text{sub-element}}\Big)_q^{T}
= \left[\idmat_4 \otimes \mat{R}^{\bigO(2)} \right] \, \vec{\bar{u}}_q.
\end{equation}
The fourth and eight order approximations are constructed analogously as
\begin{equation}
\label{eq:reco-4th-order}
\vec{\tilde{u}}_q^{\bigO(4)} = 
\Big(\underbrace{\left(\tilde{u}_1\right)_1^{\bigO(4)},\ldots,\left(\tilde{u}_4\right)_1^{\bigO(4)}}_{\bigO(4)-\text{sub-element}},
\underbrace{\left(\tilde{u}_1\right)_2^{\bigO(4)},\ldots,\left(\tilde{u}_4\right)_2^{\bigO(4)}}_{\bigO(4)-\text{sub-element}}\Big)_q^{T}
= \left[\idmat_2 \otimes \mat{R}^{\bigO(4)} \right] \, \vec{\bar{u}}_q,
\end{equation}
and
\begin{equation}
\label{eq:reco-8th-order}
\vec{\tilde{u}}_q^{\bigO(8)} = 
\Big(\underbrace{\tilde{u}_1^{\bigO(8)},\ldots,\tilde{u}_8^{\bigO(8)}}_{\bigO(8)-\text{sub-element}}\Big)_q^{T}
= \left[\idmat_1 \otimes \mat{R}^{\bigO(8)}\right]\,\vec{\bar{u}}_q,
\end{equation}
with $\mat{R}^{\bigO(4)} \in \reals^{4\times4}$ and $\mat{R}^{\bigO(8)} \in
\reals^{8\times8}$. Here we used the Kronecker product $\otimes$ in conjunction
with the identity matrix $\idmat_n = \text{diag}(1,\ldots,1) \in \reals^{n
\times n}$ in order to generate appropriate block-diagonal matrices. Fig.
\ref{fig:multilevel-blending} in appendix
\ref{sec:appendix-visu-multilevel-blending} illustrates the hierarchy of the
approximation spaces for this example. 

We correspond the approximation spaces with the respective discretizations 
\begin{align}
\label{eq:multilevel-blending-scheme}
\partial_t \vec{\bar{u}}_q^{\bigO(1)} &= -\frac{8}{\Delta x_q} \; \fvdiffmat\,\vec{f}^{*} \nonumber \\[2mm]
\partial_t \vec{\bar{u}}_q^{\bigO(2)} &= \,\;\; \frac{4}{\Delta x_q} \left\{\left[\idmat_4 \otimes \dgfvVoluOp^{\bigO(2)}\right]\vec{\tilde{f}}_q^{\bigO(2)} - \left[\idmat_4 \otimes \dgfvSurfOp^{\bigO(2)}\right]^{(0)}\vec{f}_q^* \right\}, \nonumber \\[2mm]
\partial_t \vec{\bar{u}}_q^{\bigO(4)} &= \,\;\; \frac{2}{\Delta x_q} \left\{\left[\idmat_2 \otimes \dgfvVoluOp^{\bigO(4)}\right]\vec{\tilde{f}}_q^{\bigO(4)} - \left[\idmat_2 \otimes \dgfvSurfOp^{\bigO(4)}\right]^{(0)}\vec{f}_q^* \right\}, \\[2mm]
\partial_t \vec{\bar{u}}_q^{\bigO(8)} &= \,\;\; \frac{1}{\Delta x_q} \left\{\left[\idmat_1 \otimes \dgfvVoluOp^{\bigO(8)}\right]\vec{\tilde{f}}_q^{\bigO(8)} - \left[\idmat_1 \otimes \dgfvSurfOp^{\bigO(8)}\right]^{(0)}\vec{f}_q^* \right\}, \nonumber 
\end{align}
where the DG volume fluxes are computed with the respective polynomial
\mbox{(sub-)element} reconstructions of orders $n=\{2,4,8\}$:
\begin{equation*}
    \vec{\tilde{f}}_q^{\bigO(n)} = f\left(\vec{\tilde{u}}_q^{\bigO(n)}\right).
\end{equation*}
Similar to \eqref{eq:tweak-boundary-evaluation-matrix} the DG surface operators
$\big[\idmat_{N/n} \otimes \dgfvSurfOp^{\bigO(n)}\big]$ had to be slightly
adapted in order to be compatible with the common surface flux $\vec{f}^{*}$.
We list the modified boundary evaluation matrices in the appendix
\ref{sec:appendix-boundary-eval-op}. Finally, all four candidate
discretizations are blended starting from the lowest order up to the highest
order discretization
\begin{align}
\label{eq:multilevel-blending-discretizations}
\partial_t \vec{\bar{u}}_q^{'''} &:= \partial_t \vec{\bar{u}}_q^{\bigO(1)}, \nonumber\\[2mm]
\hookrightarrow \; \partial_t \vec{\bar{u}}_q^{''}  &= \left(\vec{1}-\vec{\alpha}_q^{\bigO(2)}\right) \odot \partial_t \vec{\bar{u}}_q^{'''} + \vec{\alpha}_q^{\bigO(2)} \odot \partial_t \vec{\bar{u}}_q^{\bigO(2)}, \nonumber\\[2mm]
\hookrightarrow \; \partial_t \vec{\bar{u}}_q^{'}   &= \left(\vec{1}-\vec{\alpha}_q^{\bigO(4)}\right) \odot \partial_t \vec{\bar{u}}_q^{''}  + \vec{\alpha}_q^{\bigO(4)} \odot \partial_t \vec{\bar{u}}_q^{\bigO(4)},\\[2mm]
\hookrightarrow \; \partial_t \vec{\bar{u}}_q       &= \left(\vec{1}-\vec{\alpha}_q^{\bigO(8)}\right) \odot \partial_t \vec{\bar{u}}_q^{'}   + \vec{\alpha}_q^{\bigO(8)} \odot \partial_t \vec{\bar{u}}_q^{\bigO(8)}, \nonumber
\end{align}
where $\odot$ denotes the component-wise multiplication with the vector of
blending parameters. Each DG sub-element computes its own blending factor
according to the local data, i.e. we get 4 blending factors for the $\bigO(2)$
variant, 2 blending factors for $\bigO(4)$ and 1 for $\bigO(8)$\\
\begin{tabular}{cccc}
$\vec{\alpha}_q^{\bigO(2)} =$ & $\left(\alpha_1^{\bigO(2)},\alpha_{2}^{\bigO(2)},\alpha_{3}^{\bigO(2)},\alpha_{4}^{\bigO(2)}\right)_q^T$ &$\otimes$& $(1,1)^T$,\\[2mm]
$\vec{\alpha}_q^{\bigO(4)} =$ & $\left(\alpha_1^{\bigO(4)},\alpha_{2}^{\bigO(4)}\right)_q^T$ &$\otimes$& $(1,1,1,1)^T$,\\[2mm]
$\vec{\alpha}_q^{\bigO(8)} =$ & $\left(\alpha^{\bigO(8)}\right)_q^T$ &$\otimes$& $(1,1,1,1,1,1,1,1)^T$.
\end{tabular}\\
The actual strategy for the blending factors is described in the next section \ref{sec:blending-parameter}.

Again, as discussed in the case of a single-blending approach in section
\ref{sec:convex-blending}, special care is needed to preserve the conservation
of the final discretization. Expanding on the idea in the single-level case, we
aim to compute unique surface fluxes for each discretization with a uniquely
defined interface value at sub-element or subcell interfaces. We evaluate all
high order DG polynomials at the FV subcell interfaces $\vec{u}_q^{\pm\bigO(n)}
\in \reals^{N}$ 
\begin{align*}
\vec{u}_q^{\pm\bigO(1)} &:= \vec{\bar{u}}_q^{\bigO(1)}, &\vec{u}_q^{\pm\bigO(2)} &= \left[\idmat_4 \otimes \mat{I}^{\pm\bigO(2)} \right]\,\vec{\tilde{u}}_q^{\bigO(2)},\\[2mm]
\vec{u}_q^{\pm\bigO(4)} &= \left[\idmat_2 \otimes \mat{I}^{\pm\bigO(4)} \right]\,\vec{\tilde{u}}_q^{\bigO(4)}, &\vec{u}_q^{\pm\bigO(8)} &= \left[\idmat_1 \otimes \mat{I}^{\pm\bigO(8)}\right]\,\vec{\tilde{u}}_q^{\bigO(8)},
\end{align*}
where the interpolation operator $\mat{I}^{\pm\bigO(n)} \in \reals^{n \times
n}$ prolongates to the embedded FV subcell interfaces
$\mu^{\bigO(n)}_{i\pm1/2}$ of the respective DG sub-element, i.e.
\begin{equation}
\label{eq:subcell-fv-interpolation-operator}
\left(\mat{I}^{\pm\bigO(n)}\right)_{ij} = \ell^{\bigO(n)}_j\left(\mu^{\bigO(n)}_{i\pm1/2}\right), \quad i,j = 1,\ldots,n.
\end{equation}
In Fig. \ref{fig:multilevel-blending} in appendix \ref{sec:appendix-visu-multilevel-blending} two representatives of
$\mu^{\bigO(n)}_{i\pm1/2}$ are shown which are aligned at the dotted vertical lines, indicating the interfaces of the FV subcells.

We start with the first order interpolation and stack on top the next levels of
orders via convex blending
\begin{align*}
\left(\vec{u}_q^\pm\right)^{''''} &:= \vec{u}_q^{\pm\bigO(1)}\\[2mm]
\hookrightarrow \; \left(\vec{u}_q^\pm\right)^{'''}   &= \left[\vec{1}-\vec{\alpha}_q^{\bigO(2)}\right] \odot \left(\vec{u}_q^\pm\right)^{''''} + \vec{\alpha}_q^{\bigO(2)} \odot \vec{u}_q^{\pm\bigO(2)} \\[2mm]
\hookrightarrow \: \left(\vec{u}_q^\pm\right)^{''}    &= \left[\vec{1}-\vec{\alpha}_q^{\bigO(4)}\right] \odot \left(\vec{u}_q^\pm\right)^{'''} + \vec{\alpha}_q^{\bigO(4)}  \odot \vec{u}_q^{\pm\bigO(4)} \\[2mm]
\hookrightarrow \: \left(\vec{u}_q^\pm\right)^{'}     &= \left[\vec{1}-\vec{\alpha}_q^{\bigO(8)}\right] \odot \left(\vec{u}_q^\pm\right)^{''} + \vec{\alpha}_q^{\bigO(8)}   \odot \vec{u}_q^{\pm\bigO(8)}.
\end{align*}
To arrive at a complete vector of interface values we append the element
interface values from the left and right neighbors $q-1$ and $q+1$.
\begin{align*}
\vec{u}_q^+ &= \Big( \overbrace{\left(u^+_{N}\right)_{q-1}^{'},}^{\text{left neighbor}} \quad \overbrace{\left(u^+_{1}, \ldots, u^+_{N}\right)_q^{'}}^{\text{internal boundaries}}\Big)^T,\\[2mm]
\vec{u}_q^- &= \Big(\underbrace{\left(u^-_{1}, \ldots, u^-_{N}\right)_q^{'},}_{\text{internal boundaries}} \quad \underbrace{\left(u^-_{1}\right)_{q+1}^{'}}_{\text{right neighbor}}\Big)^T.
\end{align*}
These new common interface values allow us to evaluate a common surface flux
\eqref{eq:common-surface-flux} at each interface, which is again the key to
preserve conservation in the final discretization. This concludes the description
of the multi-level blending scheme.

\begin{changed}
\subsection{Calculation of the Blending Factor $\alpha$}
\label{sec:blending-parameter}
A good shock indicator for high order methods is supposed to recognize
discontinuities, such as weak and strong shocks, in the solution early on and
mark the affected elements for proper shock capturing. On the other hand, the
indicator should avoid to flag non-shock related fluid
features such as shear layers and turbulent flows. There are numerous
indicators available for DG schemes
\cite{persson2006sub,Dumbser2014,Dumbser2016,Vilar2019,Bohm2019aa,Ching2019,Zhu2009}. 


In this work we want to construct an a-priori indicator which relies on the
readily available information within each element provided by the multi-level
blending framework introduced in section \ref{sec:multilevel-blending}. The
principal idea is to compare a measure of smoothness of the different order reconstructions 
with each other. Smooth, well-resolved flows are expected to yield rather similar solution profiles compared to 
data that contain strong variations. 

The smoothness measure $\sigma_q^{\bigO(n)}$ within element $q$ of the $n$-th
order reconstruction $\tilde{u}_q^{\bigO(n)}$ is inferred from the $L^1$-norm
of its first derivative. For example, if we want to calculate the blending
factor $\alpha_q^{\bigO(8)}$ for the eighth order we compute
\begin{equation}
\label{eq:L1-norm-slope}
    \sigma_q^{\bigO(8)} = \int_{-\frac{1}{2}}^{\frac{1}{2}} \left| \partial_{\xi} \tilde{u}_q^{\bigO(8)}(\xi) \right|d\xi
\end{equation}
and
\begin{equation*}
    \sigma_q^{\bigO(4)}
        = \left(\sigma_1^{\bigO(4)}\right)_q + \left(\sigma_2^{\bigO(4)}\right)_q
         = \int_{-\frac{1}{2}}^{0} \left| \partial_{\xi} \left(\tilde{u}_1^{\bigO(4)}\right)_q(\xi) \right|d\xi
        + \int_{0}^{\frac{1}{2}} \left| \partial_{\xi} \left(\tilde{u}_2^{\bigO(4)}\right)_q(\xi) \right|d\xi        
\end{equation*}
utilizing information from the top-level eighth order element and its two
fourth order sub-elements. The integrals are evaluated with appropriate
quadrature rules. Then the blending factor is calculated by comparing the
smoothness of the different orders as
\begin{equation}
\label{eq:calc-blending-factor}
    \alpha_q^{\bigO(8)} = 1 - \text{cutoff}\left(0,
        \tau_A \, \left(\frac{\left|\sigma_q^{\bigO(8)}
        - \sigma_q^{\bigO(4)}\right|}{\max\big(\sigma_q^{\bigO(4)},\sigma_q^{\bigO(8)},1\big)}
        - \tau_S\right), 1\right),
\end{equation}
where $\text{cutoff}(x_{\text{lower}},x,x_{\text{upper}}) =
\min(\max(x_{\text{lower}},x),x_{\text{upper}})$ and with two parameters
$\tau_A, \tau_S > 0.$ The design of equation \eqref{eq:calc-blending-factor}
ensures that $\alpha_q^{\bigO(n)}$ is always within the unit interval and that
for very small $\sigma_q^{\bigO(n)}$ the indicator does not get hypersensitive
due to floating point truncation. If not noted otherwise we set the
amplification parameter to $\tau_A = 100$ and the sensitivity parameter $\tau_S
= 0.05$. The latter parameter has the biggest influence on the behavior of the
indicator and its value demarcates the lower bound where it does not interfere
with the convergence tests conducted in section \ref{sec:5}. In order to
capture all troubling flow features we independently apply the indicator
\eqref{eq:calc-blending-factor} on all primitive state variables and select the
smallest of the resulting blending factors. The calculations of the two fourth
order blending factors $\left(\alpha_1^{\bigO(4)}\right)_q$ and
$\left(\alpha_2^{\bigO(4)}\right)_q$ are done within each sub-element
independently together with their respective second order smoothness measures
$\left(\sigma_s^{\bigO(2)}\right)_q$, $s = 1,\ldots,4$. \\

For piecewise linear polynomials, the indicator \eqref{eq:calc-blending-factor}
is not applicable. Instead, we directly use the squeezing parameter $\beta$
obtained from the positivity limiter in \eqref{eq:positivity-limiter}, i.e.
$\alpha_s^{\bigO(2)} := \beta_s^{\bigO(2)}$, $s = 1,\ldots,4$.\\

\noindent\textbf{Remark:} The indicator \eqref{eq:calc-blending-factor} is also used for the single-level
blending scheme.
\end{changed}

\begin{changed}
\subsection{Sketch of the Algorithm}
\label{sec:sketch-of-alogrithm}
In the last part of this section, a sketch of the algorithm for the blending
framework is presented.  We present a general outline of the necessary steps
for the 1D single-blending scheme evolved with a multi-stage Runge-Kutta time
stepping method. We enter at the beginning of a Runge-Kutta cycle and do the following.
\begin{enumerate}[label=(\Roman*),start=1,itemsep=2mm,leftmargin=1cm]
\item \label{item:algo-single-level} Reconstruct the polynomial
$\vec{\tilde{u}}_q$ from the given mean values $\vec{\bar{u}}_q$ for each
element $q$ as in \eqref{eq:trafo-mv-nodal}.

\item If the reconstructed polynomial $\vec{\tilde{u}}_q$ contains
non-permissible states, see \eqref{eq:positivity-limiter-2}, then
calculate the limited version $\vec{\tilde{u}}_q^{(\beta)}$ as in
\eqref{eq:positivity-limiter}

\item If the squeezing parameter $\beta_q$ is below $\beta_L$ then set
$\alpha_q := 0$ else compute the blending factor $\alpha_q$ via
\eqref{eq:calc-blending-factor} from the \textbf{unlimited} polynomial
$\vec{\tilde{u}}_q$.

\item Compute the boundary values $u_{q}^{\pm}$ via
\eqref{eq:inner-surface-prolongation} and exchange alongside zone boundaries in
case of distributed computing.

\item Determine the common surface flux $\vec{f}_{q}^*$
via \eqref{eq:common-surface-flux}.

\item Calculate the right-hand-side $\partial_t \vec{\bar{u}}_q$:
    \begin{itemize}[label=$\triangleright$,topsep=2mm,itemsep=2mm,leftmargin=5mm]
        \item If the blending factor $\alpha_q$ above $\alpha_H$ then compute
            $\partial_t \vec{\bar{u}}_q$ with the DG-only scheme.

        \item Else if the blending factor $\alpha_q$ below $\alpha_L$ then compute
            $\partial_t \vec{\bar{u}}_q$ with the FV-only scheme.

        \item Else compute $\partial_t \vec{\bar{u}}_q$ with the single-level blending scheme \eqref{eq:dgfv-scheme-1d}.
    \end{itemize}
\item \label{item:time-int-single-level} Forward in time to the next Runge-Kutta stage and return to step \ref{item:algo-single-level}.
\end{enumerate}
The switching thresholds are set to $\alpha_H := 0.99$ and $\alpha_L := 0.01$
and the limiter threshold to $\beta_L := 0.95$. Note that the algorithm only
applies the blending procedure where necessary in order to maintain the overall
performance of the scheme.\\

This concludes the presentation of the 1D blending scheme.  The description of
the blending scheme on 3D Cartesian meshes as well as the algorithm outline for
the multi-level blending scheme can be found in appendix
\ref{sec:appendix-3D-extension} and \ref{sec:appendix-sketch-of-alogrithm}.
\end{changed}

\section{Validation}
\label{sec:5}
For the computational investigations, the multi-level algorithm with the
explicit SSP-RK(5,4) time integrator \cite{spiteri2002new} is implemented in a
Fortran-2008 prototype code with a hybrid parallelization strategy based on MPI
and OpenMP. Management of the AMR and load balancing is provided by \pforest, a
highly efficient Octree library \cite{BursteddeWilcoxGhattas11}.
\begin{changed}
The maximal time step for dimensions $d = 1,2,3$ is estimated by the CFL
condition
\begin{equation}
\label{eq:cfl-condition}
    \Delta t := CFL \, \min_q\;\frac{\Delta x_q}{2^{d-1}(2\,N-1)\,\bar{\lambda}_{q,\text{max}}},
\end{equation}
where $CFL := 0.8$, $N := 8$ is the number of mean values in each direction of
the element $q$ and $\bar{\lambda}_{q,\text{max}}$ is an estimate of the
maximum eigenvalue given in \eqref{eq:euler-max-eigenvalue}. For the numerical
interface fluxes $f^*$, we use the Harten-Lax-Leer (HLL) approximate Riemann
solver \cite{harten1983upstream} with Einfeldt signal speed estimates
\cite{einfeldt1988godunov}.
\end{changed}

\subsection{Governing Equations}
\label{sec:gov-equations}
We consider the compressible Euler equations 
\begin{align}
\label{eqn:comp-euler-conservative}
\partial_t\dens + \nabla\cdot (\dens\,\velv)   &=  0, \nonumber \\[0.2cm]
\partial_t(\dens\,\velv) + \nabla\cdot \left(\dens\,\velv\velv^T + \pres\,\idmat\right) &= \vec{0},\\[0.3cm]
\partial_t \Ener + \nabla\cdot \Big(\velv\,(\Ener + \pres)\Big) &=  0, \nonumber
\end{align}
with the vector of conserved quantities $\statevec{u} = (\dens,\dens\,\velv,
\Ener)^T$, where $\dens$ denotes the density, $\velv$ the velocity, and $\Ener$
the total energy. We assume a perfect gas equation of state and compute
the pressure as
\begin{equation}
\label{eq:equation-of-state}
p(\statevec{u}) = (\gamma - 1)\,\left[\Ener - \frac{\dens}{2}\,\velv^T\velv\right].
\end{equation}
If not stated otherwise we choose $\gamma = 1.4$. The set of permissible states is given by
\begin{equation}
    \big\{\text{permissible states}\big\} = \big\{\forall\,\vec{u} \; \big| \; \dens > 0 \wedge \pres(\vec{u}) > 0\big\}.
\end{equation}
\begin{changed}
For the CFL condition \eqref{eq:cfl-condition} the maximum eigenvalue is
evaluated on all mean values $\vec{\bar{u}}_q$ of element $q$. Given dimension $d
= \{1,2,3\}$ it reads as
\begin{equation}
\label{eq:euler-max-eigenvalue}
    \bar{\lambda}_{q,\text{max}} = \max_{i,d} \left(\big|(\bar{v}_d)_i\big| + \sqrt{\gamma\,\frac{\bar{p}_i}{\bar{\rho}_i}}\right), \quad i = 1,\ldots,N^{d_{\text{max}}}.
\end{equation}
\end{changed}

\subsection{Convergence Test}
\label{sec:conv-test}
We use the manufactured solution method \cite{roy2004verification} and validate
the 3D multi-level blending scheme on a periodic cube of unit length ($L = 1$)
where the resolution of the mesh is incrementally doubled.
We define our manufactured solution in primitive state variables as
\begin{changed}
\begin{align}
    \label{eq:manufactured-solution}
    \dens(t; x,y,z) &= 1.0 + 0.35 \, \sin\left(\frac{2\,\pi\,(x-t)}{L}\right) + 0.24 \,\cos\left(\frac{2\,\pi\,(y-t)}{L}\right) + 0.1 \,\sin\left(\frac{2\,\pi\,(z-t)}{L}\right), \nonumber \\
    \velv(t; x,y,z) &= (0.1,0.2,0.3)^T, \\
    \pres(t; x,y,z) &= 1.0 + 0.23 \, \cos\left(\frac{2\,\pi\,(x-t)}{L}\right) + 0.19\,\sin\left(\frac{2\,\pi\,(y-t)}{L}\right) + 0.2 \,\cos\left(\frac{2\,\pi\,(z-t)}{L}\right). \nonumber
\end{align}
\end{changed}
The final time of the simulation is $T = 2$ and the center of the domain is
refined to introduce non-conforming interfaces in the computational domain. We
determine the  $L^{\infty}$- and $L^2$-norms of the errors in the density and
total energy. Tables \ref{tab:eoc-order-1} to \ref{tab:eoc-order-8} list the
results for the first, second, fourth and eighth order multi-level blending
schemes. The initial conditions and the source term of our manufactured
solution \eqref{eq:manufactured-solution} is in all cases evaluated and applied
on the mean values via an appropriate quadrature rule to maintain high order.
The results confirm, that the discretizations behave as designed in this
assessment.
\begin{table}[H]
\captionsetup{width=0.8\textwidth}
\caption{Experimental order of convergence of the first-order FV variant within
the 3D multi-level blending framework.}
\centering
\begin{tabular}{c|cc|cc|cc|cc}
\toprule
resol. & \multicolumn{2}{c}{density} & \multicolumn{2}{|c}{total energy}
       & \multicolumn{2}{|c}{density} & \multicolumn{2}{|c}{total energy} \\
\midrule
$2^r$ &
    $||\epsilon||_{\infty}$ & $||\epsilon||_{2}$ &
    $||\epsilon||_{\infty}$ & $||\epsilon||_{2}$ &
    $EOC_{\infty}$ & $EOC_{2}$ &
    $EOC_{\infty}$ & $EOC_{2}$ \\
\midrule
 16 & 5.98e-03 & 4.84e-03 & 1.34e-02 & 1.03e-02 & n/a & n/a & n/a & n/a \\
 32 & 3.51e-03 & 2.83e-03 & 7.91e-03 & 6.13e-03 &  0.77 &  0.77 &  0.76 &  0.75 \\
 64 & 2.21e-03 & 1.45e-03 & 4.30e-03 & 3.13e-03 &  0.67 &  0.97 &  0.88 &  0.97 \\
128 & 1.42e-03 & 7.30e-04 & 2.26e-03 & 1.57e-03 &  0.63 &  0.99 &  0.93 &  0.99 \\
\bottomrule
\end{tabular}
\label{tab:eoc-order-1}
\end{table}
\vspace{-5mm}
\begin{table}[H]
\captionsetup{width=0.8\textwidth}
\caption{Experimental order of convergence of the second order DG variant
within the 3D multi-level blending framework.}
\centering
\begin{tabular}{c|cc|cc|cc|cc}
\toprule
resol. & \multicolumn{2}{c}{density} & \multicolumn{2}{|c}{total energy}
       & \multicolumn{2}{|c}{density} & \multicolumn{2}{|c}{total energy} \\
\midrule
$2^r$ &
$||\epsilon||_{\infty}$ & $||\epsilon||_{2}$ &
$||\epsilon||_{\infty}$ & $||\epsilon||_{2}$ &
$EOC_{\infty}$ & $EOC_{2}$ &
$EOC_{\infty}$ & $EOC_{2}$ \\
\midrule
 16 & 1.46e-03 & 5.38e-04 & 2.21e-03 & 9.35e-04 & n/a & n/a & n/a & n/a \\
 32 & 4.59e-04 & 9.29e-05 & 6.63e-04 & 1.97e-04 &  1.67 &  2.53 &  1.74 &  2.25 \\
 64 & 1.11e-04 & 1.77e-05 & 1.48e-04 & 4.26e-05 &  2.04 &  2.39 &  2.17 &  2.21 \\
128 & 2.63e-05 & 3.88e-06 & 3.09e-05 & 9.72e-06 &  2.08 &  2.19 &  2.25 &  2.13 \\
\bottomrule
\end{tabular}
\label{tab:eoc-order-2}
\end{table}
\vspace{-5mm}
\begin{table}[H]
\captionsetup{width=0.8\textwidth}
\caption{Experimental order of convergence of the fourth order DG variant
within the 3D multi-level blending framework.}
\centering
\begin{tabular}{c|cc|cc|cc|cc}
\toprule
resol. & \multicolumn{2}{c}{density} & \multicolumn{2}{|c}{total energy}
       & \multicolumn{2}{|c}{density} & \multicolumn{2}{|c}{total energy} \\
\midrule
$2^r$ &
$||\epsilon||_{\infty}$ & $||\epsilon||_{2}$ &
$||\epsilon||_{\infty}$ & $||\epsilon||_{2}$ &
$EOC_{\infty}$ & $EOC_{2}$ &
$EOC_{\infty}$ & $EOC_{2}$ \\
\midrule
 16 & 1.10e-04 & 6.45e-05 & 3.20e-04 & 1.83e-04 & n/a & n/a & n/a & n/a \\
 32 & 7.86e-06 & 4.33e-06 & 2.19e-05 & 1.33e-05 &  3.80 &  3.90 &  3.87 &  3.79 \\
 64 & 5.18e-07 & 2.59e-07 & 1.31e-06 & 7.59e-07 &  3.92 &  4.07 &  4.07 &  4.13 \\
128 & 3.45e-08 & 1.82e-08 & 7.75e-08 & 4.29e-08 &  3.91 &  3.83 &  4.07 &  4.15 \\
\bottomrule
\end{tabular}
\label{tab:eoc-order-4}
\end{table}
\vspace{-5mm}
\begin{table}[H]
\captionsetup{width=0.8\textwidth}
\caption{Experimental order of convergence of the eighth order DG variant
within the 3D multi-level blending framework. Here we set $CFL := 0.1$, as the
time integration method is only fourth order accurate.}
\centering
\begin{tabular}{c|cc|cc|cc|cc}
\toprule
resol. & \multicolumn{2}{c}{density} & \multicolumn{2}{|c}{total energy}
       & \multicolumn{2}{|c}{density} & \multicolumn{2}{|c}{total energy} \\
\midrule
$2^r$ &
$||\epsilon||_{\infty}$ & $||\epsilon||_{2}$ &
$||\epsilon||_{\infty}$ & $||\epsilon||_{2}$ &
$EOC_{\infty}$ & $EOC_{2}$ &
$EOC_{\infty}$ & $EOC_{2}$ \\
\midrule
 16 & 2.77e-07 & 1.65e-07 & 7.77e-07 & 4.55e-07 & n/a & n/a & n/a & n/a \\
 32 & 2.77e-09 & 1.00e-09 & 6.72e-09 & 2.41e-09 &  6.65 &  7.36 &  6.85 &  7.56 \\
 64 & 1.73e-11 & 3.48e-12 & 4.15e-11 & 1.04e-11 &  7.32 &  8.17 &  7.34 &  7.86 \\
128 & 8.06e-14 & 1.87e-14 & 1.85e-13 & 6.74e-14 &  7.74 &  7.54 &  7.81 &  7.27 \\
\bottomrule
\end{tabular}
\label{tab:eoc-order-8}
\end{table}

\subsection{Conservation Test}
\begin{figure}[H]
\centering
\ifnum \compilewithplots = 1
\includegraphics[width=0.6\textwidth]{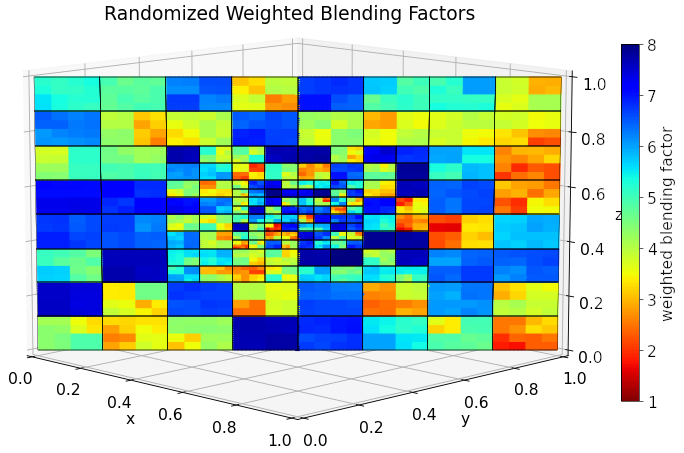}
\fi
\caption{Conservation test of the 3D multi-level blending scheme: Slice through
the computational domain showing weighted blending factor $\bar{\alpha}$
according to \eqref{eq:order-estimate} on a Cartesian non-conforming mesh with
refinement levels 3 to 5. The black lines depict the boundaries of the corresponding 8th order DG elements.}
\label{fig:3d-consistency-order-estimate}
\end{figure}

The goal in this assessment is to demonstrate that the multi-level-blending
discretization is conservative for all choices of blending factors. To do so,
we adapt the same setup as in section \ref{sec:conv-test} but deactivate the
source term. The center of the domain is refined to introduce non-conforming
interfaces in the computational domain. Additionally, as a stress test, the
blending factors are randomly chosen and changed after each Runge-Kutta stage.
As there are multiple blending factors at a given spatial location, we consider
the following weighted blending factor 
\begin{equation}
\label{eq:order-estimate}
    \bar{\alpha} = \left\{\left[\left(1-\alpha^{\bigO(2)}\right) + 2\,\alpha^{\bigO(2)}\right]\,\left(1-\alpha^{\bigO(4)}\right) + 4\,\alpha^{\bigO(4)}\right\}\,\left(1-\alpha^{\bigO(8)}\right) + 8\,\alpha^{\bigO(8)},
\end{equation}
to illustrate the distribution of the blending factors in
Fig.\ref{fig:3d-consistency-order-estimate}. Note the limiting cases
$\bar{\alpha}=1$ for pure FV and $\bar{\alpha}=8$ for the eighth order DG
scheme.

\begin{figure}[H]
\centering
\ifnum \compilewithplots = 1
\includegraphics[width=0.7\textwidth]{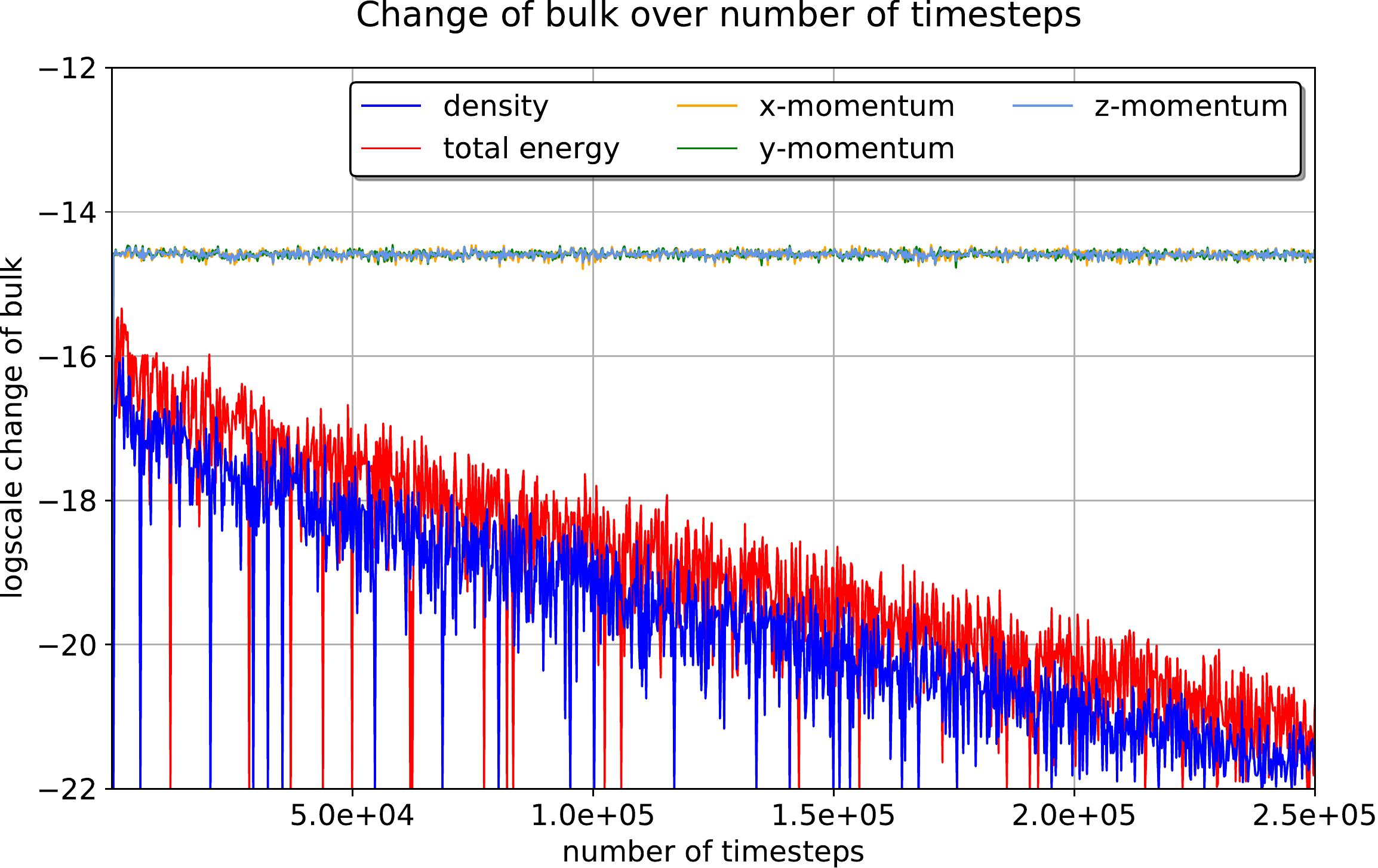}
\fi
\caption{Conservation test of the 3D multi-level blending scheme: Evolution of
the integrated change of bulk $\log_{10}(|\partial_t \block{\bar{u}}_{\text{total}}|)$ (eq.
\eqref{eq:change-of-bulk}) of each conservative state variable.}
\label{fig:3d-consistency-evolution}
\end{figure}

\begin{changed}
The simulation runs to $T=300$ performing more than a quarter million
timesteps. The result of the test is shown in Fig.
\ref{fig:3d-consistency-evolution} where we plot in log-scale the absolute
value of the \emph{change of bulk} $\partial_t
\block{\bar{u}}_{\text{total}}(t)$, integrated over the whole domain,
\begin{equation}
\label{eq:change-of-bulk}
    \partial_t \block{\bar{u}}_{\text{total}}(t) = \sum_q^Q \frac{1}{|\Omega_q|} \int_{\Omega_q} \partial_t \block{\bar{u}}(t) \, d\Omega.
\end{equation}
$Q$ is the total number of elements and $|\Omega_q| = \Delta x_q\Delta
y_q\Delta z_q$ is the volume per element. The results show that the conservation error lies within the
range of 64 bit (double precision) floating point truncation and hence confirm that the multi-level blending discretization is fully conservative up to machine precision errors. 
\end{changed}


\begin{changed}
\subsection{1D Shock Tube Problems}
\label{sec:1d-shock-tubes}
We validate the multi-level blending scheme
with three well-established 1D problems, namely Sod, Lax and Shu-Osher shock
tubes. But first, in order to gain insights into the individual contributions
of the different schemes, we introduce a more intuitive reformulation of the
blending factors $\alpha^{\bigO(n)}$. The multi-level blending operation
\eqref{eq:multilevel-blending-discretizations} within element $q$ can be
interpreted as a linear superposition of four numerical schemes. For that, we
collapse the following relations
\begin{align*}
\partial_t \vec{\bar{u}}_q^{'''} &:= \partial_t \vec{\bar{u}}_q^{\bigO(1)}, \nonumber\\[2mm]
\hookrightarrow \; \partial_t \vec{\bar{u}}_q^{''}  &= \left(\vec{1}-\vec{\alpha}_q^{\bigO(2)}\right) \odot \partial_t \vec{\bar{u}}_q^{'''} + \vec{\alpha}_q^{\bigO(2)} \odot \partial_t \vec{\bar{u}}_q^{\bigO(2)}, \nonumber\\[2mm]
\hookrightarrow \; \partial_t \vec{\bar{u}}_q^{'}   &= \left(\vec{1}-\vec{\alpha}_q^{\bigO(4)}\right) \odot \partial_t \vec{\bar{u}}_q^{''}  + \vec{\alpha}_q^{\bigO(4)} \odot \partial_t \vec{\bar{u}}_q^{\bigO(4)},\\[2mm]
\hookrightarrow \; \partial_t \vec{\bar{u}}_q       &= \left(\vec{1}-\vec{\alpha}_q^{\bigO(8)}\right) \odot \partial_t \vec{\bar{u}}_q^{'}   + \vec{\alpha}_q^{\bigO(8)} \odot \partial_t \vec{\bar{u}}_q^{\bigO(8)}, \nonumber
\end{align*}
into one single term and define \emph{blending weights} $\theta^{\bigO(n)}$,
i.e.
\begin{equation*}
\partial_t \vec{\bar{u}}_q = \vec{\theta}_q^{\bigO(1)} \odot \partial_t \vec{\bar{u}}_q^{\bigO(1)} +
    \vec{\theta}_q^{\bigO(2)} \odot \partial_t \vec{\bar{u}}_q^{\bigO(2)} +
    \vec{\theta}_q^{\bigO(4)} \odot \partial_t \vec{\bar{u}}_q^{\bigO(4)} +
    \vec{\theta}_q^{\bigO(8)} \odot \partial_t \vec{\bar{u}}_q^{\bigO(8)},
\end{equation*}
where
\begin{align}
\label{eq:blending-weights}
\vec{\theta}_q^{\bigO(1)} &:= \left(\vec{1}-\vec{\alpha}_q^{\bigO(8)}\right) \odot \left(\vec{1}-\vec{\alpha}_q^{\bigO(4)}\right) \odot\left(\vec{1}-\vec{\alpha}_q^{\bigO(2)}\right), \nonumber\\[2mm]
\vec{\theta}_q^{\bigO(2)} &:= \left(\vec{1}-\vec{\alpha}_q^{\bigO(8)}\right) \odot \left(\vec{1}-\vec{\alpha}_q^{\bigO(4)}\right) \odot \vec{\alpha}_q^{\bigO(2)}, \nonumber\\[2mm]
\vec{\theta}_q^{\bigO(4)} &:= \left(\vec{1}-\vec{\alpha}_q^{\bigO(8)}\right) \odot \vec{\alpha}_q^{\bigO(4)}, \\[2mm]
\vec{\theta}_q^{\bigO(8)} &:= \vec{\alpha}_q^{\bigO(8)}. \nonumber
\end{align}
By construction the blending weights have the following property
\begin{equation}
\label{eq:blending-weights-unit-property}
    \vec{\theta}_q^{\bigO(1)} + \vec{\theta}_q^{\bigO(2)} + \vec{\theta}_q^{\bigO(4)} + \vec{\theta}_q^{\bigO(8)} = \vec{1},
\end{equation}
and thus give a proper fraction of each contribution.\\

The first shock tube problem is the Sod shock tube \cite{sod1978survey}. It is
defined on the unit interval $\Omega = [0,1]$ with a diaphragm located at $x_D
= 0.5$. The initial condition in primitive state variables reads
\begin{equation*}
    \big(\rho_0(x),v_0(x),\pres_0(x)\big) = \begin{lcases}
    \big(1, 0, 1\big), \quad x < x_D,\\
    \big(0.125, 0, 0. 1\big), \quad x \geq x_D.
\end{lcases}
\end{equation*}
The resolution is set to 32 elements of $N = 8$ mean values each. This amounts
to 256 total DOF. The result for the density profile (top row) at the final
simulation time $T = 0.2$ is presented in Fig. \ref{fig:1d-sod-shock-tube}
together with the exact solution. It shows the correct approximation of
the rarefaction wave, contact discontinuity and the forward facing shock front.
As designed, only at the shock front the blending scheme gets triggered, which
is visible in the bottom row of the plot. The stacked bar chart directly
corresponds to the blending weights $\theta^{\bigO(n)}$ in
\eqref{eq:blending-weights}. The vertical dimension of the stacked bars
completely fill the unit interval mirroring property 
\eqref{eq:blending-weights-unit-property}.\\
\begin{figure}[H]
\centering
\ifnum \compilewithplots = 1
\includegraphics[width=0.9\textwidth]{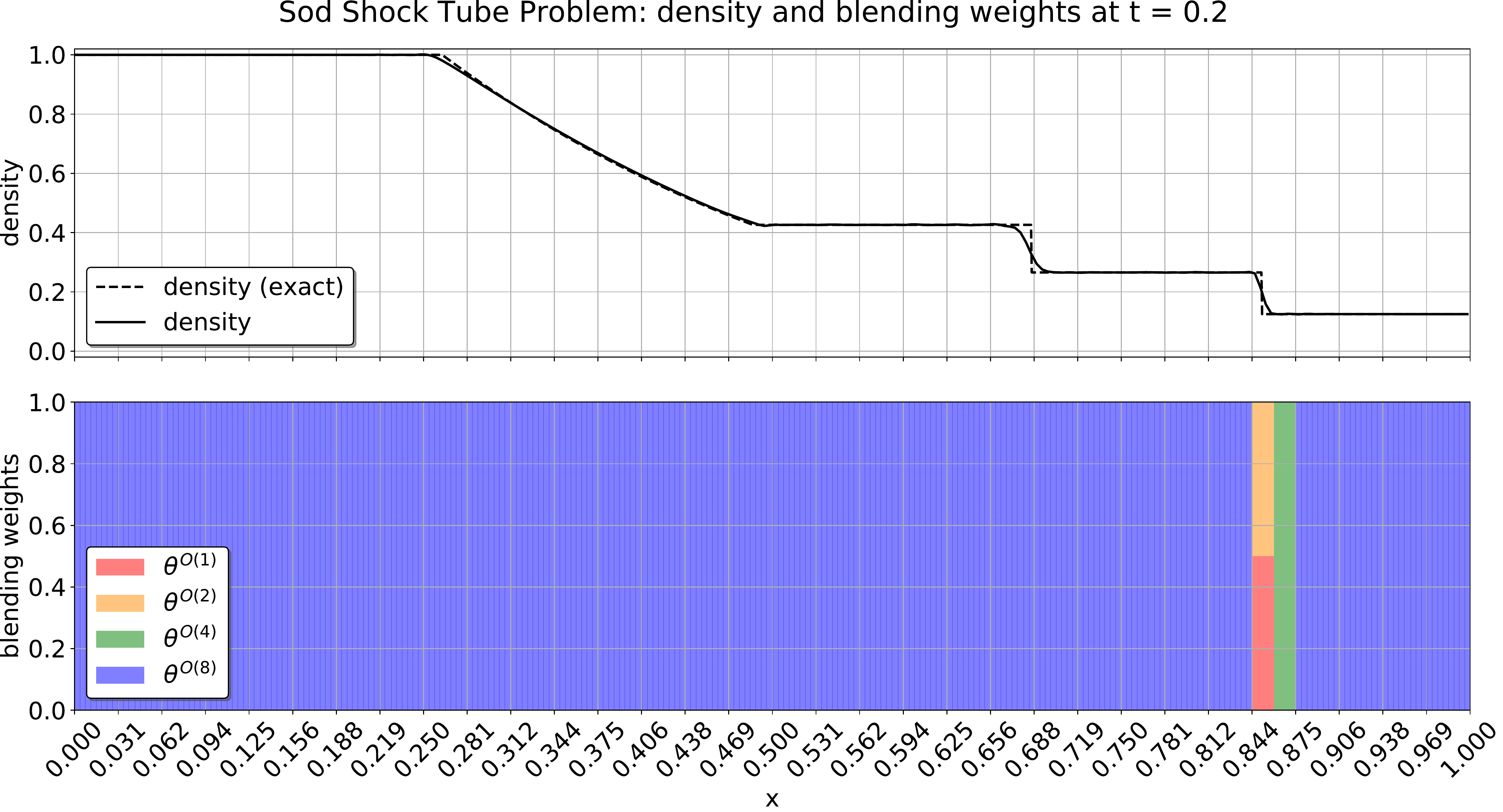}
\fi
\caption{Sod Shock Tube Problem: Numerical solution of the density profile (top
row) with the DGFV8 multi-level blending scheme together with the exact
solution at final simulation time $T = 0.2$. The bottom row shows the blending
weights $\theta^{\bigO(n)}$ (eq. \eqref{eq:blending-weights}) encoded by
stacked bars in different colors. The light grey vertical grid lines depict the element
boundaries.}
\label{fig:1d-sod-shock-tube}
\end{figure}
The second test case is the Lax shock tube \cite{lax1954weak} with initial
data set to
\begin{equation*}
    \big(\rho_0(x),v_0(x),\pres_0(x)\big) = \begin{lcases}
    \big(0.445, 0.689,3.528\big), \quad x < x_D,\\
    \big(0.5,0,0.571\big), \quad x \geq x_D.
\end{lcases}
\end{equation*}
All other simulation parameters as well as the resolution are the same as in
the Sod test case. The result for the density profile (top row) at final
simulation time $T = 0.15$ is presented in Fig. \ref{fig:1d-lax-shock-tube}
together with reference solution on a finer grid of 1024 DOF with the
third-order piece-wise parabolic method (PPM) \cite{colella1984piecewise}
implemented in the astrophysics code FLASH (version 4.6, March 2019), see e.g.
\cite{fryxell2000flash}. For comparison, we also included the solution of the
PPM with the same grid resolution of 256 DOF. All flow features are resolved
correctly and the blending scheme is only triggered in the region around the
forward facing shock. Furthermore, this example clearly reveals the adaptive
blending on the sub-element level.
\begin{figure}[H]
\centering
\ifnum \compilewithplots = 1
\includegraphics[width=0.9\textwidth]{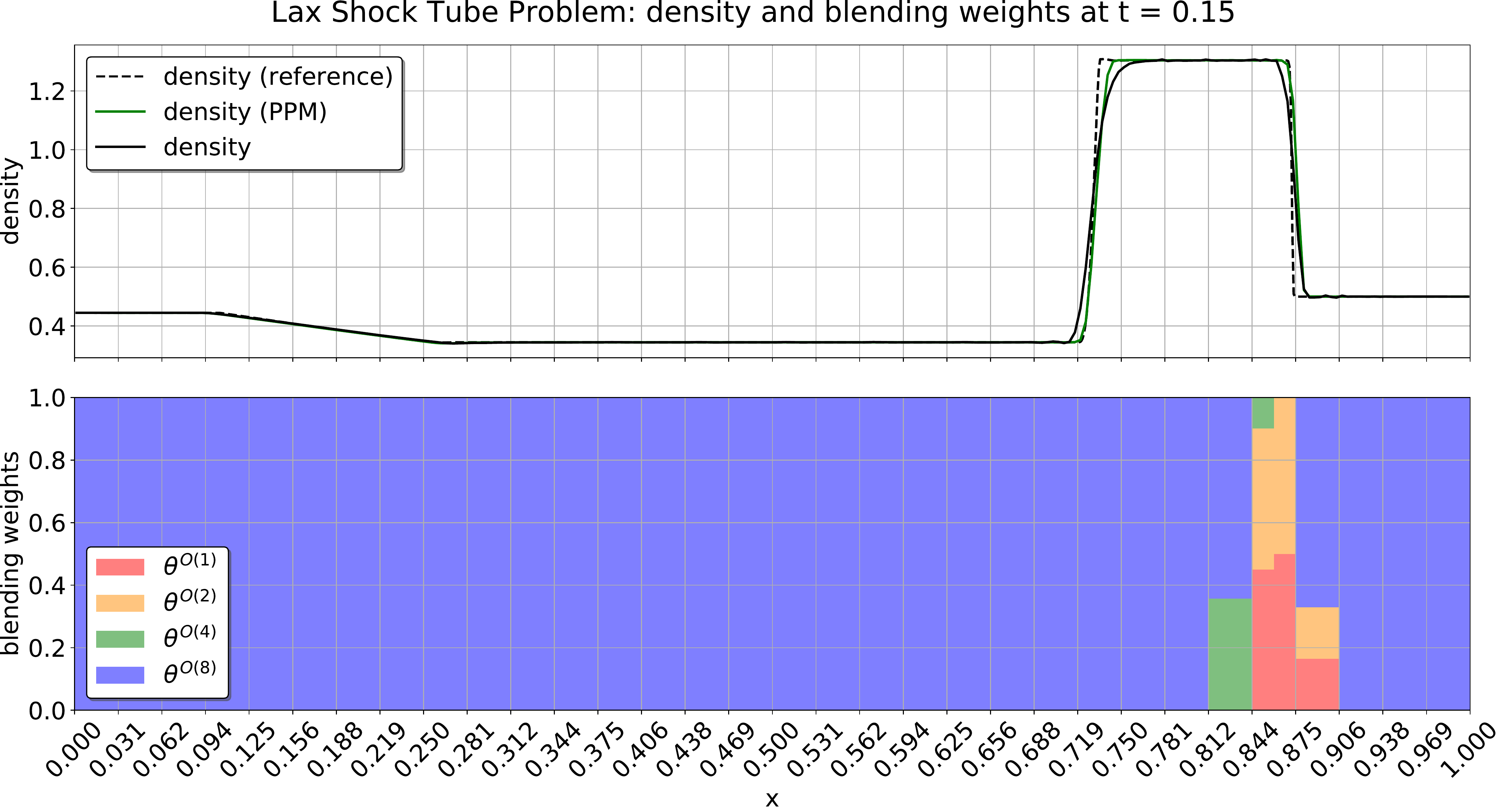}
\fi
\caption{Lax Shock Tube Problem: Numerical solution of the density profile (top
row) with the DGFV8 (256 DOF) multi-level blending scheme together with the
reference solution (PPM, 1024 DOF) and PPM (256 DOF) at final simulation time
$T = 0.2$. The bottom row shows the blending weights $\theta^{\bigO(n)}$ (eq.
\eqref{eq:blending-weights}) encoded by stacked bars in different colors. The
vertical grid lines depict the element boundaries.}
\label{fig:1d-lax-shock-tube}
\end{figure}
The third shock tube, the Shu-Osher test \cite{shu1988efficient}, is a Mach $3$
shock interacting with a sinusoidal density wave. It reveals the scheme's
capability of capturing both, discontinuous and smooth parts of the flow. The
computational domain in 1D is set to $\Omega = [-4.5,4.5]$, the final
simulation time is $T = 1.8$, and the primitive variables are initialized as
\begin{equation*}
    \big(\rho_0(x),v_0(x),\pres_0(x)\big) = \begin{lcases}
    \big(3.8567143, 2.629369, 10.33333\big), \quad x < -4,\\
    \big(1 + 0.2\,\sin(5\,x), 0, 1\big), \quad x \geq -4.
\end{lcases}
\end{equation*}
Here, we compare the results of the multi-level blending DGFV8 to PPM. In order
to investigate the importance of the multi-level approach on the accuracy of
the DGFV8 result, we additionally perform a simulation where we intentionally
deactivate the sub-elements, i.e.  $\vec{\alpha}^{\bigO(2)} := 0$ and
$\vec{\alpha}^{\bigO(4)} := 0$. This is identical to the single-level blending
approach presented in section \ref{sec:convex-blending}, where only the first
order FV scheme is blended with the eighth order DG scheme. The resolution is
as before 256 total DOF whereas the reference solution is computed with PPM on
a much finer grid of 2048 DOF.  The numerical experiments shown in Fig.
\ref{fig:1d-shu-osher-density} nicely demonstrate the benefit of using the multi-level approach. Whereas the shocks are about equally resolved,
the multi-level blending variant gives the best results in the smooth parts of
the solution. Since the resolution is not very high, non-shock related small
scale flow features cannot always be resolved by the eighth order DG scheme.
This is especially visible in the range $x = [1.688,2.25]$ where the lower
order scheme has to completely or at least partially take over.
\begin{figure}[H]
\centering
\ifnum \compilewithplots = 1
\includegraphics[width=0.9\textwidth]{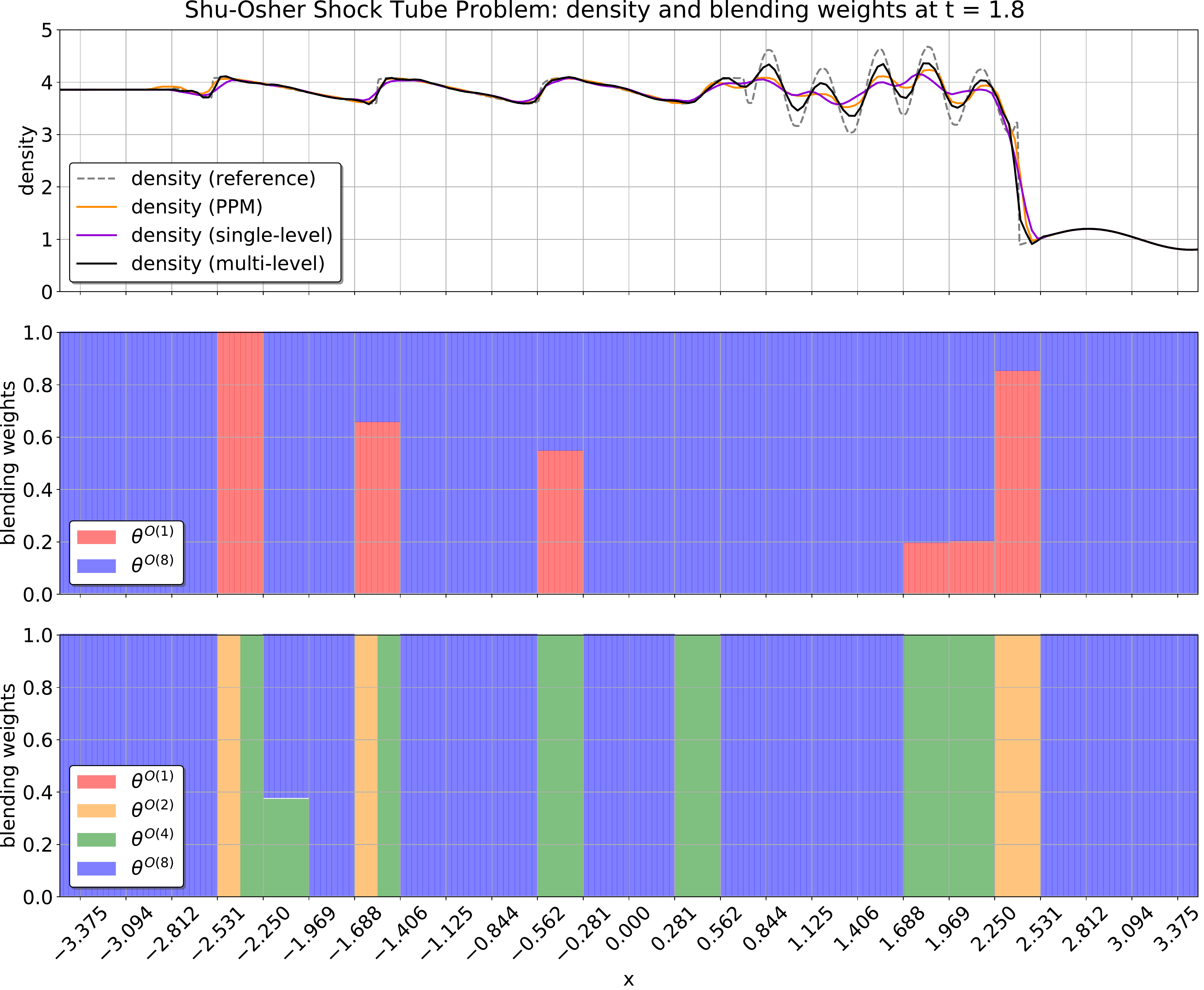}
\fi
\caption{Shu-Osher Shock Tube Problem: Numerical solution of the density
profile (top row) with the DGFV8 multi-level blending scheme together with the
reference solution (PPM on 2048 DOF), PPM and single-level blending scheme all
using 256 DOF. The center row shows the blending weights $\theta^{\bigO(n)}$
(eq. \eqref{eq:blending-weights}) for the single-level blending while the
bottom row shows the multi-level blending. The vertical grid lines
depict the element boundaries.}
\label{fig:1d-shu-osher-density}
\end{figure}
\end{changed}

\begin{changed}
\subsection{2D Riemann Problems}
\label{sec:2d-riemann-problems}
In this section we present a selection of three 2D Riemann problems
\cite{lax1998solution,schulz1993classification,krais2020flexi}. The domain for
all simulations is set to $\Omega = [-0.5,0.5]^2$ with a uniform grid
resolution of $64^2$ elements, respectively $512^2$ DOF. The setup consists of
the four quadrants each initialized with their own constant states. The exact
parameters of each Riemann problem are given in \cite{lax1998solution} where
they are addressed with a fixed configuration number. For brevity we omit the
setup parameters and refer to \cite{lax1998solution}. In this work we show
configurations 3, 4 and 6. We are interested in the change of the numerical
solution when we incrementally add one level of blending order from one run to
the next. Hence, we present three solutions for each 2D Riemann problem denoted
by their respective multi-level schemes: DGFV2, DGFV4 and DGFV8. Note that the
schemes are implemented in such a way that they operate on the same element
size of $N^2 = 8^2$ mean values. The results are shown in Fig
\ref{fig:2d-riemann-problem-dens-blend-3} to
\ref{fig:2d-riemann-problem-dens-blend-6}. The left column
shows the density contour and the right column the weighted blending factors as
in \eqref{eq:order-estimate}. The general observation is that with increasing
order there is more structure visible in the density plots and the blending
patterns get more nuanced in tracing the flow structures.
\begin{figure}[H]
\centering
\ifnum \compilewithplots = 1
\includegraphics[width=0.9\textwidth]{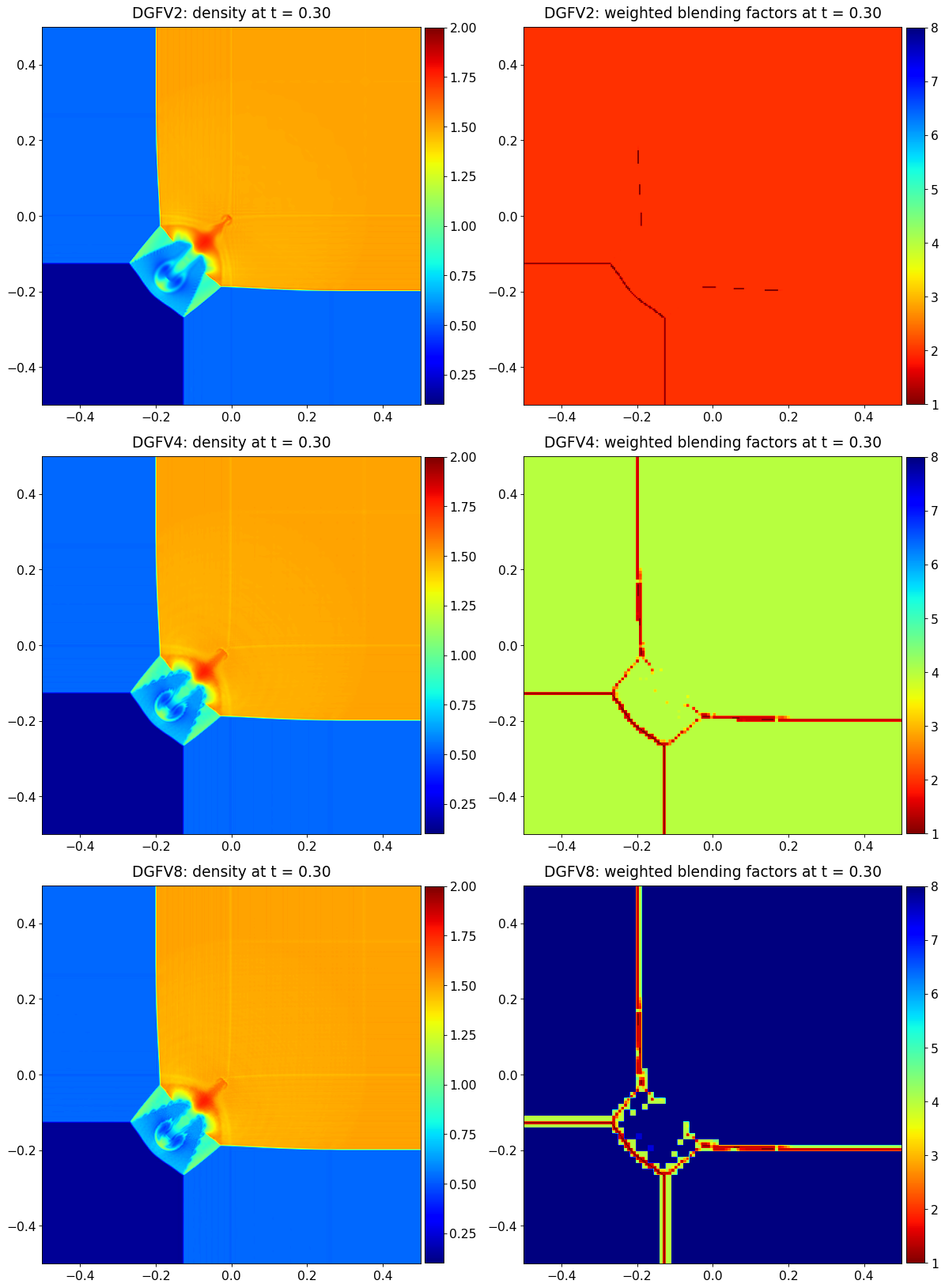}
\fi
\caption{2D Riemann problem (configuration 3, see \cite{lax1998solution})
computed with the (from top to bottom) DGFV2, DGFV4 and DGFV8 scheme. The
resolution of $512^2$ DOF is the same for all runs.  Left: density contours at
final time $T = 0.3$.  Right: Weighted blending factors as defined in
\eqref{eq:order-estimate}. Dark blue represents the full eighth order DG and
dark red the full first order FV scheme.}
\label{fig:2d-riemann-problem-dens-blend-3}
\end{figure}

\begin{figure}[H]
\centering
\ifnum \compilewithplots = 1
\includegraphics[width=0.9\textwidth]{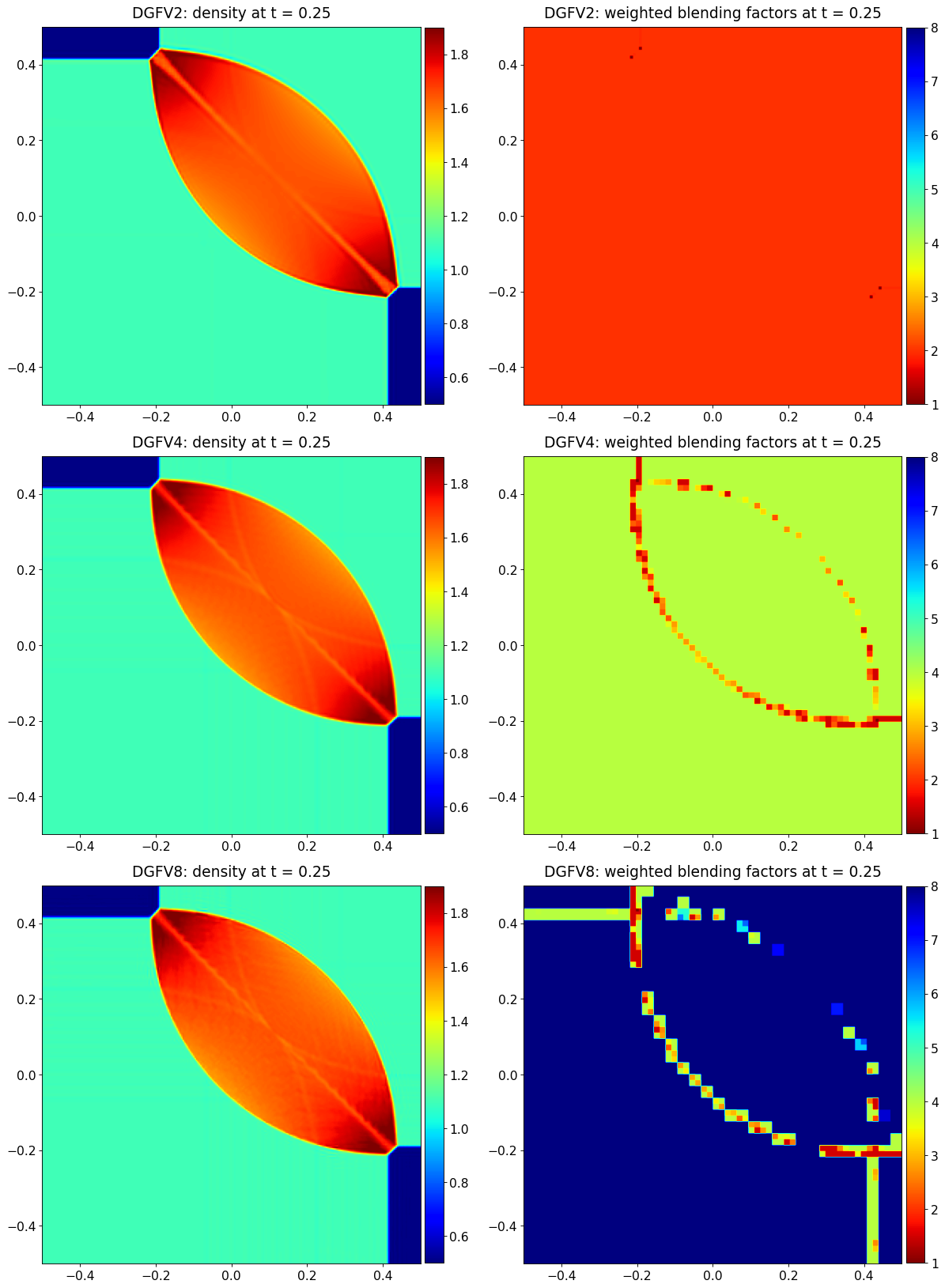}
\fi
\caption{2D Riemann problem (configuration 4, see \cite{lax1998solution})
computed with the (from top to bottom) DGFV2, DGFV4 and DGFV8 scheme. The
resolution of $512^2$ DOF is the same for all runs.  Left: density contours at
final time $T = 0.25$.  Right: Weighted blending factors as defined in
\eqref{eq:order-estimate}. Dark blue represents the full eighth order DG and
dark red the full first order FV scheme.}
\label{fig:2d-riemann-problem-dens-blend-4}
\end{figure}

\begin{figure}[H]
\centering
\ifnum \compilewithplots = 1
\includegraphics[width=0.9\textwidth]{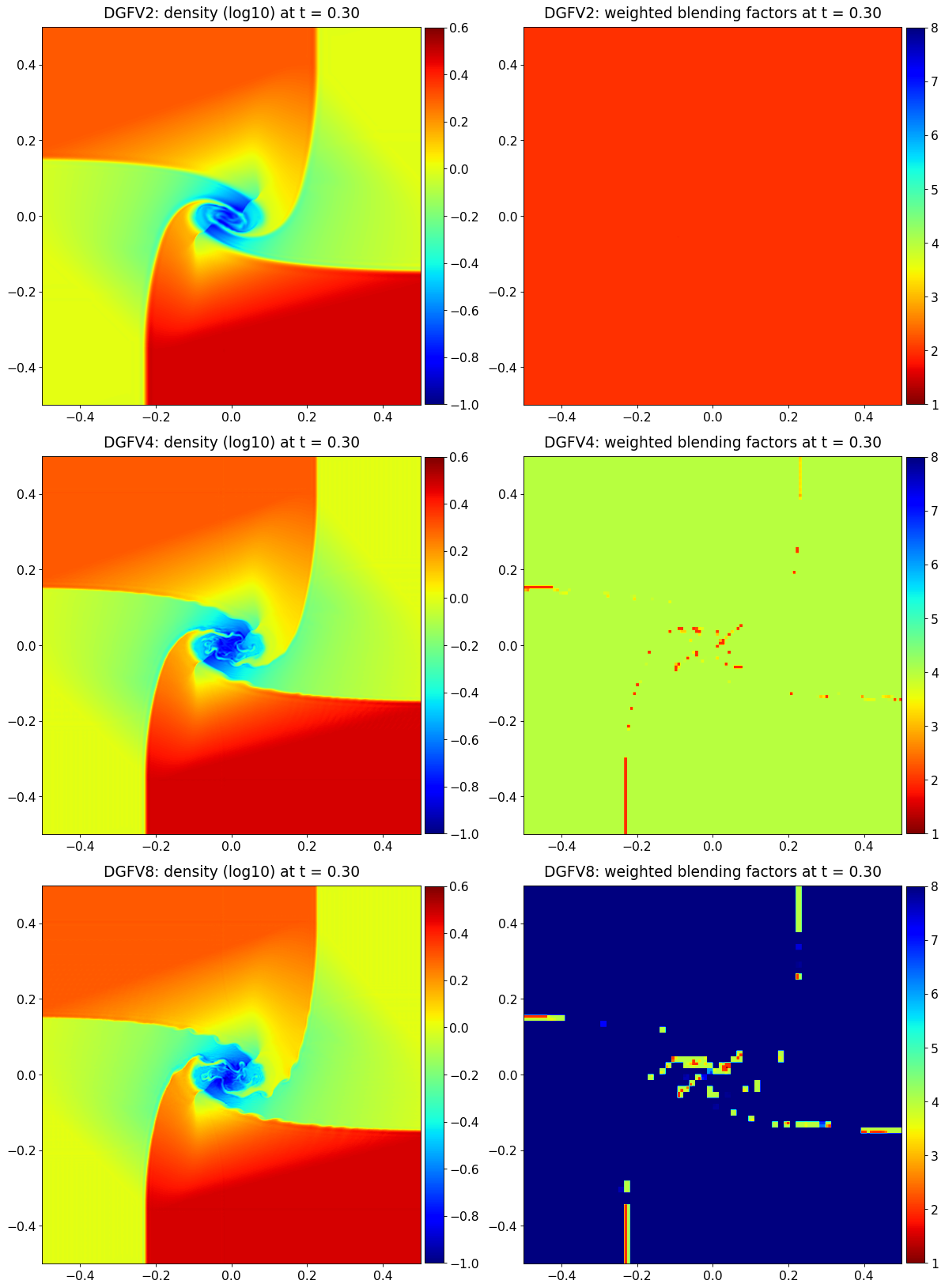}
\fi
\caption{2D Riemann problem (configuration 6, see \cite{lax1998solution})
computed with the (from top to bottom) DGFV2, DGFV4 and DGFV8 scheme. The
resolution of $512^2$ DOF is the same for all runs.  Left: density contours in
log-scale at final time $T = 0.3$.  Right: Weighted blending factors as defined
in \eqref{eq:order-estimate}. Dark blue represents the full eighth order DG and
dark red the full first order FV scheme.}
\label{fig:2d-riemann-problem-dens-blend-6}
\end{figure}
\end{changed}

\subsection{2D Sedov Blast}
\label{sec:2d-sedov-blast}
\begin{figure}[H]
\centering
\ifnum \compilewithplots = 1
\includegraphics[width=0.95\textwidth]{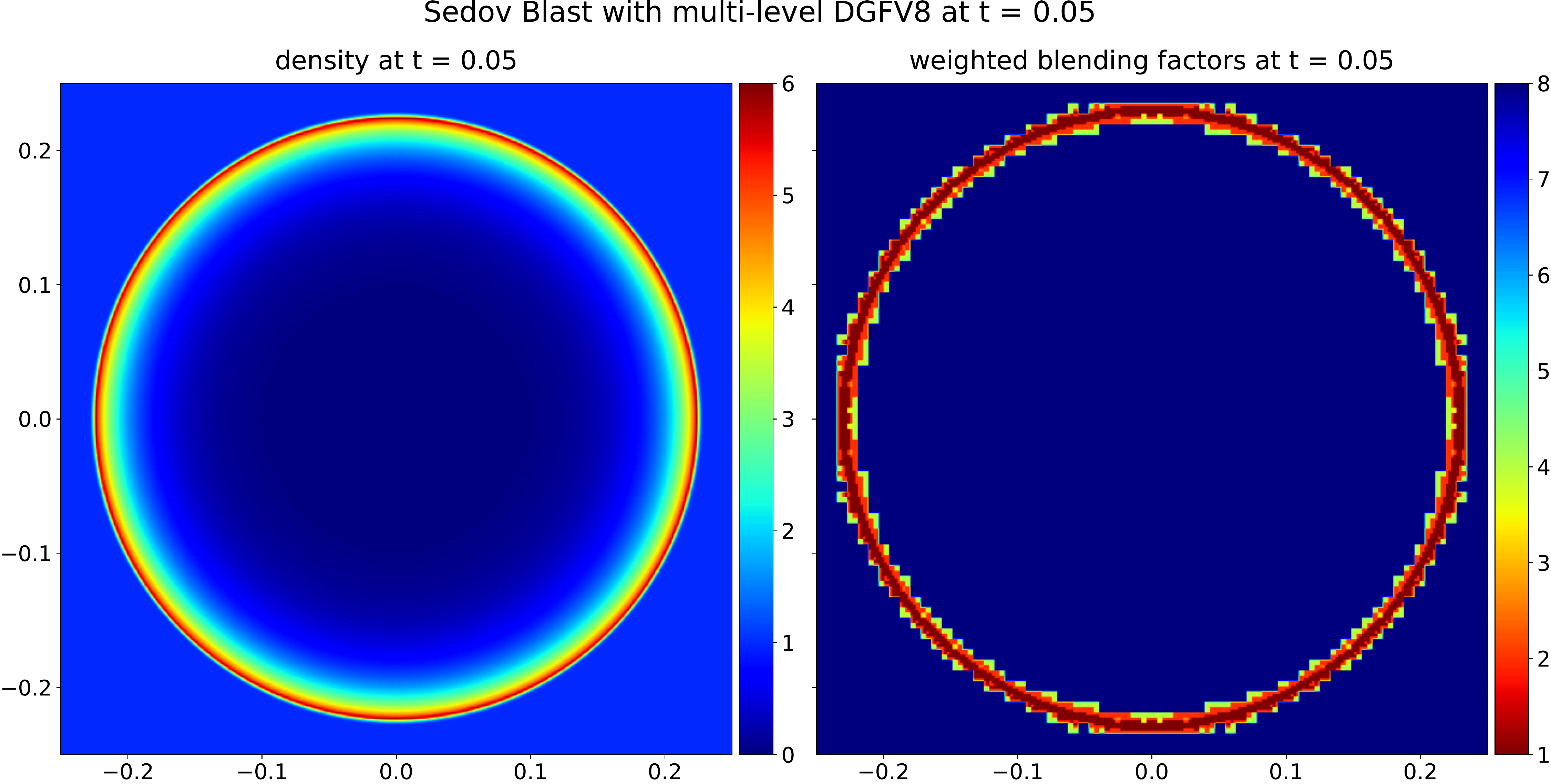}
\fi
\caption{2D Sedov Blast: Numerical solution computed with DGFV8 at final
simulation time $T = 0.05$. Left: density contours.
Right: Weighted blending factors as defined in \eqref{eq:order-estimate}. Dark
blue represents the full eighth order DG and dark red the full first order FV
scheme.}
\label{fig:2d-sedov-blast-dgfv8}
\end{figure}

\begin{figure}[H]
\centering
\ifnum \compilewithplots = 1
\includegraphics[width=0.8\textwidth]{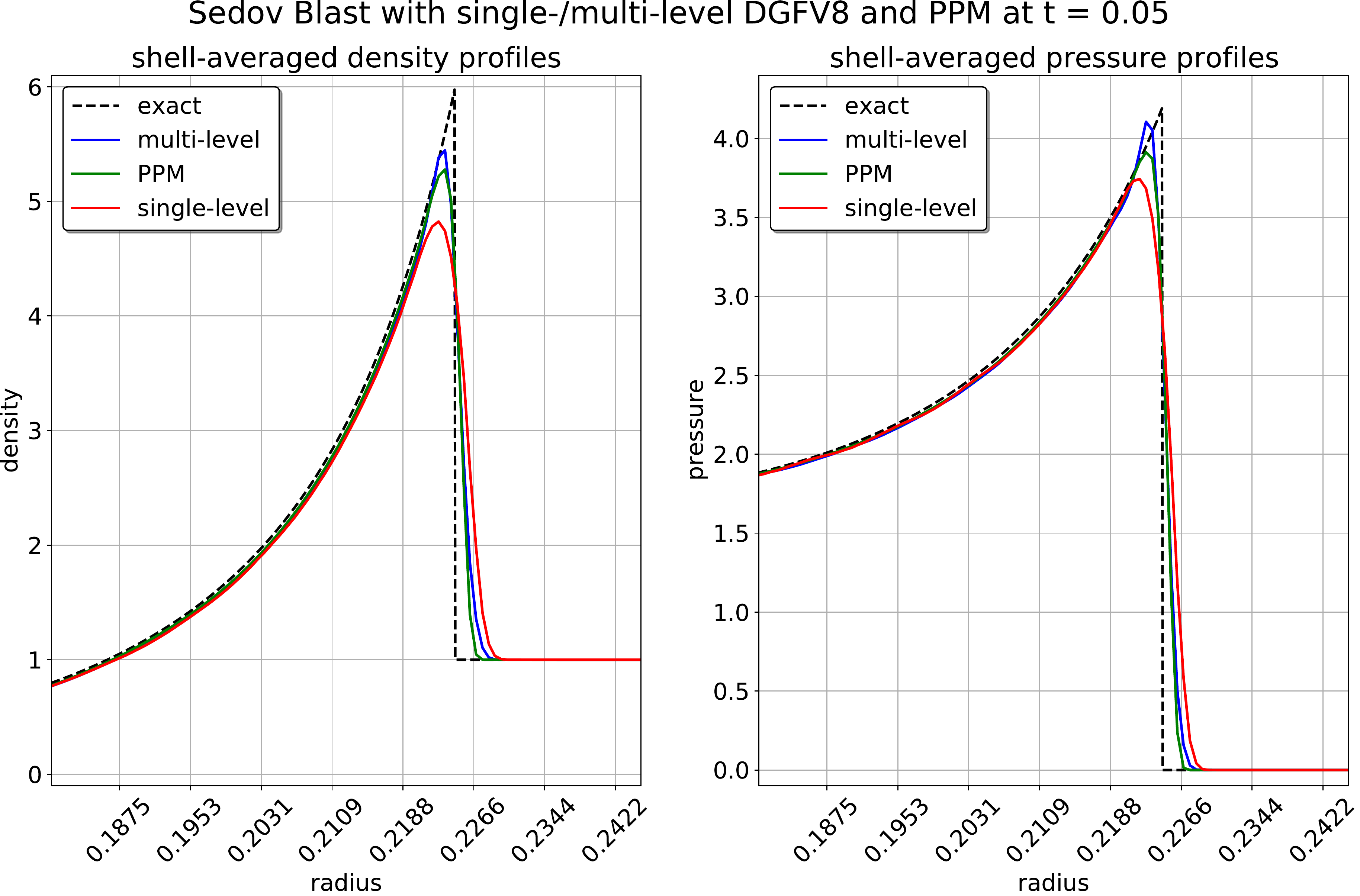}
\fi
\caption{2D Sedov Blast: shell-averaged density and pressure profiles at final
simulation time $T = 0.05$ computed with the single-/multi-level blending
scheme (DGFV8) and PPM on a FV equivalent uniform grid resolution of $512^2$
DOF. The vertical grid lines depict the element boundaries.}
\label{fig:2d-sedov-blast-density}
\end{figure}

\begin{figure}[H]
\centering
\ifnum \compilewithplots = 1
\includegraphics[width=0.8\textwidth]{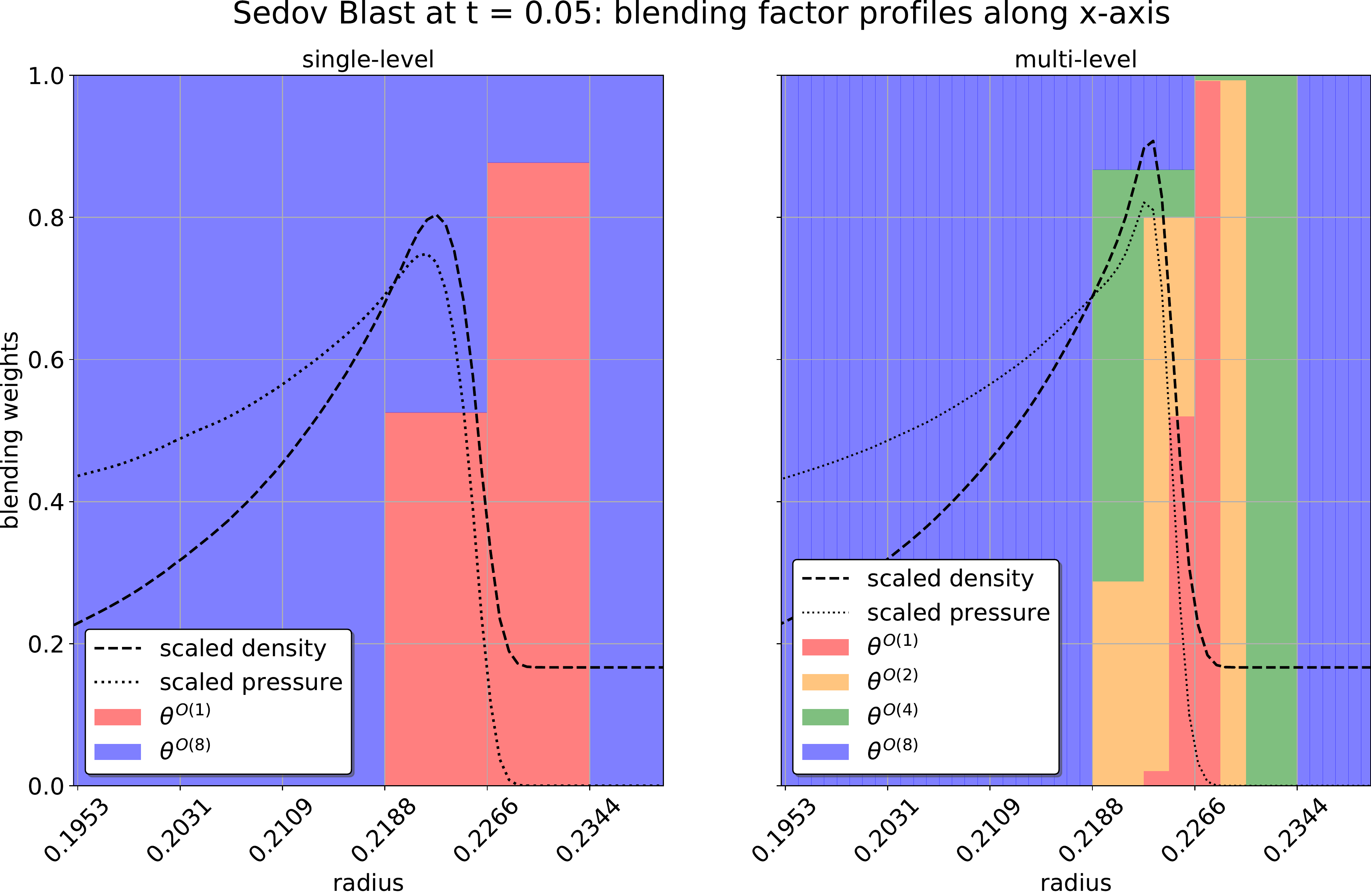}
\fi
\caption{2D Sedov Blast: We compare the blending factor profiles along the
x-axis of the single- and multi-level blending schemes at final time $T =
0.05$. The blending weights $\theta^{\bigO(n)}$ (eq.
\eqref{eq:blending-weights}) are encoded with stacked bars in their respective
colors. The vertical grid lines highlight the eighth order DG element
boundaries. For reference the scaled density ($\rho/6$) and scaled
pressure ($p/5$) are included.}
\label{fig:2d-sedov-blast-dgfv8-blending-profile}
\end{figure}

The Sedov blast problem \cite{zhang2010positivity,zhang2012positivity,dumbser2018efficient}
describes the self-similar evolution of a radially symmetrical blast wave from
an initial pressure point (delta distribution) at the center into the
surrounding, homogeneous medium. The analytical solution is given
in \cite{sedov19591,korobeinikov1991problems}. In our setup, we approximate the
initial pressure point with a smooth Gaussian distribution
\begin{equation}
\label{eq:sedov-init-energy}
\Ener_0(\vec{x}) = \frac{p_0}{\gamma-1} + \frac{E}{(2\,\pi\,\sigma^2)^{d/2}} \,
    \exp{\Big(-\frac{1}{2}\,\frac{\vec{x}^T\vec{x}}{\sigma^2}\Big)},
\end{equation}
with the spatial dimension $d = 2$, the blast energy $E = 1$ and the width
$\sigma$ such that the initial Gaussian is reasonably resolved. The surrounding
medium is initialized with $\rho_0 = 1$ and $\pres_0 = 10^{-14}$. Dimensional
analysis \cite{sedov19591} reveals that the analytical solution of the density
right at the shock front is determined by
\begin{equation*}
    \rho_{\text{shock}} = \frac{\gamma + 1}{\gamma -1}\,\rho_0.
\end{equation*}
With the adiabatic coefficient $\gamma=1.4$, we investigate how close the
numerical results match $\rho_{\text{shock}} = 6$. The Cartesian mesh has a FV
equivalent uniform grid resolution of $512^2$ DOF, i.e. when computing with the
full eighth order DG approximation space, the total number of DG elements is
$64^2$. The spatial domain is $\Omega \in [-0.25,0.25]^2$ with the initial
blast width $\sigma = 5\times 10^{-3}$. \begin{changed}We compare the accuracy of the results
obtained with the single-level and the multi-level (DGFV8) blending discretizations as well as the PPM. Fig. \ref{fig:2d-sedov-blast-density} shows the shell-averaged density
and pressure profiles at final time $T = 0.05$. Fig.
\ref{fig:2d-sedov-blast-dgfv8} presents the numerical solution over the whole
domain as computed with the multi-level DGFV8. To further illustrate the
behavior of the multi-level and single-level blending approach we also show in
Fig. \ref{fig:2d-sedov-blast-dgfv8-blending-profile} the weighted blending
factors along the x-axis. The shock front is much sharper for the multi-level
blending compared to the single-level blending.  It can be observed how the
weighted blending factor is dominated by the first order FV scheme
($\bar{\alpha}\approx 1$) directly at the shock, but transitions quickly to a
blended discretization ($1<\bar{\alpha}<8$) up to the full eighth order DG
($\bar{\alpha}\approx 8$) away from the shock, even \emph{within a single DG
element}. This behavior demonstrates the sub-element adaptivity of our novel
approach.  Again, similarly to the shock tube section \ref{sec:1d-shock-tubes}
and the 2D Riemann problem section \ref{sec:2d-riemann-problems} the results
for the 2D Sedov blast wave show the benefit of the multi-level approach, with
the numerical profiles even slightly closer to $\rho_{\text{shock}} = 6$
compared to the PPM with equal resolution of $512^2$ DOF.
\end{changed}

\section{Simulation of a Young Supernova Remnant}
\label{sec:7}
Supernova models have been analyzed and discussed for many decades and since
they unite a broad range of features such as strong shocks, instabilities and
turbulence, they resemble a good test bed for our novel shock capturing
approach in combination with AMR.

The general sequence of events of the presented supernova simulation is like
this: We start with a constant distribution of very low density resembling
interstellar media (ISM) that typically fills the space between stars. When a
star explodes by turning into a supernova it ejects its own mass at very high
speeds into the ISM preceded by a strong shock front heating up the ISM. The
ejected mass is rapidly decelerated by the swept-up ISM giving rise to a so
called reverse shock that travels backwards to the center. The interface, or
more precisely the contact discontinuity, between shocked ejecta and shocked
ISM is unstable and leads to a layer of slowly growing
Rayleigh-Taylor instabilities. This gradually expanding layer, called supernova
remnant, is of special interest since this is where astronomical observations
reveal a lot of ongoing physics and chemistry, especially driven by mixing
and turbulence. 

\newcommand{\Msun}{\text{M}_{\odot}}
\newcommand{\pc}{\text{pc}}
\newcommand{\yr}{\text{yr}}

\begin{figure}[!htbp]
\centering
\ifnum \compilewithplots = 1
    \includegraphics[width=0.6\textwidth]{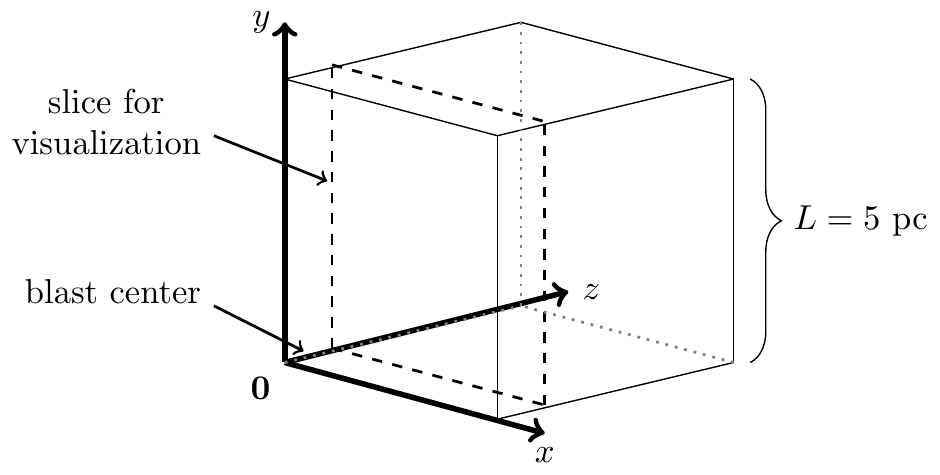}
\fi
\caption{Computational domain (cubic box) covering one octant of the supernova
model. The faces at the coordinate axes are set to reflecting walls while
the opposite sides are set to outflow.}
\label{fig:supernova-setup}
\end{figure}

We adapt the setup descriptions in
\cite{chevalier1982self,fraschetti2010simulation,ferrand2012three} where we
have the initial (internal) blast energy $E$ and the ejecta mass $M$ given in
cgs (centimeter-gram-seconds) units. It is beneficial to convert the
given units to convenient simulation units reflecting characteristic dimensions
of the physical model at hand. Table \ref{tab:cgs-to-supernova-units} lists the
conversion between cgs and simulation units and table
\ref{tab:cgs-to-supernova} lists the initial parameters used in this
simulation.
\begin{table}[!htbp]
\captionsetup{width=0.95\textwidth}
\caption{Conversion from cgs units to simulation units.}
\centering
\begin{tabular}{c|c|l}
\toprule
\multicolumn{1}{c|}{quantity} & \multicolumn{1}{c|}{cgs units} & \multicolumn{1}{c}{simulation units} \\
\midrule
mass        &    $\text{g}$        & $\Msun = 1.989\times10^{33}\,\text{g}$ \\
time        &    $\text{s}$        & $\yr = 3.154\times10^{7}\,\text{s}$ \\
length      &    $\text{cm}$       & $\pc = 3.086\times10^{18}\,\text{cm}$ \\
temperature &    $\text{K}$        & $\text{K}$ \\
\bottomrule
\end{tabular}
\label{tab:cgs-to-supernova-units}
\end{table}
\vspace{-5mm}
\begin{table}[!htbp]
\captionsetup{width=0.95\textwidth}
\caption{Hydrodynamical parameters in cgs units and simulation units.}
\centering
\begin{tabular}{l|l|l}
\toprule
\multicolumn{1}{c|}{description} & \multicolumn{1}{c|}{cgs units} & \multicolumn{1}{c}{simulation units} \\
\midrule
domain size         & $L = 1.543\times10^{19}\,\text{cm}$       & $L = 5 \, \pc$ \\
blast energy        & $E = 10^{51} \text{erg}$                  & $E = 5.2516\times10^{-5}\,\Msun \frac{\pc^2}{\yr^2}$ \\
ejecta mass         & $M = 2.7846\times10^{33}\,\text{g}$       & $M = 1.4 \, \Msun$ \\
ambient density     & $n_H = 0.13 \, \text{cm}^{-3}$            & $\rho_{a} = 2.4539\times10^{-3} \,\frac{\Msun}{\pc^3}$ \\
ambient temperature & $T = 10^4\,\text{K}$                      & $p_{a} = 2.1309\times10^{-13} \,\frac{\Msun}{\pc \, \yr^2}$ \\
\bottomrule
\end{tabular}
\label{tab:cgs-to-supernova}
\end{table}
The ambient density $\rho_a$ is related to the mono-atomic particle (hydrogen)
density $n_H$ via $\rho_a = m_u \, n_H$, where $m_u$ is $\frac{1}{12}$ of the mass
of a carbon-12 atom. The ambient pressure $p_a$ is calculated from the ideal
gas law, i.e.  $p_a = n_H\,k_B\,T$ with the Boltzmann constant $k_B$. There is
some ambiguity regarding the ambient gas temperature $T$. The literature
mentions a warm and neutral interstellar medium which is attributed to
temperatures between $6\times10^3\,\text{K}$ and $10^4\,\text{K}$. The heat capacity ratio is
fixed to $\gamma = 5/3$. The simulation time spans a period from $t_0 = 10\,\yr$
to $T = 500\,\yr$. The expansion of the forward shock (eq.
\eqref{eq:expansion-law}) is then expected to approximately reach
$R_{\text{FS}} = 5\,\pc$, which determines the size of the computational domain
$L := 5\,\pc$. Fig. \ref{fig:supernova-setup} depicts a schematic of the
simulation setup. Due to the rotational symmetry of the setup it is sufficient
to simulate just one octant of the supernova.
The following formulas have been derived in \cite{chevalier1982self}
and were adapted to the current setup, i.e. power-law indices of $(s,n) =
(0,7)$. The self-similar solution at initial time $t_0 = 10\,\yr$ within the
power law region, respectively the blast center, is given by
\begin{equation}
r(t) = t \; \sqrt{\frac{5}{3}\,\frac{E}{M}}
\quad\text{and}\quad
\dens(t) = \frac{25}{21\,\pi} \; \frac{E^2}{M} \, t^4 \, r(t)^{-7}.
\end{equation}
We initialize the density as
\begin{equation}
\label{eqn:supernova-init-density}
\dens_0(\vec{x}) =  \dens_a + \dens(t_0) \cdot
\begin{lcases}
     & \quad\;\;\; 1, &                                  \quad |\vec{x}| \leq r(t_0), \\
     & \left(\frac{|\vec{x}|}{r(t_0)}\right)^{-7}, & \quad |\vec{x}| > r(t_0),
\end{lcases}
\end{equation}
and the total energy with \eqref{eq:sedov-init-energy} where $\pres_0 =
\pres_a$, $d = 3$ and $\sigma = \frac{3}{4}\,r(t_0)$.  The initial momentum is
$(\rho\,\velv)_0 = \vec{0}$. Since we are only interested in the evolution of
the instability layer of the supernova we apply the following rules for mesh
refinement and coarsening. The expansion radius of the forward shock over time
is given by
\begin{equation}
\label{eq:expansion-law}
R_{\text{FS}}(t) = 1.06 \left(\frac{E^2}{M\,\rho_{a}}\right)^{1/7} \, t^{4/7}.
\end{equation}
This allows us to assign an adaptive, high resolution shell of maximal refined
elements following the remnant as it expands into the computational domain. The
inner and outer radii of the shell are estimated as
\begin{equation}
\label{eq:adaptive-shell}
R_{\text{inner}}(t) = 0.7\,R_{\text{FS}}(t) \quad\text{and}\quad R_{\text{outer}}(t) = 1.15\,R_{\text{FS}}(t),
\end{equation}
which have been found adequate via numerical experimentation. Moreover, up to
$t = 200\,\yr$ we enforce $R_{\text{inner}} = 0$, which ensures that the first
phase of the explosion is well resolved in any case. The refinement levels
range from 2 to 6, which translates to a FV equivalent resolution from
$2^2\cdot 8 = 16$ up to $2^6\cdot 8 = 512$ cells in each spatial direction. 

\begin{changed}
We perform three simulations with the 3D multi-level blending schemes DGFV2,
DGFV4 and DGFV8, analog to the 2D Riemann problems discussed in section
\ref{sec:2d-riemann-problems}. Fig. \ref{fig:3D-supernova-all-orders} shows the
density slice (left column) sketched in Fig. \ref{fig:supernova-setup} at the
final simulation time of $T = 500\,\yr$, while Fig.
\ref{fig:3D-supernova-density-contour} shows the corresponding 3D density
contour rendering of the DGFV8 solution.
The shock front partially left the domain, which is not considered a problem
since the region of interest, namely the instability layer, is still completely
covered. The areas of no interest, i.e. the outer region as well as the center,
are only coarsely resolved by the AMR scheme. Clearly, an increase in order
leads to a much more detailed remnant structure emphasizing the advantage of higher order
schemes in resolving small scale turbulence driven by Rayleigh-Taylor instability, see \cite{chevalier1982self}. The weighted
blending factor is shown in the right column of Fig.
\ref{fig:3D-supernova-all-orders}. The band of highly refined elements is
clearly visible following the remnant as intended. Two distinctive lines of
blending activity trace the front and reverse shocks. The plots also show a
strong qualitative difference of the amount of scales in the instability layer
for DGFV2 and DGFV4. Clearly, the DGFV4 result features more scales and finer
structures as the DGFV2 simulation with the same FV equivalent grid resolution.
The difference in scales between DGFV4 and DGFV8 is less pronounced which is
probably related to the extensive blending activity with O4 inside the
instability layer of the DGFV8 simulation. In this turbulent part we suspect
that the standard collocation eighth order DG scheme might face aliasing
instability issues as described in the introduction. There are techniques to
reduce aliasing issues available, such as filtering, consistent integration and
split-forms. However, this aspect is not the main topic of the paper and it is
interesting to see, that the multi-level blending automatically adjusts to cope
with these issues as well. 
\begin{figure}[H]
\centering
\ifnum \compilewithplots = 1
\includegraphics[width=0.9\textwidth]{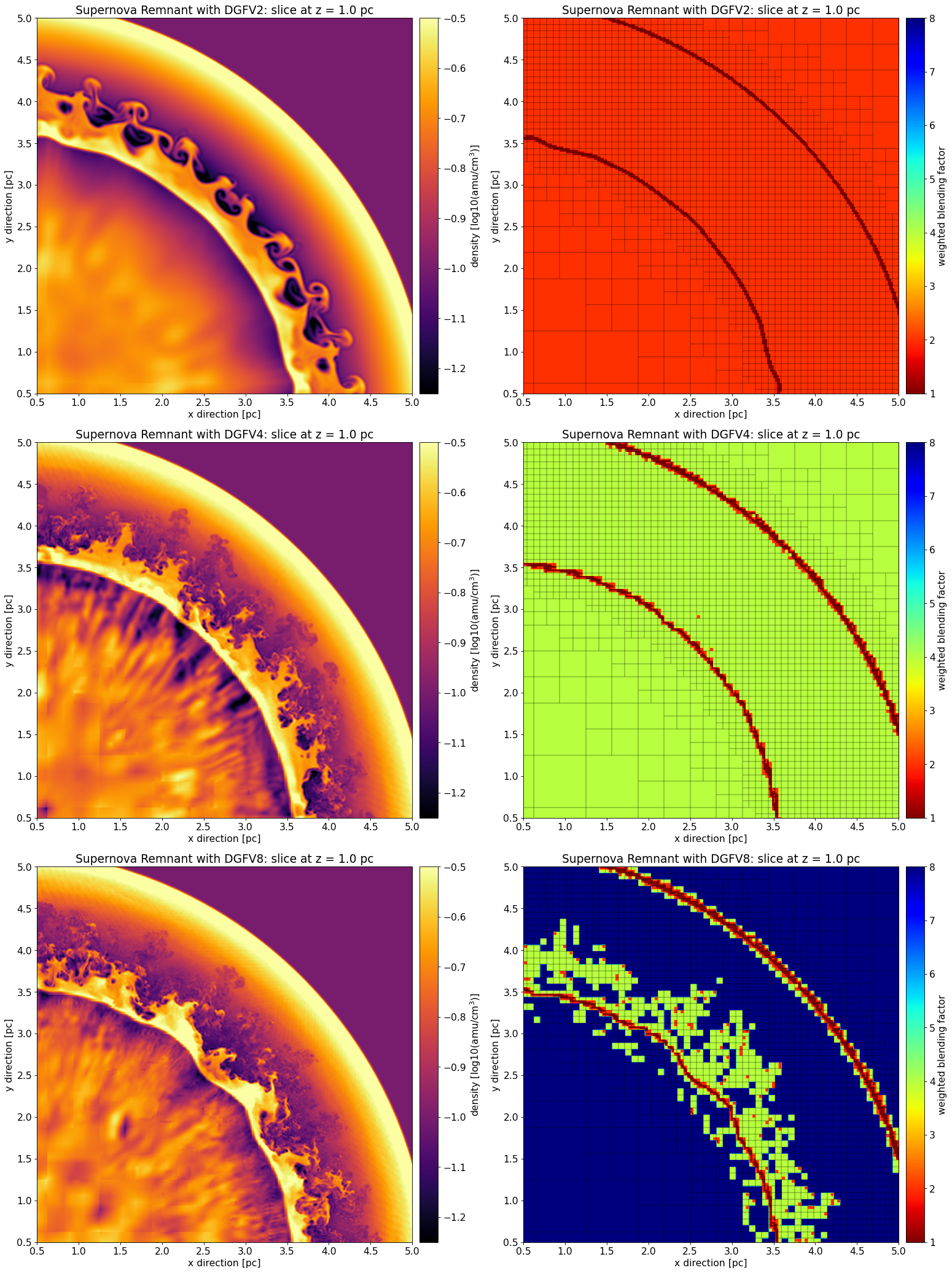}
\fi
\caption{\textit{left column:} 2D density slice (see figure
\ref{fig:supernova-setup}) showing the instability region respectively remnant
of the supernova at $T = 500\,\yr$ simulated with the (from top to bottom)
DGFV2, DGFV4 and DGFV8 scheme. \textit{right column:} 2D slice of the weighted
blending factors of the respective blending schemes at $T = 500\,\yr$. The black
lines correspond to the element boundaries of the Cartesian non-conforming
mesh.}
\label{fig:3D-supernova-all-orders}
\end{figure}
\begin{figure}[H]
\centering
\ifnum \compilewithplots = 1
\includegraphics[width=0.7\textwidth]{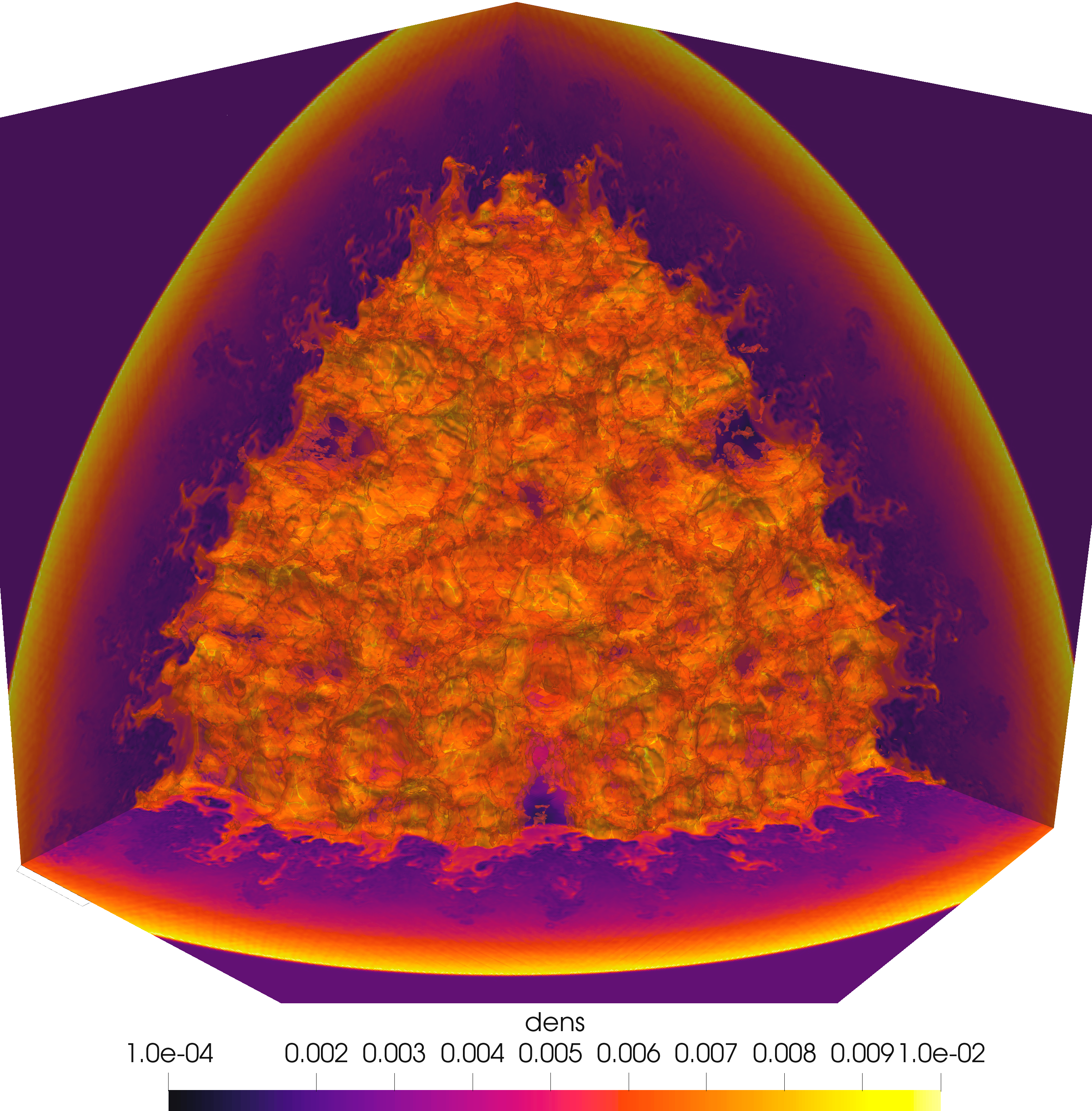}
\fi
\caption{3D density contour rendering showing the instability region
respectively remnant of the supernova at final simulation time $T = 500\,\yr$
simulated with the DGFV8 scheme.}
\label{fig:3D-supernova-density-contour}
\end{figure}
\end{changed}

\section{Conclusion}

In this work, we introduce an adaptive sub-element based shock capturing
approach for DG methods. We interpret an element first as a collection of
piecewise constant data. This set of piecewise constant data can be interpreted
and reconstructed in a variety of ways. We focus on piecewise polynomial
approximation spaces, starting from a pure piecewise constant FV type
interpretation (no reconstruction) up to the fully maximum order polynomial
reconstruction and every combination in-between, e.g. piecewise linear and
piecewise cubic reconstructions. In a second step, we link the data
interpretation to a corresponding (high order) DG discretization.  Thus, we get
a hierarchy of discretizations acting on the same set of data. 

The idea is then to adaptively blend this hierarchy of different
discretization, where the low order variants are chosen close to
discontinuities and the high order variants in smooth or turbulent parts of the
simulation. Instead of having an element based troubled cell indicator
approach, we use the different data interpretations again to compute
sub-element localized indicators, which allows for a sub-element adaptive
blending of the discretizations. When the blending can change throughout the
element, special care is necessary to preserve exact conservation of the
resulting multi-level blended discretization. We achieve exact conservation, by
introducing unique, blended reconstruction states at subcell and sub-element
interfaces. 

In our prototype implementation, we further demonstrate that a natural
combination of this sub-element adaptive approach is with adaptive mesh
refinement.  We base the AMR implementation on the \pforest Octree library,
which allows for a straight forward parallelization of the whole computational
framework.   

We show standard numerical test cases to validate the new approach and a
simplified model of a supernova remnant to highlight its high accuracy for
challenging test cases with strong shocks and turbulence like structures.

\section*{Acknowledgements}
JM acknowledges funding through the Klaus-Tschira Stiftung via the project
"DG$^2$RAV". GG thanks the Klaus-Tschira Stiftung and the European Research
Council for funding through the ERC Starting Grant “An Exascale aware and
Un-crashable Space-Time-Adaptive Discontinuous Spectral Element Solver for
Non-Linear Conservation Laws” (EXTREME, project no. 71448).  SW thanks the
Klaus-Tschira Stiftung, and acknowledges the Deutsche Forschungsgemeinschaft
(DFG) for funding the sub-project C5 in the SFB956 and the European Research
Council for funding the ERC Starting Grant “The radiative interstellar medium”
(RADFEEDBACK, project no. 679852). This work was performed on the Cologne High
Efficiency Operating Platform for Sciences (CHEOPS) at the Regionales
Rechenzentrum K\"oln (RRZK). Furthermore, we thank Dr. Andrés Rueda-Ramírez for
his help with the generation of the 3D supernova plot in Fig
\ref{fig:3D-supernova-density-contour}.

\section{APPENDIX}

\subsection{Four levels of data interpretation of $8$ mean values}
\label{sec:appendix-visu-multilevel-blending}
\begin{figure}[H]
\centering
\ifnum \compilewithplots = 1
    \includegraphics[width=0.7\textwidth]{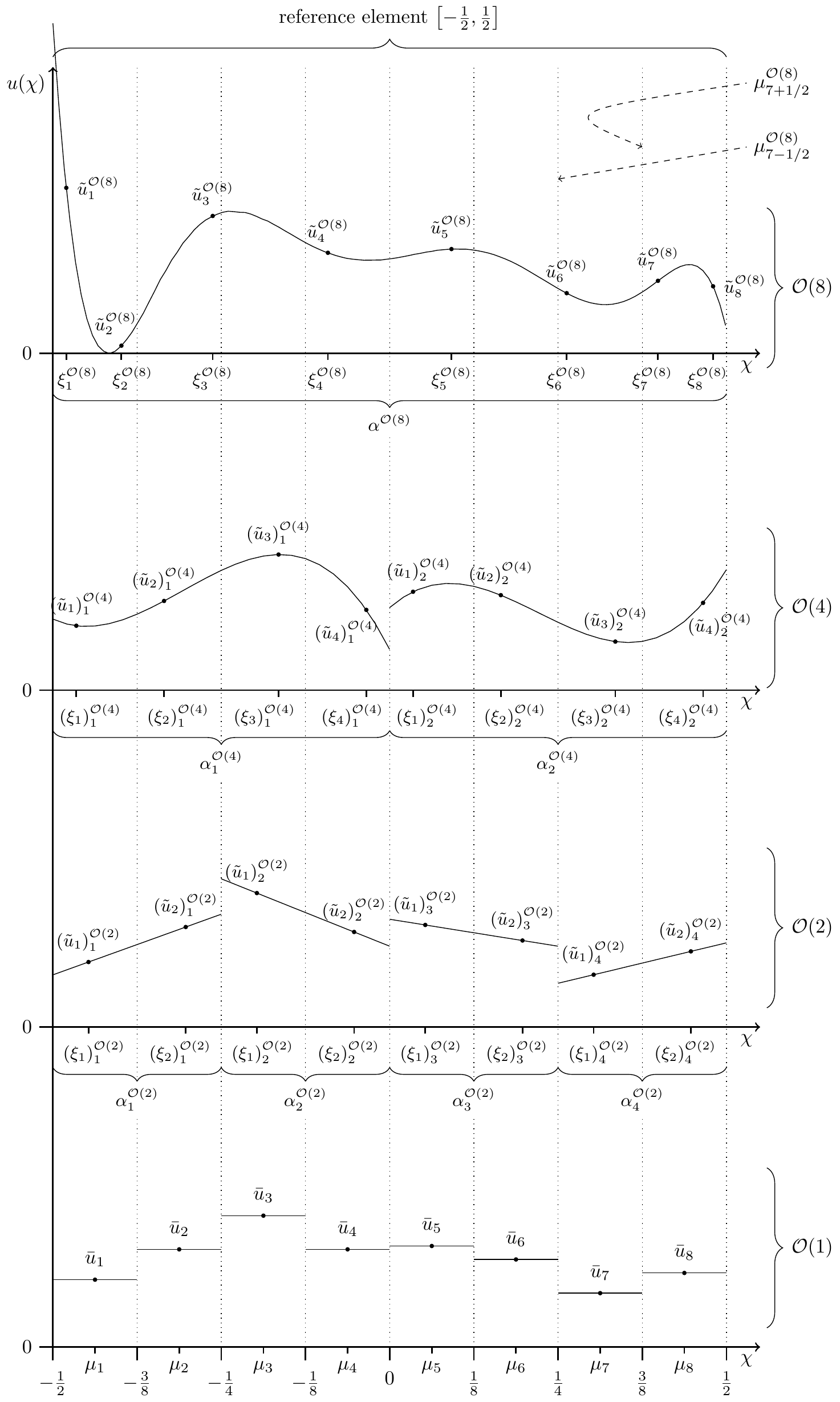}
\fi
\caption{Schematic of eight ($N = 8$) mean values and its four levels of data interpretation.
The values $\mu_{7\pm1/2}^{\bigO(8)}$ indicate the location of the interpolation nodes
in \eqref{eq:subcell-fv-interpolation-operator} of the FV subcell interpolation operator $\mat{I}^{\pm\bigO(n)}$,
$n = 2,4,8$.
}
\label{fig:multilevel-blending}
\end{figure}

\subsection{Boundary Evaluation Operators for the Multi-level Blending Scheme \eqref{eq:multilevel-blending-scheme}.}
\label{sec:appendix-boundary-eval-op}
\begin{equation*}
\left[\idmat_4 \otimes \dgfvSurfOp^{\bigO(2)}\right]^{(0)} =
\left(
\begin{array}{*{9}{@{}C{\mycolwd}@{}}}
\hat{B}^{\bigO(2)}_{1\,1} & 0 & \hat{B}^{\bigO(2)}_{1\,3} & 0 & 0 & 0 & 0 & 0 & 0 \\[2mm]
\hat{B}^{\bigO(2)}_{2\,1} & 0 & \hat{B}^{\bigO(2)}_{2\,3} & 0 & 0 & 0 & 0 & 0 & 0 \\[2mm]
0 & 0 & \hat{B}^{\bigO(2)}_{1\,1} & 0 & \hat{B}^{\bigO(2)}_{1\,3} & 0 & 0 & 0 & 0 \\[2mm]
0 & 0 & \hat{B}^{\bigO(2)}_{2\,1} & 0 & \hat{B}^{\bigO(2)}_{2\,3} & 0 & 0 & 0 & 0 \\[2mm]
0 & 0 & 0 & 0 & \hat{B}^{\bigO(2)}_{1\,1} & 0 & \hat{B}^{\bigO(2)}_{1\,3} & 0 & 0 \\[2mm]
0 & 0 & 0 & 0 & \hat{B}^{\bigO(2)}_{2\,1} & 0 & \hat{B}^{\bigO(2)}_{2\,3} & 0 & 0 \\[2mm]
0 & 0 & 0 & 0 & 0 & 0 & \hat{B}^{\bigO(2)}_{1\,1} & 0 & \hat{B}^{\bigO(2)}_{1\,3} \\[2mm]
0 & 0 & 0 & 0 & 0 & 0 & \hat{B}^{\bigO(2)}_{2\,1} & 0 & \hat{B}^{\bigO(2)}_{2\,3}
\end{array}\right) \; \in \reals^{8 \times 9}
\end{equation*}
\begin{equation*}
\left[\idmat_2 \otimes \dgfvSurfOp^{\bigO(4)}\right]^{(0)} =
\left(
\begin{array}{*{9}{@{}C{\mycolwd}@{}}}
\hat{B}^{\bigO(4)}_{1\,1} & 0 & 0 & 0 & \hat{B}^{\bigO(4)}_{1\,5} & 0 & 0 & 0 & 0 \\[2mm]
\hat{B}^{\bigO(4)}_{2\,1} & 0 & 0 & 0 & \hat{B}^{\bigO(4)}_{2\,5} & 0 & 0 & 0 & 0 \\[2mm]
\hat{B}^{\bigO(4)}_{3\,1} & 0 & 0 & 0 & \hat{B}^{\bigO(4)}_{3\,5} & 0 & 0 & 0 & 0 \\[2mm]
\hat{B}^{\bigO(4)}_{4\,1} & 0 & 0 & 0 & \hat{B}^{\bigO(4)}_{4\,5} & 0 & 0 & 0 & 0 \\[2mm]
0 & 0 & 0 & 0 & \hat{B}^{\bigO(4)}_{1\,1} & 0 & 0 & 0 & \hat{B}^{\bigO(4)}_{1\,5} \\[2mm]
0 & 0 & 0 & 0 & \hat{B}^{\bigO(4)}_{2\,1} & 0 & 0 & 0 & \hat{B}^{\bigO(4)}_{2\,5} \\[2mm]
0 & 0 & 0 & 0 & \hat{B}^{\bigO(4)}_{3\,1} & 0 & 0 & 0 & \hat{B}^{\bigO(4)}_{3\,5} \\[2mm]
0 & 0 & 0 & 0 & \hat{B}^{\bigO(4)}_{4\,1} & 0 & 0 & 0 & \hat{B}^{\bigO(4)}_{4\,5}
\end{array}\right) \; \in \reals^{8 \times 9}
\end{equation*}
\begin{equation*}
\left[\idmat_1 \otimes \dgfvSurfOp^{\bigO(8)}\right]^{(0)} =
\left(
\begin{array}{*{9}{@{}C{\mycolwd}@{}}}
\hat{B}^{\bigO(8)}_{1\,1} & 0 & 0 & 0 & 0 & 0 & 0 & 0 & \hat{B}^{\bigO(8)}_{1\,9} \\[2mm]
\hat{B}^{\bigO(8)}_{2\,1} & 0 & 0 & 0 & 0 & 0 & 0 & 0 & \hat{B}^{\bigO(8)}_{2\,9} \\[2mm]
\hat{B}^{\bigO(8)}_{3\,1} & 0 & 0 & 0 & 0 & 0 & 0 & 0 & \hat{B}^{\bigO(8)}_{3\,9} \\[2mm]
\hat{B}^{\bigO(8)}_{4\,1} & 0 & 0 & 0 & 0 & 0 & 0 & 0 & \hat{B}^{\bigO(8)}_{4\,9} \\[2mm]
\hat{B}^{\bigO(8)}_{5\,1} & 0 & 0 & 0 & 0 & 0 & 0 & 0 & \hat{B}^{\bigO(8)}_{5\,9} \\[2mm]
\hat{B}^{\bigO(8)}_{6\,1} & 0 & 0 & 0 & 0 & 0 & 0 & 0 & \hat{B}^{\bigO(8)}_{6\,9} \\[2mm]
\hat{B}^{\bigO(8)}_{7\,1} & 0 & 0 & 0 & 0 & 0 & 0 & 0 & \hat{B}^{\bigO(8)}_{7\,9} \\[2mm]
\hat{B}^{\bigO(8)}_{8\,1} & 0 & 0 & 0 & 0 & 0 & 0 & 0 & \hat{B}^{\bigO(8)}_{8\,9}
\end{array}\right) \; \in \reals^{8 \times 9}
\end{equation*}

\subsection{Extension to 3D Cartesian Grids}
\label{sec:appendix-3D-extension}
The extension to higher spatial dimensions with Cartesian conforming grids is
relatively straight forward via a tensor-product strategy. In this section we
want to present the most important building blocks for the 3D single-level
blending scheme. The steps of the 3D multi-level blending are analogous as in
section \ref{sec:multilevel-blending} and not detailed any further. Consider
the general 3D conservation law
\begin{equation}
\label{eq:conservation-law-2d}
\partial_t u + \partial_x f(u) + \partial_y g(u) + \partial_z h(u) = 0.
\end{equation}
We now have an element $q$ with $N^3$ mean values as available data
\begin{equation*}
    \block{\bar{u}}_q \in \reals^{N \times N \times N},
\end{equation*}
and the element sizes $\Delta x_q$, $\Delta y_q$ and $\Delta z_q$.
With the tensor-product of Lagrange polynomials $\ell$ on Legendre-Gauss nodes
$\xi_i$, $i = 1,\ldots,N$,
\begin{equation}
\label{eq:tensor-product-ansatz}
    u(t;\xi,\eta,\zeta) \approx \sum_{i,j,k=1}^{N} \tilde{u}_{ijk}(t) \, \ell_i(\xi)\,\ell_j(\eta)\,\ell_k(\zeta),
\end{equation}
and the reconstruction operator $\recoMat$, which we introduced in section
\ref{sec:reconstruction}, we get the nodal coefficients as
\begin{equation}
\label{eq:3D-nodal-reconstruction}
    \block{\tilde{u}}_q = \block{\bar{u}}_q \bigtimes_{d=x}^z \recoMat = \block{\bar{u}}_q \times_x \recoMat \times_y \recoMat \times_z \recoMat \quad \in \reals^{N \times N \times N}.
\end{equation}
Here we make use of the \emph{n-mode} product notation \cite{kolda2009tensor}.
A definition is given in appendix \ref{sec:appendix-n-node-product}. In general
all vector and matrix operations presented so far can be directly expanded to
3D via the n-mode product. The blending scheme \eqref{eq:dgfv-scheme-1d} in 3D
reads as
\begin{align}
\label{eq:dgfv-scheme-3d}
\block{\partial_t \bar{u}}_q =
        - & \, (1-\alpha_q)\,\frac{N}{\Delta x_q} \, \block{\bar{f}}_q^* \times_x  \fvdiffmat + \, \alpha_q \,\frac{1}{\Delta x_q}\,\Big\{\block{\tilde{f}}_q \times_x \dgfvVoluOp \times_y \projMat \times_z \projMat - \block{\bar{f}}_q^* \times_x \dgfvSurfOp \Big\}\nonumber \\[2mm]
        - & \, (1-\alpha_q)\,\frac{N}{\Delta y_q} \, \block{\bar{g}}_q^* \times_y  \fvdiffmat + \, \alpha_q \,\frac{1}{\Delta y_q}\,\Big\{\block{\tilde{g}}_q \times_y \dgfvVoluOp \times_z \projMat \times_x \projMat - \block{\bar{g}}_q^* \times_y \dgfvSurfOp \Big\}\\[2mm]
        - & \, (1-\alpha_q)\,\frac{N}{\Delta z_q} \, \block{\bar{h}}_q^* \times_z  \fvdiffmat + \, \alpha_q \,\frac{1}{\Delta z_q}\,\Big\{\block{\tilde{h}}_q \times_z \dgfvVoluOp \times_x \projMat \times_y \projMat - \block{\bar{h}}_q^* \times_z \dgfvSurfOp \Big\}.\nonumber
\end{align}
The DG volume fluxes are calculated from the node values as
\begin{equation*}
    \left(\block{\tilde{f}}_q\right)_{ijk} = f\left(\left(\block{\tilde{u}}_q\right)_{ijk}\right), \;
    \left(\block{\tilde{g}}_q\right)_{ijk} = g\left(\left(\block{\tilde{u}}_q\right)_{ijk}\right) \; \text{and} \;
    \left(\block{\tilde{h}}_q\right)_{ijk} = h\left(\left(\block{\tilde{u}}_q\right)_{ijk}\right), \; i,j,k = 1,\ldots,N.
\end{equation*}

\noindent\textbf{Remark:} If necessary, the $\beta$-reconstruction
\eqref{eq:positivity-limiter} ensures permissible states for the reconstructed
polynomial, i.e.
\begin{equation}
\label{eq:3D-nodal-beta-reconstruction}
    \block{\tilde{u}}_q^{(\beta)} = \block{\bar{u}}_q \bigtimes_{d=x}^z \recoMat^{(\beta)}.
\end{equation}
Again, for $\block{\bar{f}}_q^* \in \reals^{(N+1) \times N \times N}$,
$\block{\bar{g}}_q^* \in \reals^{N \times (N+1) \times N}$ and
$\block{\bar{h}}_q^* \in \reals^{N \times N \times (N+1)}$, we need to compute
common surface fluxes in order to preserve conservation. Therefore, we define
two selection operators,
\begin{equation*}
    \vec{s}^- = (1,0,\ldots,0,0)^T \in \reals^{N} \quad \text{and} \quad \vec{s}^+ = (0,0,\ldots,0,1)^T \in \reals^{N},
\end{equation*}
mapping the outermost mean values of $\block{\bar{u}}_q$ to the respective
interface.
\begin{changed}
Then the 3D analogue of the prolongation procedure
\eqref{eq:inner-surface-prolongation} along direction $d$ reads as
\begin{equation}
\label{eq:3D-prolongation}
    \bar{\face{u}}_{q\pm\frac{1}{2}}^{\pm d} = 
        (1-\alpha_q) \;\;\underbrace{\block{\bar{u}}_q \times_{d} \vec{s}^{\pm}}_{\mathclap{\substack{\text{face centered}\\\text{mean values}}}}
        \quad + \quad 
        \alpha_q \;\; \underbrace{\projMat \overbrace{\left(\block{\tilde{u}}_q \times_{d} \vec{b}^{\pm}\right)}^{\substack{\text{nodal boundary}\\\text{values}}}\projMat^T}_{\mathclap{\substack{\text{face centered}\\\text{mean values}}}} \;\:\;\; \in \reals^{N \times N}.
\end{equation}
In Fig. \ref{fig:conforming-interface} the two kinds of boundary values are
illustrated.
\begin{figure}[H]
\centering
\ifnum \compilewithplots = 1
    \includegraphics[width=0.7\textwidth]{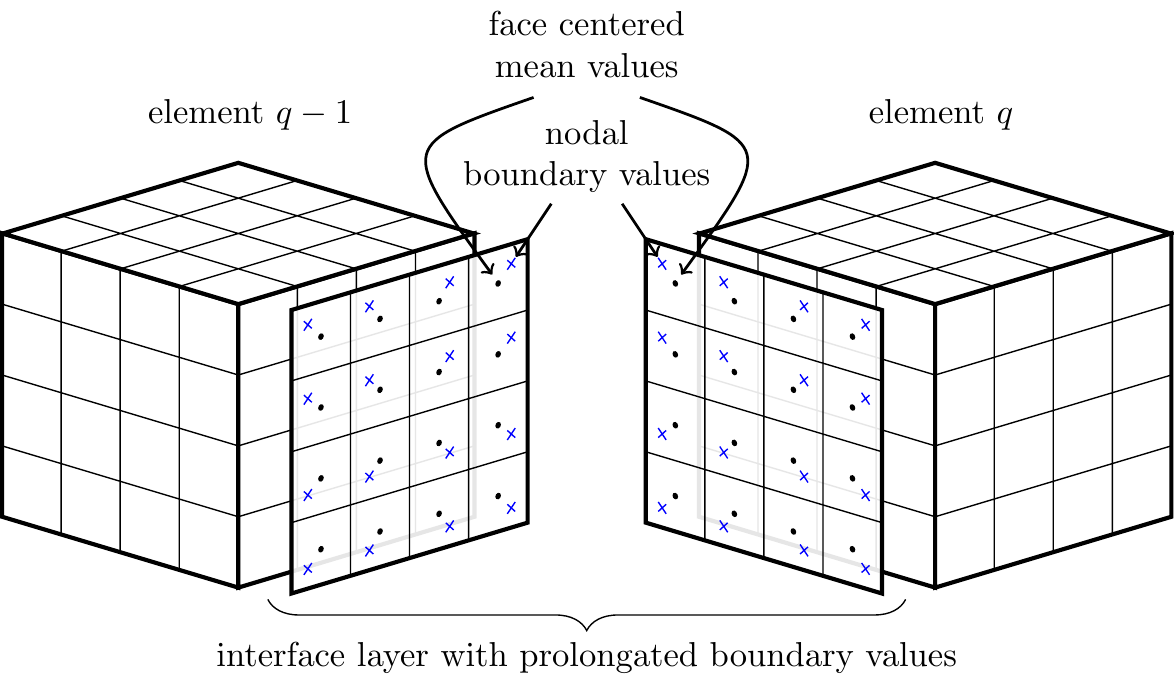}
\fi
\caption{Schematic of conforming interface $q-\frac{1}{2}$ with 1:1 adjacent
elements $q-1$ and $q$ for $N = 4$. The two kinds of boundary values occurring
in \eqref{eq:3D-prolongation} are illustrated.}
\label{fig:conforming-interface}
\end{figure}
Note that only the face centered mean values
$\bar{\face{u}}_{q\pm\frac{1}{2}}^{\pm d}$ together with the blending factor $\alpha_q$ are
supposed to be exchanged between processes in case of distributed
computing. \\

Next, we reconstruct again a nodal representation from
$\bar{\face{u}}_{q\pm\frac{1}{2}}^{\pm d}$
\begin{equation}
\label{eq:3D-boundary-value-recon}
    \tilde{\face{u}}_{q\pm\frac{1}{2}}^{\pm d} = 
        \underbrace{\recoMat\;\bar{\face{u}}_{q\pm\frac{1}{2}}^{\pm d}\;\recoMat^T}_{\mathclap{\substack{\text{nodal boundary}\\\text{values}}}} \;\:\;\; \in \reals^{N \times N}.
\end{equation}
From the two representations of boundary values $\bar{\face{u}}_{q\pm\frac{1}{2}}^{\pm x}$ and
$\tilde{\face{u}}_{q\pm\frac{1}{2}}^{\pm x}$ we calculate two candidate fluxes at the interface
$q-\frac{1}{2}$. In x-direction it reads as
\begin{align}
\label{eq:3D-interface-fluxes}
      \bar{\face{f}}_{q-\frac{1}{2}}^{*'} = f^*\left(  \bar{\face{u}}_{q-\frac{1}{2}}^{+x},  \bar{\face{u}}_{q-\frac{1}{2}}^{-x}\right) \;\:\;\; \in \reals^{N \times N}, \\
    \tilde{\face{f}}_{q-\frac{1}{2}}^{*}  = f^*\left(\tilde{\face{u}}_{q-\frac{1}{2}}^{+x},\tilde{\face{u}}_{q-\frac{1}{2}}^{-x}\right) \;\:\;\; \in \reals^{N \times N}.
\end{align}
Next, we determine a common
surface flux $\bar{\face{f}}_{q-\frac{1}{2}}^{*}$ by blending the two candidate
interface fluxes
\begin{equation}
\label{eq:3D-blend-interface-fluxes}
    \bar{\face{f}}_{q-\frac{1}{2}}^{*} = 
        (1-\alpha_{q-\frac{1}{2}}) \;\; \bar{\face{f}}_{q-\frac{1}{2}}^{*'}
        \;\;+ \;\;
        \alpha_{q-\frac{1}{2}} \;\; \projMat \; \tilde{\face{f}}_{q-\frac{1}{2}}^{*} \; \projMat^T \;\:\;\; \in \reals^{N \times N},
\end{equation}
where $\alpha_{q-\frac{1}{2}} = \frac{\alpha_{q-1}+\alpha_q}{2}$.
The final step is to calculate the inner interface fluxes 
analogous to \eqref{eq:common-surface-flux}
\begin{equation*}
    \left(\block{\bar{f}}_{q}^*\right)_{ijk} = f^*(\bar{u}_{i-1jk},\bar{u}_{ijk}), \quad i = 2,\ldots,N, \; j,k = 1,\ldots,N
\end{equation*}
and insert the common surface fluxes $\bar{\face{f}}_{q\pm\frac{1}{2}}^{*}$
\begin{align*}
    \left(\block{\bar{f}}_q^*\right)_{ijk} &= \left(\bar{\face{f}}_{q-\frac{1}{2}}^{*}\right)_{jk}, \quad i = 1, \; j,k = 1,\ldots,N \\
    \left(\block{\bar{f}}_q^*\right)_{ijk} &= \left(\bar{\face{f}}_{q+\frac{1}{2}}^{*}\right)_{jk}, \quad i = N+1, \; j,k = 1,\ldots,N.
\end{align*}
The y- and z-direction are treated analogously. The computation of the blending
parameter $\alpha$ for 3D follows the same steps as in section
\ref{sec:blending-parameter} where the integrals \eqref{eq:L1-norm-slope} are
rewritten in 3D form. The presented parameters $\tau_a$ and $\tau_s$ stay the
same. This concludes the 3D blending scheme with conforming interfaces.\\

We extend the 3D single-blending scheme to non-conforming grids where we assume
4:1 transitions only. See Fig. \ref{fig:non-conforming-interface}.  The steps
for the 3D multi-level blending scheme are analogous and not detailed any
further.  First, we define four matrix operators which allow us to construct
refinement and coarsening procedures within the blending framework. We require
that $N = 2^l,\, l \in \naturals$, and define the following expansion and
compression operators
\begin{equation}
    \label{eq:3D-expansion-mat}
    \mat{E}^{(N \rightarrow 2N)} := \idmat_N \otimes \vec{1}_{2} \;\; \in \reals^{2N \times N}
\end{equation}
and
\begin{equation}
    \label{eq:3D-compression-mat}
    \mat{C}^{(N \rightarrow \tfrac{N}{2})} := \frac{1}{2} \, \idmat_N \otimes \vec{1}_2^T \;\; \in \reals^{\tfrac{N}{2} \times N}.
\end{equation}
\begin{figure}[H]
\centering
\ifnum \compilewithplots = 1
    \includegraphics[width=0.7\textwidth]{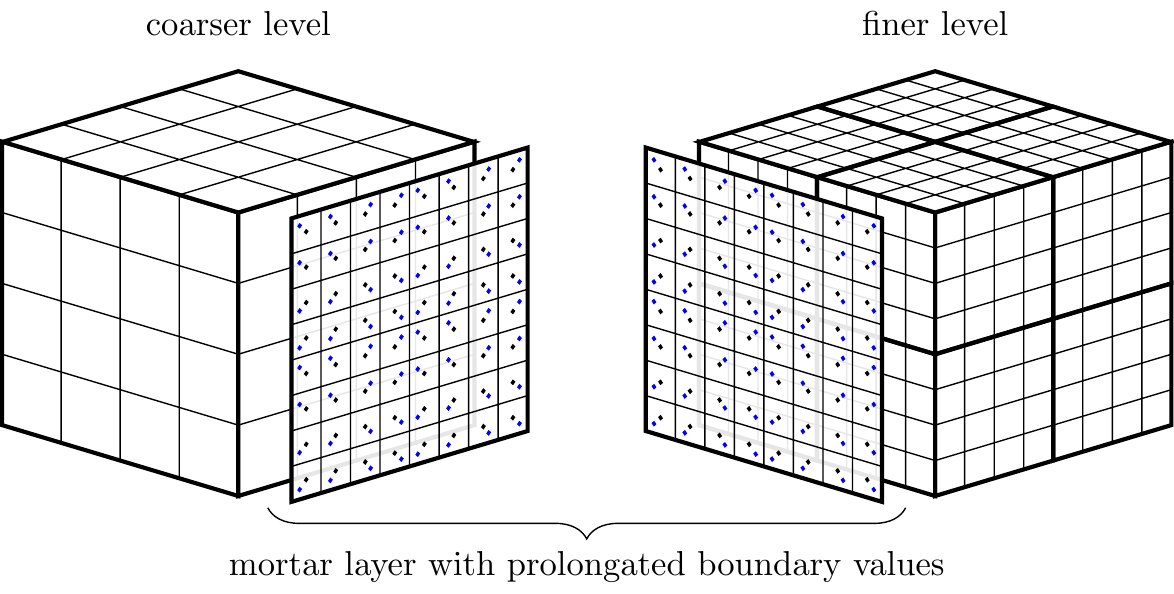}
\fi
\caption{Schematic of a non-conforming interface with 4:1 adjacent cells for $N
= 4$. The two kinds of boundary values occurring in
\eqref{eq:3D-prolongation-nonconforming} are illustrated.  Compare with Fig.
\ref{fig:conforming-interface}.}
\label{fig:non-conforming-interface}
\end{figure}
Moreover, we need the projection operator $\projMat^{(N \rightarrow M)}$
introduced in section \ref{sec:projection} for mapping $N$ node values to $M$ mean
values. We define:
\begin{equation*}
    \underbrace{\mat{E} := \mat{E}^{(N \rightarrow 2N)}, \quad \mat{C} := \mat{C}^{(N \rightarrow \tfrac{N}{2})}}_{\text{mean values to mean values}}
    \quad \text{and} \quad
    \underbrace{\mat{Z} := \projMat^{(N \rightarrow 2N)}, \quad \mat{Q} := \projMat^{(N \rightarrow \tfrac{N}{2})}}_{\text{node values to mean values}}.
\end{equation*}
We refine the parent element $q$ as
\begin{equation}
\label{eq:3D-refinement}
\block{\bar{r}}_q = (1-\alpha_q) \, \block{\bar{u}}_q \, \bigtimes_{d=x}^z \,\mat{E}
    + \alpha_q \, \block{\tilde{u}}_q \,\bigtimes_{d=x}^z \,\mat{Z} \;\; \in \reals^{2N \times 2N \times 2N}.
\end{equation}
The refined block $\block{\bar{r}}_q$ is then split into $8$ child elements. The
reverse operation, namely coarsening, is done by compressing each child
element $\block{\bar{u}}_{r}, \; r = 1,\ldots,8$, separately:
\begin{equation}
\label{eq:3D-coarsening}
\block{\bar{c}}_{q} = (1-\alpha_q) \, \block{\bar{u}}_{q} \, \bigtimes_{d=x}^z \,\mat{C}
    + \alpha_q \, \block{\tilde{u}}_{q} \,\bigtimes_{d=x}^z \,\mat{Q} \;\; \in \reals^{\tfrac{N}{2} \times \tfrac{N}{2} \times \tfrac{N}{2}}.
\end{equation}
Afterwards the family of $8$ blocks $\block{\bar{c}}_{q}$ is glued together.
These blending operations are especially useful for refining or coarsening of
highly oscillatory data or at pronounced discontinuities.\\
\end{changed}

For the treatment of non-conforming 4:1 interfaces we need a procedure which
maps the boundary values from the coarser element to a so-called mortar layer
\cite{Koprivaetal2000} that has the same resolution as the adjacent four
smaller elements. See Fig. \ref{fig:non-conforming-interface}. To compute the
mortar layer of the coarser element $q$, we formulate
\begin{equation}
\label{eq:3D-prolongation-nonconforming}
\left(\face{\bar{u}}_{q\pm\frac{1}{2}}^{\pm d}\right)_{\text{coarse}}
    = (1-\alpha_q)\,\underbrace{\mat{E}\,\big(\block{\bar{u}}_q \times_d \vec{s}^{\pm}\big)\,\mat{E}^T}_{\mathclap{\substack{\text{face centered}\\\text{mean values}}}}
    + \; \alpha_q \,\underbrace{\mat{\tilde{Z}} \overbrace{\big( \block{\tilde{u}}_q \times_d \vec{b}^{\pm} \big)}^{\substack{\text{nodal boundary}\\\text{values}}} \mat{\tilde{Z}}^T}_{\mathclap{\substack{\text{face centered}\\\text{mean values}}}}  \;\: \in \reals^{2N \times 2N}.
\end{equation}
Note the similarity to the prolongation procedure for conforming interfaces
\eqref{eq:3D-prolongation}. The coarse side of the mortar layer
$\left(\face{\bar{u}}_{q\pm\frac{1}{2}}^{\pm d}\right)_{\text{coarse}}$ is
split up to match the faces of the four finer elements $r = 1,\ldots,4$. The
boundary values of the smaller elements are constructed by applying
\eqref{eq:3D-prolongation} individually. Only the mortar layers together with
the associated blending factors $\alpha_r$, $r = 1,\ldots,4$, are supposed to
be exchanged between processes in case of distributed computing.  The separate
computation of the interface fluxes $\face{\bar{f}}^{*}_{\text{fine},r} \in \reals^{N
\times N}$, $r = 1,\ldots,4$, follows the formulas
\eqref{eq:3D-boundary-value-recon}-\eqref{eq:3D-blend-interface-fluxes}.  While
the resulting four fluxes can be directly copied back to the smaller faces
without further treatment, they have to be mapped to the coarser side via
$L^2$-projection
\cite{Koprivaetal2000}
\begin{align*}
&  \int_{-\frac{1}{2}}^{0}\int_{-\frac{1}{2}}^{0}\left(f^{*}_{\text{coarse}} - f^{*}_{\text{fine},1} \right)\phi\,d\xi\,d\eta
+ \int_{0}^{\frac{1}{2}}\int_{-\frac{1}{2}}^{0}\left(f^{*}_{\text{coarse}} - f^{*}_{\text{fine},2} \right)\phi\,d\xi\,d\eta\\[2mm]
 + &\int_{-\frac{1}{2}}^{0}\int_{0}^{\frac{1}{2}}\left(f^{*}_{\text{coarse}} - f^{*}_{\text{fine},3} \right)\phi\,d\xi\,d\eta
+ \int_{0}^{\frac{1}{2}}\int_{0}^{\frac{1}{2}}\left(f^{*}_{\text{coarse}} - f^{*}_{\text{fine},4} \right)\phi\,d\xi\,d\eta = 0.
\end{align*}
The integrals are evaluated in the reference space of the coarse face. Applying
exact quadrature rules gives
\begin{align*}
\left(\face{\tilde{f}}^{*}_{\text{coarse}}\right)_{ij} =
    & \sum_{kl}^N \left(\face{\tilde{f}}^{*}_{\text{fine},1}\right)_{kl}
        \ell_i\left(\frac{1}{2}\,\xi_k  - \frac{1}{4}\right) \frac{\omega_k}{2\,\omega_i}
        \ell_j\left(\frac{1}{2}\,\eta_l - \frac{1}{4}\right) \frac{\omega_l}{2\,\omega_j} + \\[2mm]
    & \sum_{kl}^N \left(\face{\tilde{f}}^{*}_{\text{fine},2}\right)_{kl}
        \ell_i\left(\frac{1}{2}\,\xi_k  + \frac{1}{4}\right) \frac{\omega_k}{2\,\omega_i}
        \ell_j\left(\frac{1}{2}\,\eta_l - \frac{1}{4}\right) \frac{\omega_l}{2\,\omega_j} + \\[2mm]
    & \sum_{kl}^N \left(\face{\tilde{f}}^{*}_{\text{fine},3}\right)_{kl}
        \ell_i\left(\frac{1}{2}\,\xi_k  - \frac{1}{4}\right) \frac{\omega_k}{2\,\omega_i}
        \ell_j\left(\frac{1}{2}\,\eta_l + \frac{1}{4}\right) \frac{\omega_l}{2\,\omega_j} + \\[2mm]
    & \sum_{kl}^N \left(\face{\tilde{f}}^{*}_{\text{fine},4}\right)_{kl}
        \ell_i\left(\frac{1}{2}\,\xi_k  + \frac{1}{4}\right) \frac{\omega_k}{2\,\omega_i}
        \ell_j\left(\frac{1}{2}\,\eta_l + \frac{1}{4}\right) \frac{\omega_l}{2\,\omega_j},
\end{align*}
where $i,j = 1,\ldots,N$ and $\xi_k$, $\eta_l$ are the collocation nodes of the
tensor-product \eqref{eq:tensor-product-ansatz}. We translate above equation
into matrix notation and get
\begin{equation}
\face{\tilde{f}}^{*}_{\text{coarse}} =\;
    \mat{L}_{-}\:\face{\tilde{f}}^{*}_{\text{fine},1}\:\mat{L}_{-}^T+
    \mat{L}_{+}\:\face{\tilde{f}}^{*}_{\text{fine},2}\:\mat{L}_{-}^T+
    \mat{L}_{-}\:\face{\tilde{f}}^{*}_{\text{fine},3}\:\mat{L}_{+}^T+
    \mat{L}_{+}\:\face{\tilde{f}}^{*}_{\text{fine},4}\:\mat{L}_{+}^T \; \in \reals^{N \times N},
\end{equation}
where $\left(\mat{L}_{\pm}\right)_{ij} = \ell_i\left(\frac{1}{2}\,\xi_j \pm
\frac{1}{4}\right) \frac{\omega_j}{2\,\omega_i}$, $i,j = 1,\ldots,N$.
Analogous to \eqref{eq:3D-blend-interface-fluxes} we determine a common surface
flux for the coarse side
\begin{equation}
\label{eq:3D-blend-interface-fluxes-non-conforming}
    \bar{\face{f}}_{q-\frac{1}{2}}^{*} = 
        (1-\alpha_{q-\frac{1}{2}}) \;\; \mat{C}\;\bar{\face{f}}_{q-\frac{1}{2}}^{*'} \; \mat{C}^T
        \;\;+ \;\;
        \alpha_{q-\frac{1}{2}} \;\; \projMat \; \left(\face{\tilde{f}}_{\text{coarse}}^{*}\right)_{q-\frac{1}{2}} \; \projMat^T \;\:\;\; \in \reals^{N \times N},
\end{equation}
where $\alpha_{q-\frac{1}{2}} = \frac{1}{2}\big(\alpha_{q-1} +
\min_{r=1,\ldots,4} \; (\alpha_r)_q\big)$. The face centered fluxes
$\bar{\face{f}}_{q-\frac{1}{2}}^{*'} \in \reals^{2 N \times 2 N}$ are the
result of the appropriate glue operation of the four first order fluxes
$\left(\bar{\face{f}}_r^{*'}\right)_{q-\frac{1}{2}} \in \reals^{N \times N}$,
$r = 1,\ldots 4$. This concludes the description of the treatment of
non-conforming 3D Cartesian grids within the blending framework.

\subsection{Sketch of the Algorithm for the multi-level blending scheme}
\label{sec:appendix-sketch-of-alogrithm}
\begin{changed}
The algorithm of the multi-level blending scheme is more involved but still
follows the same general sequence of steps as in the single-level scheme
outlined in section \ref{sec:sketch-of-alogrithm}. Moreover, we use the notation
for the 3D scheme (appendix \ref{sec:appendix-3D-extension}). At each
Runge-Kutta stage we do:
\begin{enumerate}[label=(\Roman*),start=1,itemsep=2mm,leftmargin=1cm]
\item \label{item:algo-multi-level} Initialize two tracing variables: $n_{q,\text{min}} = n_{q,\text{max}} = 8$.
\item Loop from highest to lowest order: $n_q = 8,4,2$.
    \begin{itemize}[label=$\triangleright$,topsep=2mm,itemsep=2mm,leftmargin=5mm]
        \item Set $n_{q,\text{min}} = n_q$ of element $q$.

        \item For each sub-element $s$ in element $q$ at level $n_q$:
            \begin{itemize}[topsep=2mm,itemsep=2mm,leftmargin=5mm]
                \item Reconstruct the polynomial
                    $\left(\block{\tilde{u}}_s\right)_q^{\bigO(n)}$ from given mean values
                    $\block{\bar{u}}_q$ via \eqref{eq:reco-2nd-order}, \eqref{eq:reco-4th-order} or
                    \eqref{eq:reco-8th-order}.

                \item If the reconstructed polynomial $\left(\block{\tilde{u}}_s\right)_q^{\bigO(n)}$ contains
                non-permissible states, see \eqref{eq:positivity-limiter-2}, then
                calculate the limited version $\left(\left(\block{\tilde{u}}_s\right)_q^{\bigO(n)}\right)^{(\beta)}$ as in
                \eqref{eq:positivity-limiter}.

                \item If the squeezing parameter $(\beta_s)_q^{\bigO(n)}$ is below $\beta_L$ then set $(\alpha_s)_q^{\bigO(n)} := 0$ else
                compute the blending factor $(\alpha_s)_q^{\bigO(n)}$ via \eqref{eq:calc-blending-factor} from the \textbf{unlimited} polynomial
                $\left(\block{\tilde{u}}_s\right)_q^{\bigO(n)}$ .
            \end{itemize}

            \item If all blending factors $(\alpha_s)_q^{\bigO(n)}$ are above $\alpha_H$ then break the loop.

            \item If all blending factors $(\alpha_s)_q^{\bigO(n)}$ are below $\alpha_L$ then set $n_{q,\text{max}} = n-1$.

    \end{itemize}

\item Compute the multi-level blended boundary values
$\bar{\face{u}}_{q\pm\frac{1}{2}}^{\pm d}$ analog to \eqref{eq:3D-prolongation}
within the minimum and maximum bounds given by $n_{q,\text{min}}$ and
$n_{q,\text{max}}$. Exchange them together with all involved blending factors
$\alpha_q^{\bigO(n)}$ alongside zone boundaries in case
of distributed computing.

\item Determine the multi-level blended common surface flux
$\bar{\face{f}}_{q-\frac{1}{2}}^{*}$ analog to
\eqref{eq:3D-blend-interface-fluxes} within the minimum and maximum bounds
given by $n_{q,\text{min}}$ and $n_{q,\text{max}}$.

\item \label{item:multi-level-scheme} Compute the multi-level blended
right-hand-side $\partial_t \block{\bar{u}}_q$ as described in section
\ref{sec:multilevel-blending} but within the minimum and maximum bounds 
given by $n_{q,\text{min}}$ and $n_{q,\text{max}}$.

\item \label{item:time-int-multi-level} Forward in time to the next Runge-Kutta stage
and return to step \ref{item:algo-multi-level}.
\end{enumerate}
The switching thresholds are set to $\alpha_H := 0.99$ and $\alpha_L := 0.01$
and the limiter threshold to $\beta_L := 0.95$. Note that the algorithm only
applies the blending procedure where necessary in order to maintain the overall
performance of the scheme. The two bounds $n_{q,\text{min}}$ and
$n_{\text{max}}$ even avoid redundant computation at levels sorted out by the
shock indicator.
\end{changed}

\subsection{N-mode product}
\label{sec:appendix-n-node-product}
Given tensor $\block{u} \in \reals^{I_1 \times I_2 \times \ldots \times I_D}$
with a matrix $\mat{A} \in \reals^{J \times I_d}$ than the \emph{n-mode}
product is defined as
\begin{equation}
    \left(\block{u} \times_d \mat{A}\right)_{i_1 \ldots i_{d-1} \; j \; i_{d+1}\ldots i_D} = \sum_{i_d}^{I_d} u_{i_1\ldots i_d \ldots i_D} \, A_{j i_d}.
\end{equation}
The resulting tensor has following dimensions
\begin{equation*}
    \block{u} \times_d \mat{A} \in \reals^{I_1 \times \ldots \times I_{d-1} \times J \times I_{d+1} \times \ldots \times I_D}.
\end{equation*}
Given the same tensor $\block{u}$ and a vector $\mat{b} \in \reals^{I_d}$ than
the \emph{n-mode} product acts as a contraction of $\block{u}$ with $\vec{b}$
along dimension $d$. We write
\begin{equation}
    \left(\block{u} \times_d \vec{b}\right)_{i_1 \ldots i_{d-1} \; i_{d+1}\ldots i_D} = \sum_{i_d}^{I_d} u_{i_1\ldots i_d \ldots i_D} \, b_{i_d}.
\end{equation}
The resulting tensor has following dimensions
\begin{equation*}
    \block{u} \times_d \vec{b} \in \reals^{I_1 \times \ldots \times I_{d-1} \times I_{d+1} \times \ldots \times I_D}.
\end{equation*}

\section*{Conflict of Interest}
On behalf of all authors, the corresponding author states that there is no
conflict of interest. 

\bibliographystyle{plain}
\bibliography{citations}

\end{document}